\documentclass[bj,preprint,12pt]{imsart}

\RequirePackage[utf8]{inputenc}
\RequirePackage{natbib} %possible option: [numbers]
\RequirePackage[colorlinks,citecolor=blue,urlcolor=blue]{hyperref}
\RequirePackage[normalem]{ulem} %  \usepackage{pdftricks}

\def\journal@name{}
\setattribute{copyright}  {text} {}
\usepackage{pagecolor}

\usepackage{datetime}
\usepackage{todonotes}	
\usepackage{booktabs}
\usepackage[inline]{enumitem}	
\usepackage{fullpage}	
  \usepackage{amsthm,amsmath}
  \usepackage{amsfonts} 
  \usepackage{empheq}
  \usepackage{mathrsfs}                  %paquet pour les formules maths
\usepackage{amssymb} 
\usepackage{MnSymbol}%needs to be after amssymb
  \usepackage{dsfont}
  \usepackage{cases}
  \usepackage{booktabs}
  \usepackage{float}

  \theoremstyle{plain}
\startlocaldefs

 \newcommand\permanentcomments[2][Notes]{\rule{\linewidth}{0.4pt}\\\noindent\begin{pdfsidelinecomment}[color=blue,linewidth=2bp,linesep=.3cm]{#1}{
{\footnotesize
 #2
}

\noindent\makebox[\linewidth]{\rule{\linewidth}{0.4pt}}
}\end{pdfsidelinecomment}\\\vspace{-.1cm}}

\newcommand\todomargin[1]{{\todo{\footnotesize #1}}}%
\newcommand\comments[1]{}
\newcommand\todof[1]{}%
\newcommand\delete[1]{{\color{blue}\sout{ #1}}}%
\renewcommand\delete[1]{}%

\def\size{\mathrm{cardinality}}

\def\translated{\hookrightarrow}
\newcommand{\Extension}[1]{\mathbf{\mathring{\text{$#1$}}}}
\newcommand{\Comp}[1]{\mathbf{\bar{\text{$#1$}}}}

\def\isogeq{\geqclosed}
\def\isog{\gtrclosed}
\def\isoeq{\closedequal}

\def\restrict{\Big|}
	\newtheorem{example}{Example}
    \newtheorem{examples}[example]{Examples}
    \newtheorem{theorem}{Theorem}

	\newtheorem{definition}{Definition}
	\newtheorem{definitions}[definition]{Definitions}
	\newtheorem{propertydefinition}[definition]{Property-Definition}
	
	\newtheorem{property}[theorem]{Property}

	\newtheorem{remark}{Remark}

    \DeclareSymbolFont{lettersA}{U}{txmia}{m}{it} %symbol for pi, the real number
 	\DeclareMathSymbol{\piup}{\mathord}{lettersA}{"19}

    \def\delindex#1#2{{}}

    \newcommand{\rangeformat}[1]{\mathscr{#1}}
    \newcommand{\range}[1]{\rangeformat{#1}}
    \newcommand{\realisationformat}[1]{#1}
    \newcommand{\deterministicfunctionformat}[1]{\mathbf{#1}}
    
	\newcommand{\allfunc}[2]{#1\to#2}
	
	\def\class{\mathrm{class}}
	
	\def\setsep{:}
	\def\pperp{\perp \!\!\! \perp}
	
	\def\varindep{\perp}
	\def\indep{\pperp}
	
	\def\dependsonly{\prec}
	\def\dominatedby{\ll}
	\def\interclass{\sqcap}
	
    \def\somesubset{A}
    \def\Atrandomize{\mathrm{Atrandomize}}
    \def\somesubsetB{B}
    \def\Borel{\mathrm{Borel}}
	\def\codomain{\mathrm{codomain}}
	\def\Design{\Delta}
	\def\design{\realisationformat{\delta}}
    \def\Dirac{\mathrm{Dirac}}
    \def\setindicator{\mathds{1}}
	\def\domain{\mathrm{domain}}
    \def\derive{\mathrm{d}} %Density operator
	\def\spaceE{\range{E}}
	\def\expected{E}

	\def\density{\mathrm{f}} %Density operator
	\def\spaceF{\range{F}}
	\def\spaceG{\range{G}}
    \def\funcc{g}
    \def\func{h}
    \def\Sampleind{I}

    \def\image{\mathrm{image}}
    \def\individual{k}
	\def\Sampleindexset{L}
	\def\sampleindex{\ell}
	\def\Popsize{N}
	\def\popsize{N^\star}
	\def\Samplesize{n}
	\def\samplesize{n^\star}
	
    \def\priordist{Q}%mind\priordistset=capital\priordist=B
    \def\priordistset{\range{\priordist}}%mind\priordistset=capital\priordist=B
	\def\Samplepopmap{R}
	\def\samplepopmap{r}
	
	\def\Sufstat{S}
	\def\sufstatf{\deterministicfunctionformat{s}}
	\def\sufstat{\realisationformat{s}}
	\def\Sampleset{S}
	
	\def\Transfo{T}
	\def\transfo{\realisationformat{t}}
	\def\pop{\realisationformat{u}}
	
	\def\Pop{U}
	\def\powerset{\mathrm{powerset}}
	\def\superpop{\range{U}}
	\def\Rv{V}
	\def\Rvv{W}
	\def\rv{\realisationformat{v}}
	\def\rvv{\realisationformat{w}}
	\def\elt{x}
	\def\Obs{X}
	\def\obs{\realisationformat{x}}
    \def\obsf{\deterministicfunctionformat{x}}
	\def\eltt{y}
	\def\Signal{Y}
	\def\signal{\realisationformat{y}}
	\def\Designvar{Z}
	\def\likelihood{\mathscr{L}}
	\def\designvar{\realisationformat{z}}

	\def\Id{\operatorname{Id}}
	
    \def\dominant{\zeta}

	\def\biclass{\class}
	
	\def\ignoreset{\Phi}
	
    \def\tribu{\mathfrak{S}}
    
    \def\param{{\theta}}
    \def\paramf{{\boldsymbol\param}}
 	\def\xiparam{{\Comp\param}}
 	
 	\def\Thetaset{\Theta}
 	\def\Xiset{\Comp\Theta}
    
    \def\inclusionprob{\pi}%\mathbbold{\Pi}
    \def\inclusionexpec{\upsilon}%\mathbbold{\Pi}
\def\Inclusionexpec{\Upsilon}%\mathbbold{\Pi}

    \def\Inclusionprob{\Pi}%\mathbbold{\Pi}
    \def\countingmeasure{\mu}

    \DeclareMathOperator*{\argmax}{arg\,max}
	\DeclareMathOperator*{\argmin}{arg\,min}
		
    \def\d{\operatorname{d}\!} %pour differentiate operator in integrals

    \makeatletter
    
\def\env@sqcases{%
  \let\@ifnextchar\new@ifnextchar
  \left|
  \def\arraystretch{1.2}%
  \array{@{}l@{\quad}l@{}}%
}
\makeatother

\usetikzlibrary{backgrounds}
\makeatletter

\tikzset{%
  fancy quotes/.style={
    every path/.style={double=black,color=black},
    text width=\fq@width pt,
    align=justify,
    inner sep=1em,
    anchor=north west,
    text=black,
    color=black,
    minimum width=\linewidth,
  },
  fancy quotes width/.initial={.8\linewidth},
  fancy quotes marks/.style={
    scale=4,
    text=black,
    inner sep=0pt,
  },
  fancy quotes opening/.style={
    fancy quotes marks,
  },
  fancy quotes closing/.style={
    fancy quotes marks,
  },
  fancy quotes background/.style={
    show background rectangle,%
    inner frame xsep=0pt,
    text=black,
    background rectangle/.style={
      fill=gray!25,
      rounded corners,
    },
  }
}

\newcommand{\fancyquote}[2]{%
\noindent
\tikzpicture[fancy quotes background];
\node[fancy quotes opening,anchor=north west] (fq@ul) at (0,0) {``};
\tikz@scan@one@point\pgfutil@firstofone(fq@ul.east)
\pgfmathsetmacro{\fq@width}{\linewidth - 2*\pgf@x}
\node[fancy quotes] (fq@txt) at (fq@ul.north west) {#2};
\node[overlay,fancy quotes closing,anchor=east] at (fq@txt.south east) {''};
\node[fancy quotes] (fq@author) at (fq@txt.south west)   {#1};
\endtikzpicture}

\makeatother

\endlocaldefs

\usepackage[subject={Top1},author={Josef}]{pdfcomment}

\begin{document}
\begin{frontmatter}
\title{On the definition of an informative vs. an ignorable nuisance process.}

\runtitle{On the definition of informative selection}
\thankstext{T1}{}
\thankstext{T2}{}
\begin{aug}
\author{\fnms{Daniel}
				\snm{Bonn\'ery}
				\thanksref{ad1,e1}%<-do not suppress this percent
				\ead[label=e1,mark]{dbonnery@umd.edu}% <-do not suppress this percent
				\ead[label=u1,url]{jpsm.umd.edu/facultyprofile/Bonn\'ery/Daniel}},
\author{\fnms{Joseph}
				\snm{Sedransk}
				\thanksref{ad2,e2}%<-do not suppress this percent
				\ead[label=e2,mark]{jxs123@case.edu}% <-do not suppress this percent
				\ead[label=u2,url]{jpsm.umd.edu/facultyprofile/Sedransk/Joseph}},

				%\thankstext{t1}{Some comment}
%\thankstext{t2}{First supporter of the project}
%\thankstext{t3}{Second supporter of the project}
\runauthor{D. Bonn\'ery, J. Sedransk}
\affiliation{University of Maryland\thanksmark{ad1}}

\address[ad1]{\printead{e1}, \printead{u1}}
\address[ad2]{\printead{e2}, \printead{u2}}
%\printead{add} 
\end{aug}

	\begin{abstract}

This paper is an early version.

We propose to generalise the notion of "ignoring" a random process as well as the notions of informative and ignorable random processes in a very general setup and for different types of inference (Bayesian or frequentist), and for different purposes (estimation, prediction or testing). We then confront the definitions we propose to mentions or definitions of informative and ignorable processes found in the litterature. To that purpose, we provide a very general statistical framework for survey sampling in order to define precisely the notions of design and selection, and to serve to illustrate and discuss the notions proposed.

	\end{abstract}



\end{frontmatter}

{\huge This is the version of \shortdate\today}

\tableofcontents

%\newpage
%\comments{Temporary comments in red}
%
%\permanentcomments{
%Comments to be kept for us in some version but deleted in the version submitted to referees in yellow
%
%}

%\todof{Things to do in green:

% 1. Rewrite section 1
% 
% 2.  Write notations index
%}

\newpage

\section{Introduction and motivation}

It is common in survey sampling to decompose the random process behind the observations into distinct random processes: the variable of interest generation, the design variable generation, the sample selection according to the design, the non response, the measurement error, ...etc. This has been explained by \cite{Pfeffermann1998a} as well as \cite{Skinner1994}. \cite{Rubin1976} explains what it means to ignore a particular random process (the one that causes missing data) for a particular type of inference (likelihood based inference). \cite{SugdenSmith1984} oppose the notions of informative selection \cite{Scott1975} and the notion of ignorable missing data mechanism. We have not found a unique definition of informative process, but rather different ad-hoc definitions, that can be applied in different statistical frameworks. 
\cite{Scott1977}, commenting on \cite{Godambe1966a}, suggests that in a Bayesian framework that a selection is "informative" if the posterior distribution of some estimates depends only on the sample drawn and not on the design used to use it.
\cite{Rubin}  defines what ignoring a missing data mechanism means and gives sufficient conditions under which the likelihood based inference will not be affected after ignoring the missing data mechanism. \cite{Pfeffermann1988a} gives a heuristic definition of informative selection as a process that has to be taken into account in the inference.\cite{Pfeffermann1998a} also defines the sample distribution by opposition to the population distribution, and focuses on one of the possible effect of informative selection: under informative selection, the sample distribution is different from the population distribution.
There exists other definitions in the litterature that focus on the possible causes of informative selection: \cite{CasselSarndalWretman1977} explains that selection is informative when the design variable or the design is dependent on the study variable. Other authors have a very restrictive and convoluted definition, as \cite{Fuller2009} for whom informative selection occurs when the distribution of the inclusion probabilities conditionnally to the design variable and the study variable depends on the study variable.
There are issues with the two last definitions:  \cite{CasselSarndalWretman1977}  provides a definition in a fixed population model for design-based inference, under which design and study variables are not random, so are stochastically independent, so the definition has no interest. Even when transposed to a model based framework, this definition is not satisfying.
Indeed, when sampling with replacement, even independent on the study variable, and the population index on the sample is not observed, the early papers by \cite{Scott1977}, or \cite{Rubin1976} do not consider this case, but if we follow their reasoning, sampling with replacement should be considered as informative. The definition of \cite{Fuller2009} is too restrictive, and wrong: for example, in cluster sampling where all clusters have the same probability to be selected, the inclusion probabilities are constant, but the selection is informative as it induces dependence among the observations on the sample (see \cite{Bonnerytheseen}). 
So there is a need to clarify what informative selection means. One can size this opportunity to define, in general, what an informative vs ignorable process is. This could be applied to 
a coarsening  process for example (see \citet{Heitjan1991}). The task of defining what an ignorable vs. informative process is reminds the attempts to define ancillarity and sufficiency in the presence of a nuisance parameter as the framework is the same: we are in presence of two random processes, one is the process of interest, and the second a nuisance process. Any distribution in a model can be viewed as the combinaision of a marginal distribution of the process of interest and of the distribution of the nuisance process conditionnally on the process of interest. This framework is used in \cite{Rubin1976}.
Nevertheless, it is not straighforward to propose a definition as it has be applicable to different statistical frameworks, and has to take into account the type of inference(Bayesian, frequentist), the goal of the inference (testing, estimation, prediction, model selection), the nature of the random process with respect to the observations (is the ignorable random process observed or latent), the criteria to define equivalent models (models under which estimators have the same bias, the same distribution, ...), as well as what we will call the injectivity of the model separation into the marginal model for the process of interest and the conditional model on the nuisance process. The two last aspects require to address issues that were not raised in the papers by \cite{Scott1975} and \cite{Rubin1976}: both of them consider the sample index as observed, but in practice (for example for web surveys) it may happen that one cannot identify the units for each observation of the sample, and one may not be able to identify duplicates for example: the sample index is considered as latent, and \cite{Rubin1976} lists as a necessary condition for ignorability  the separability of the model into the process of interest marginal model by the nuisance process conditional model.

We can follow the early papers by \cite{Godambe1966a},\cite{Scott1975} and \cite{Rubin1976}, to setup a definition of informative or ignorable process. These are the steps to follow:
\begin{enumerate*}
\item Provide a general framework and formalise the condition that the observations are the outcome of two random processes,
\item Consistently with \cite{Rubin1976}, generalise the notion of ignoring one process introduced by to this general framework, for Bayesian or frequentist inference.
\item Consistently with \cite{Rubin1976}, explain what equivalent inferences from different models mean, especially detail what it means when one or the two processes are latent processes.
\item Define an ignorable process (and by opposition an informative process) as any process whose ignorance will lead to an equivalent inference.
\end{enumerate*}
After we get a definition for informative vs. ignorable process, we will critically examine the existing definitions found in the litterature.

In section 2, we give some mathematical background for the different notions that will be defined. In section 3, we propose a general framework for survey sampling. In section 4, we propose to generalise the notion of ignoring a random process. In section 5, we give the definition of ignorable versus informative process, and we characterise such processes. Section 6 is a discussion on different topics, including a critical examination of the different definitions of informative selecction that can be found in the scientific litterature,  we also discuss the link between the notions of informativity and information, sufficiency and ancillarity in presence of nuisance parameters, as well as the debate on the likelihood principle. In appendix, we provide notes on a selection of papers that are link more or less closely to the current topic.

\section{Mathematical notations and tools}

For clarity, it is imperative to remove any notation or concepts ambiguity. Especially with respect with the term dependence when applied to deterministic functions. 
Quoting \cite{NeymanPearson1936}, ``it is inevitable, [...] that a paper dealing with [the] problem [of what conclusions regarding sufficient statistics may be drawn from the existence of uniformly most powerful tests, or vice versa] should bear some mark of the theory of functions, in spite of its concern with statistical questions''. This effort was visible in \cite{BarndorffNielsen1973} who cautiously defined ``variation independence'', in \cite{Rubin1976} who cautiously defined ``distinct parameters'', or in \cite{Basu1977} who reproduced definitions usefull for his paper. We have listed all the concepts required to define an informative latent or observed process in this section.

We tried to follow this rule for notations:
non random elements are denoted with lowercase roman letters, random variables with uppercase roman letters, domains and codomains with calligraphic uppercase roman letters, functions on the model distributions will be bold greek letters, and the values they take  greek non bold letters.

\subsection{Set theory: notation conventions and definitions}
The following conventions on functions will be used in this document. Notations 1. to 6. below can also be found in \citet[(Section 0.1. Basics)]{DummitFoote2003}.
\begin{enumerate*}[mode=unboxed]
\item Given two sets $\spaceE$ and $\spaceF$, $\allfunc{\spaceE}{\spaceF}$ is the set of all functions from $\spaceE$ to $\spaceF$.
\item The statement ``let $\func:\allfunc{\spaceE}{\spaceF}$'' means ``let $\spaceE$ and $\spaceF$ be two sets'' if they are not already defined  and ``let $\func\in \allfunc{\spaceE}{\spaceF}$''.
\item The domain of $\func$ will be denoted $\domain(\func)$, its codomain: $\codomain(\func)$, and its image or range $\image(\func)$.
\item The statement ``$\func:\allfunc{\spaceE}{}$'' means ``let $\func$ be a function defined on $\spaceE$'' (the codomain is not always specified).
\item The statement ``$\func:{\mathbb{N}}\to{\mathbb{N}},:\elt\mapsto \elt+1$'' means let $\func$ be the function mapping $\mathbb{N}$ to $\mathbb{N}$ defined by the formal expression $\func(\elt )=\elt +1$. 
\item Let $\func:\spaceE\to \spaceF$, $\func':\spaceF\to \spaceG$, then $\func'(\func)$ is the function $\func'\circ\func:\spaceE\to\spaceG,\elt \mapsto \func'(\func(\elt ))$.
\item Let $\func:\spaceE\to \spaceF$, $\func':\spaceE\to (\allfunc{\spaceF}{\spaceG})$, $v\in \spaceF$, then $\func'[\func]$ is the function $\func'[\func]:\spaceE\to\spaceG,\elt \mapsto (\func'(\elt ))(\func(\elt ))$, and $\func'[v]$ is the function $\func'[v]:\spaceE\to\spaceG,\elt \mapsto (\func'(\elt ))(v)$.
\item Let $\func:\spaceE\to \spaceF$, $\spaceE'\subseteq \spaceE$. Then $\func\restrict_{\spaceE'}$ denotes the restriction of $\func$ to $\spaceE'$ defined as $\func\restrict_{\spaceE'}:\spaceE'\to\spaceF,\elt \mapsto\func(\elt )$. 
\item Let $\func:\spaceE\to \{\text{subsets of }\spaceF\}$, $\func':\spaceE\to (\allfunc{\spaceF}{\spaceG})$,  then \begin{equation*}\func'\restrict_{[\func]}:\begin{array}[t]{lcl}
\spaceE&\to&\bigcup_{\spaceF'\subseteq\spaceF}(\spaceF'\to\spaceG),\\\elt &\mapsto& \func'\restrict_{[\func]}(\elt ):\begin{array}[t]{lcl}
\func(\elt )&\to&\spaceG,\\\eltt&\mapsto& (\func'(\elt ))(\eltt).\end{array}\end{array}\end{equation*}
\item 
 For $\func:\allfunc{\spaceE}{}$, let $\sim_\func$ be the equivalence relation on 
 $\spaceE$ defined by $\elt\sim_\func \elt'\leftrightarrow \func(\elt)=\func(\elt')$.
 The set of classes of $\spaceE$ for $\sim_\func$ is denoted $\spaceE/\func$, the class of $\elt$ for $\sim_\func$ is $\class_\func(\elt)=\func^{-1}\left(\{\func(\elt)\}\right)$.\item The identity on a set $\spaceE$ is denoted $\Id_\spaceE$, $\mathds{1}_\spaceE=:\spaceE\to{1},\elt\mapsto 1$, and for $\somesubset\in\spaceE$, $\setindicator_\somesubset:\somesubset\to\{0,1\},\elt\mapsto 1$ if $\elt\in\somesubset$, $0$ otherwise.
\item Given $\func,~\func':\spaceE\to$, then
$\func$ and $\func'$ are variation independent and one notes
$\func\varindep\func'$ if and only $\image(\func,\func')=\image(\func)\times\image(\func')$,
\item $\func:\spaceE\to$ (deterministically) depends only on $\func':\spaceE\to$ and one notes $\func\dependsonly\func'$ if and only if $\exists\funcc:\to$ such that $\func=\funcc\circ\func'$.
\end{enumerate*}

\begin{propertydefinition}[Complement and distinct complement of a function]
 Let $\func: \allfunc{\spaceE}{}$, then a complement of $\func$ is any function $\Comp\func:\allfunc{\spaceE}{}$ such that  
 $(\class_{\func},\class_{\Comp\func}):\spaceE\to (\spaceE/\func)\times(\spaceE/\Comp\func), \elt\mapsto (\class_\func(\elt),\class_{\Comp\func}(\elt)))$ is injective.
 If a function $\Comp\func$ satisfies this property, it will be called a complement of $\func$, and we will write $(\func,\Comp\func)\isogeq\spaceE$. 
 If $\biclass_{\func,\Comp\func}$ is bijective (e.g. if in addition $\func\varindep\Comp\func$), then $\Comp\func$ is called a distinct or variation independent complement, and we will write $(\func,\Comp\func)\isoeq\spaceE$, if not, we will write $(\func,\Comp\func)\isog\spaceE$.
 For a complement  $\Comp\func$ of $\func$, define the function $\interclass_{\func,\Comp\func}:\{\func(\elt),\Comp\func(\elt)\setsep  \elt\in \spaceE\} \to \spaceE$ as the inverse of $(\class_{\func},\class_{\Comp\func})\restrict^{\image((\class_{\func},\class_{\Comp\func}))}$. There always exists a complement (for example, $\Id_\spaceE$ is a complement of $\func$), and a sufficient condition for a distinct complement of a function $\func:\spaceE\to$ to exist, is that all the classes $\class_\func(\elt),\elt\in\spaceE$  are in bijection one with another.
\end{propertydefinition}

\begin{remark}[Link with distinct parameters, variation independence, and variation independent complement.]
The definition of variation independence of parameters is given by \citet[p.~2.1.2.]{BarndorffNielsen1973}, \citet[Definition 3]{Rubin1976} uses the term distinct parameters, whereas \citet[p.~357]{Basu1977} uses the term variation independent complement of a parameter.
% A set  indexes another set if there exists a bijection  between the two. 
% If two functions $\paramf,~ \paramf':\range{P}\to$ are such that $(\paramf,\paramf'):P\mapsto(\paramf(P),\paramf'(P))$ is injective, then $(\paramf,\paramf')\restrict^{\image(\paramf,\paramf')}$ is a bijection between $\range{P}$ and $\image(\paramf,\paramf')\subseteq\image(\paramf)\times \image(\paramf')$, so $\image(\paramf,\paramf')$ indexes $\range{P}$. \citet[]{BarndorffNielsen1973} defines two parameters $\param=\paramf(P)$, $\param'=\paramf'(P)$ as variation independent when $(\paramf,\paramf')$ are distinct complements. \citet[Definition 3]{Rubin1976} defines two parameters $\param=\paramf(P)$, $\param'=\paramf'(P)$ as distinct when $\image(\paramf,\paramf')=\image(\paramf)\times \image(\paramf')$.% The injectivity is implicit The fact that they parametrize $\range{P}$ implicitely implies that $(\paramf,\paramf')$ is injective. 
% %So $\param,\param'$ being distinct parameters corresponds to the case where $(\paramf,\paramf')$ are distinct complements. 
% Following the definition of variation independence of \citet[p.~2.1.2.]{BarndorffNielsen1973}, \citet[p.~357]{Basu1977} uses the term variation independent complement [of $\param$] for $\param'$ when $(\paramf,\paramf')\isoeq\range{P}$.
\end{remark}

\begin{examples}
 \begin{enumerate*}
  \item Let $\spaceE=\{(0,0),(0,1),(1,0)\}$, and let $\func:(\elt_1,\elt_2)\mapsto \elt_1$. Then $\func\isog\spaceE$ : the bijection is only possible if the cardinal of $\spaceE$ is the product of the cardinals of 
  $\spaceE/\func$ and $\spaceE/\Comp\func$, and 3 is not a multiple of 2. Consider the complement $\func:(\elt_1,\elt_2)\mapsto \elt_2$ , then $\Comp\func(\elt)=1\Rightarrow \func(\elt)=0)$, whereas $\Comp\func(\elt)=0\nRightarrow \func(\elt)=0$ and $\Comp\func(\elt)=0\nRightarrow \func(\elt)=1$. When $\func\isog\spaceE$, the knowledge of $\Comp\func(\elt)$ may or may not restrict the set of possible values for $\func(\elt)$. More generally, $\func\isoeq\spaceE$ if all the elements of $\spaceE/\func$ are in bijection with each other.
  \item Let $\spaceE=\{(0,0),(0,1),(1,0),(1,1)\}$, and let $\func$, $\Comp\func$ be formally defined as above. Then $(\func,\Comp\func)\isoeq\spaceE$: in this case, the knowledge of 
  $\Comp\func(\elt)$ does not restrict the set of possible values for $\func(\elt)$. 
  \item  Consider the familiy of prabability distributions $\mathscr{P}=\left(\mathrm{Normal}(\mu,\sigma^2)\right)_{(\mu,\sigma^2)\in\mathbb{R}\times\mathbb{R}^+}$, then $(\func,\Comp\func)\isoeq\range{P}$ for $\func:P\mapsto\int x\derive P(x)$, $\Comp{\func}:P\mapsto\int (x-\func(P))^2~\derive P(x)$.
  \item  Consider a statistical model $(\Omega,\tribu_\Omega,\range{P})$, a random variable $\Obs:(\Omega,\tribu_\Omega)\to$, $\paramf:\left(P\mapsto P^{\Id_\Omega\mid \Obs}\right)$ and $\Comp\paramf:\left(P\mapsto P^{\Obs}\right)$ then $(\paramf,\Comp\paramf)\isogeq\range{P}$. If furthermore $\forall P,P'\in\range{P}$, $(P'':\tribu_\Omega\to[0,1],A\mapsto \int P'(A\mid \Obs=x)\d P^\Obs(x))\in\mathscr{P}$, then $(\paramf,\Comp\paramf)\isoeq\range{P}$.
 \end{enumerate*}

\end{examples}
% 
% \begin{definition}[Deterministic independence to $\func$ given $\Comp\func$]
% Given two functions $\func:\spaceE\to$, $\func':\spaceE\to $, one says that $\func$ is deterministically independent of $\func'$ and one notes
% $\func\varindep\func'$ if and only if $\exists \Comp\func'$ a distinct complement of $\func'$ such that
% \begin{equation}
% \forall u,u'\in\spaceE,\ \left[\Comp\func'(\elt)=\Comp\func'(u')\Rightarrow\func(\elt)=\func(u')\right] \label{eq:deterninisticdependence}
% \end{equation}
% \end{definition}

% \begin{examples}
% \item $\func:\allfunc{\mathbb{R}^3}{},(\elt_1,\elt_2,\elt_3)\mapsto \elt_1$ is deterministically independent of $\func':\allfunc{\mathbb{R}^3}{},(\elt_1,\elt_2,\elt_3)\mapsto \elt_2$ (consider the distinct complement of $\func'$: $\Comp\func':(\elt_1,\elt_2,\elt_3)\mapsto (\elt_1,\elt_3)$).
% \end{examples}

\subsection{Probability and Statistics Theory}
%\subsubsection{Notations conventions for general concepts}

The following concepts from Probability theory will be used in the paper.
\begin{enumerate*}[mode=unboxed]
\item Given a space $\Omega$, $\tribu_{\Omega}$ always denotes a sigmafield on $\Omega$.
\item The set of all measurable functions from $(\spaceE,\tribu_{\spaceE})$ to $(\spaceF,\tribu_{\spaceF})$ is denoted $\allfunc{(\spaceE,\tribu_{\spaceE})}{(\spaceF,\tribu_{\spaceF})}$.
\item The statement ``let $\Obs:\allfunc{(\spaceE,\tribu_{\spaceE})}{(\spaceF,\tribu_{\spaceF})}$'' means ``let $(\spaceE,\tribu_{\spaceE})$ and $(\spaceF,\tribu_{\spaceF})$ be two measurable spaces '' if they are not already defined  and ``let $\Obs$ be a measurable function from the measurable space $(\spaceE,\tribu_{\spaceE})$ to the measurable space $(\spaceF,\tribu_{\spaceF})$''. $(\spaceF,\tribu_{\spaceF})$ may be omitted.
\item The sub-sigmafield of $\tribu_{\Omega}$ generated by a random variable $\Obs:\allfunc{(\Omega,\tribu_{\Omega})}{(\spaceF,\tribu_{\spaceF})}$ is denoted $\tribu(\Obs)$ and is defined by $\tribu(\Obs)=\{\Obs^{-1}(A)\setsep  A\in\tribu_\range{\Obs}\}$. 
\item Given a measure $P$ on $(\Omega, \tribu_{\Omega})$ and $\Obs:\allfunc{(\Omega,\tribu_{\Omega})}{(\range{\Obs},\tribu_{\range{\Obs}})}$, $P^\Obs$ is the measure on $(\range{\Obs},\tribu_{\range{\Obs}})$ defined by:
$P^\Obs(A)=P(\Obs^{-1}(A))$.

\item Given a subsigmafield $\tribu'_\Omega$ of $\tribu_\Omega$, the conditional probability of $A$ ``given the most accurate desctiption of $\omega$ by means of statements of $\tribu_\Omega'$''\citep[See~][p.~]{Sverdrup1966Present} is denoted $P(A\mid \tribu_\Omega', \omega)$. For all $A\in\tribu_\Omega$, the function $\omega\mapsto P(A\mid\tribu',\omega)$ does always exist and is $P$-almost uniquely defined.
\item Given two random variables $\Rv_1$, $\Rv_2$, $P^{\Rv_1\mid\Rv_2}$ is a function: $(\somesubset,\omega)\mapsto 
P(\Rv_1^{-1}(\somesubset)\mid \tribu(\Rv_2),\omega)$. It satisfies:
$P^{\Rv_1\mid{\Rv_2}}\dependsonly\left((\somesubset,\omega)\mapsto(\somesubset,\Rv_2(\omega))\right)$, which defines $P^{\Rv_1\mid\Rv_2=\rv_2}:A\mapsto P^{\Rv_1\mid\Rv_2}(A,\omega)$ for $\omega\in{\Rv_2}^{-1}(\{\rv_2\})$, and $\rv_2\mapsto P^{\Rv_1\mid\Rv_2=\rv_2}$ is defined $P^{\Rv_2}$-almost uniquely. Given $\somesubset\in\tribu_\Omega$,
define $P^{\Rv\mid\somesubset}=P^{\Rv\mid\setindicator_{A}=1}$.
\item Let $(\Omega,\tribu_\Omega)$, $(\Omega',\tribu_\Omega')$ be a measurable space. A function
$P:\Omega'\to(\tribu_\Omega\to[0,1]),(\omega')\mapsto (\somesubset\mapsto P_{\omega'}(\somesubset))$
is a transition probability function from $(\Omega,\tribu_\Omega)$, $(\Omega',\tribu_\Omega')$ to $\forall \omega'\in\Omega'$ if and only if
$P_{\omega'}$ is a probability measure and if $\forall \somesubset\in\tribu_\Omega$, $\omega'\to P_{\omega'}(\somesubset)$ is measurable.
\item Given a measurable space $(\Omega,\tribu_\Omega)$, a probability space $(\Omega',\tribu_{\Omega'}, P')$ and a transition probability function $P$ from $(\Omega',\tribu_{\Omega'})$ to $(\Omega,\tribu_\Omega)$, the probability induced by $P$ and $P'$ is 
 $\left<\omega'\mapsto P_{\omega'}\mid P'\right>_{\Omega,\Omega'}: \tribu_\Omega\to, A\mapsto \int_{\Omega'} P_{\omega'}(A)~\derive P'(\omega')$.
\item Given two measures $P, P'$ on $\tribu_\Omega$, $\tribu_\Omega'$, respectively $P\otimes P'$ is the measure on the sigma-field on $\Omega\times\Omega'$ generated by the products of elements $A$ of $\tribu_\Omega$ and $B$ of  $\tribu_\Omega'$ such that for any such elements, $(P\otimes P')(A\times B)= P(A)\times P'(B)$.

\item Given $\somesubset\in \tribu_{\Omega}$, $P^{\Rv\mid \somesubset}$ is the probability distribution $P^{\Rv\mid \mathds{1}_\somesubset=1}:\somesubsetB\mapsto P(\Rv^{-1}(\somesubsetB)\cap \somesubset)$.
\item For a set of probability distributions $\mathscr{P}$ on $(\Omega,\tribu_{\Omega})$, 
$$E:
\begin{array}[t]{lcl}
    \range{P}&\to& ,\\
    P&\mapsto&\expected_P:
\begin{array}[t]{lcl}
    \left((\Omega,\tribu_{\Omega})\to \right)\cup \left((\Omega,\tribu_{\Omega})\to \right)^2&\to     & \left((\Omega,\tribu_{\Omega})\to \right),\\
    \quad\quad\quad\quad\quad\quad\quad(\Obs_1,\Obs_2)                  &\mapsto & \expected_P(\Obs_1\mid\Obs_2)=\int \obs_1\derive P^{\Obs_1\mid\Obs_2},\\
    \quad\quad\Obs                  &\mapsto & \expected_P(\Obs)=\int \Obs\derive P\end{array}\end{array}$$
    so $E[\Obs]:P\mapsto \expected_P(\Obs)$. The quantity $\expected_P(\Obs)$ may take the value ``undefined''.
\item When we write conditional probabilities or expected values, we implicitely assume that they are defined.
 Given a probability $ P$ on a measurable space ($\Omega,\tribu_\Omega$), random variables $\Obs_1$,$\Obs_2$, $\Obs_3$, 
 $ P^{\Obs_1\mid \Obs_2}=  P^{\Obs_1\mid \Obs_3}$ means that for all $f$ a positive measurable function on $(\Obs_1(\Omega),  P^{\Obs_1})$, $\expected_P(f(\Obs_1)\mid \Obs_2)\stackrel{ P-a.s.}{=}\expected_P(f(\Obs_1)\mid \Obs_3)$.
 \item Given a statistical model $(\Omega$, $\tribu_{\Omega},\left\{P\in\range{P}\right\})$, when writing density and conditional densities, we implicitely assume all necessary conditions for these densities to be defined, including that $\range{P}$ is dominated by a $\sigma$-finite measure $\dominant$.
 If $P$ is absolutely continuous with respect to $\dominant$, then for any random variable $\Rv$, 
 $P^\Rv$ is absolutely continuous with respect to $\dominant^\Rv$ (\citep[See][p.~229]{HalmosSavage1949}).
 The expression $f_{\Rv;P}$ will denote a density of $P^\Rv$ with respect to $\dominant^\Rv$.
 Given an additional random variables $\Rvv$, $f_{\Rv\mid\Rvv=\rvv;P}(\rv)$ is defined as as $f_{\Rv\mid\Rvv=\rvv;P}(\rv)=f_{\Rv,\Rvv;P}(\rv,\rvv)/f_{\Rv;P}(\rv)$ when $f_{\Rv;P}(\rvv)>0$, $0$ otherwise.
 We assume that densities, conditional densities and conditional distributions are chosen to be compatible, 
 so that the following properties are satisfied:
 \begin{enumerate*}
  \item $f_{(\Rv,\Rvv);P}(\rv,\rvv)=f_{\Rv\mid\Rvv=\rvv;P}(\rv)f_{\Rv;P}(\rv)=f_{\Rvv\mid\Rv=\rv;P}(\rvv)f_{\Rvv;P}(\rvv)$
  \item $P^{\Rv\mid\Rvv=\rvv}(A)=\int_A f_{\Rv\mid\Rvv=\rvv;P}(\rv)\derive \dominant^{\Rv\mid \Rvv=\rvv}(\rv)$.
  \item $f_{\Rv;P}(\rv)=$
  $\int_{\Rvv(\Omega)} f_{\Rvv,\Rv;P}(\rvv,\rv)\derive \dominant^{\Rv\mid \Rvv=\rvv}(\rv)$.
 \item $\int_{\Rvv(\Omega)} f_{\Rvv\mid\Rv=\rv;P}(\rvv)\derive \dominant^{\Rvv\mid \Rv=\rv}(\rvv)=1$.
 \end{enumerate*}
 %Given a separation of $\Omega$ by $(\Obs,\Comp\Obs)$, such that the distribution $\zeta$ is dominated by the distribution $\dominant^\Obs\otimes\dominant^{\Comp\Obs}$, 
 \item Given compatible densities $\{f_{.\mid.;P}\setsep  P\in\range{P}\}$, the likelihood 
of $\paramf$ in $\param$  given $\Obs_1=\obs_1$ conditionnally to $\Obs_2=\obs_2$ is defined as  
$\likelihood(\paramf;\Obs_1\mid\Obs_2)(\theta;\obs_1\mid\obs_2)=\sup\left\{\left.f_{\Obs_1\mid\Obs_2=\obs_2;P}(\obs_1)\right| P\in \paramf^{-1}(\theta)\right\}.$
The maximum likelihood estimator of $\paramf(P)$ is the (non necessarily non empty or reduced to one element) set $\argmax\left\{\mathscr{L}(\paramf;\Obs)(\theta;x)\setsep  \theta\in\paramf(\range{P})\right\}$.

\end{enumerate*}

%\todof{Define class of equivalence for a.s. equality}

% \begin{propertydefinition}[Transition probabilities](see \citet[p.~]{Sverdrup1966Present})\todomargin{Give page}
% Consider a measured space $(\Omega,\tribu,P)$, and a subsigmafield $\tribu'$ of $\tribu$, then define the conditional probability of $A$ given ''the most accurate description of $y$ by means of statements from $\tribu'$ as the value un $(A,\omega)$ of the function: $\tribu\times \Omega\to[0.1], (A,\omega)\mapsto P(A\mid \tribu_Y, \omega)$, which $P(A\mid \tribu_Y, .):\Omega\to[0,1],\omega\mapsto P(A\mid \tribu_Y, \omega)$ is the class of almost unique $\tribu'$-measurable functions of $\omega$ which satisfies $[\forall B\in \tribu'],\ P(A\cap B)=\int_B P(A\mid \tribu',\omega) \derive P(\omega)$. The function $P(A\mid\tribu',.)$ does always exist and is almost uniquely defined.
% \end{propertydefinition}
% 

%\subsubsection{Sufficiency, ancillarity, ignorability}

\section{A generic model for survey sampling and other areas of Statistics.}\label{sec:modelTYZ}
In Section \ref{sec:definitions}, definitions of ``at random'', ``R-ignorable'' and ``non informative'' will be given in a general statistical model. To illustrate this definitions with real-life models including many common models used when analysing survey data, we introduce in Section \ref{sec:statisticalframework.1} a more particular model that will be called the ``$\Obs=\obsf(\Transfo,\Transfo[\Signal],\Designvar)$''-model.
In section \ref{sec:statisticalframework.2}, we present usual models used in survey sampling and show how they can be written as particular cases of the  ``$\Obs=\obsf(\Transfo,\Transfo[\Signal],\Designvar)$''-model. We show that the ``$\Obs=\obsf(\Transfo,\Transfo[\Signal],\Designvar)$''-model can be used for design-based, model based or bayesian inference, missing data models, models with auxiliary information, selection with or without replacement ...etc.
Doing so we define a series of common random variables or parameters used in survey sampling, as the sample, the sample indicators, the design, the design variables, the inclusion probabilities of different orders, ...etc,  
In section \ref{sec:statisticalframework.3}, we present real-life statistical models that account for transformations of the signal not necessarily limited to survey sampling such as noise addition, that also fit the ``$\Obs=\obsf(\Transfo,\Transfo[\Signal],\Designvar)$''-model.

\subsection{The ``$\Obs=\obsf(\Transfo,\Transfo[Y],\Designvar)$'' model }\label{sec:statisticalframework.1}

%In this paper, we propose to define key concepts as at random, or informative in a general way. 

Define the ``$\Obs=\obsf(\Transfo,\Transfo[\Signal],\Designvar)$ model as the statistical model $\left(\Omega,\tribu_\Omega, P)_{P\in \range{P}}\right)$, where $\Omega$ is a subset of the product of three spaces: $\Omega\subseteq\left(\range{\Obs}\times\range{Z}\times\range{T}\right)$, for which $\range{T}\subseteq\left(\allfunc{\range{Y}}{}\right)$ and on which three random variables are defined, that correspond to the projections on $\range{\Obs}$,$\range{Z}$,$\range{T}$ with respect to the previous space product:
the {\bf signal} $Y:(\Omega,\tribu_\Omega)\to (\range{Y},\tribu_\range{Y})$, 
 the {\bf transformation conditioning variable} $Z:(\Omega,\tribu_\Omega)\to (\range{Z},\tribu_\range{Z})$, and the {\bf transformation mechanism} of the signal $\Transfo:(\Omega,\tribu_\Omega)\to(\range{\Transfo},\tribu_{\range{T}})$. Assume that  
 \begin{equation}\forall  P\in\range{P},\  P^{T\mid Y,Z}= P^{T\mid Z}.\label{eq:A1}\end{equation}

% Let $\Thetaset$ and $\Xiset$ be index sets in bijection with $\left\{ P^Y\mid  P\in\range{P}\right\}$, and 
%            $\left\{ P^{Z,T\mid Y}\mid  P\in\range{P}\right\}$ respectively, and denote $ P^Y_\theta$, $ P^{Z,T\mid Y}_\xiparam$ the indexed elements. 
%            Let $\Gamma$ be a subset of $\Thetaset\times \Xiset$ such that 
%            $\left\{( P_\theta^Y, P_\xiparam^{Z,T\mid Y})\mid (\theta,\xiparam)\in\Gamma\right\}=\left\{( P^Y, P^{Z,T\mid Y})\mid  P\in\range{P}\right\}$. The set
%            $\left\{ P^{Y,Z,T}\mid  P\in \range{P}\right\}$ is naturally indexed by $\Gamma$ and let denote by $ P_{\theta,\xiparam}^{Y,Z,T}$ its elements. 
%            We do not impose any assumption on the separability of $\Gamma$ into the product of a subset of $\Thetaset$ by a subset of $\Xiset$. Let $\Transfo[Y]$ denote the random variable $\Transfo[Y]:\Omega\to \range{X},\ \omega\mapsto (T(\omega)(Y(\omega))$.
The observations consist of the image by a deterministic and known function ``$\obsf$'' of the triplet $(\Transfo[Y], T, Z)$: $\Obs=\obsf(\Transfo[Y], T, Z)$, and the statistician works with the model $$\left(\codomain(\Obs),\tribu_{\codomain(\Obs)},\{P^\Obs\mid P\in\range{P}\}\right).$$

When doing Bayesian inference, a set $\{\priordist\in \priordistset\}$ of prior distributions is defined on $\range{P}$.

\subsection{Inference}

Let $\paramf:(\range{P}\times\range{\Obs})\to$, such that $\exists \paramf':\left(\left\{P^Y\mid P\in\range{P}\right\}\times \left\{Y(\omega)\mid \omega\in\Omega\right\}\right)\to$ for which $\paramf=\paramf'\circ\left(P\mapsto P^Y,Y\right)$. 

The target of the inference is $\paramf$: when $\paramf(P,\omega)$ is free of $\omega$ (e.g. $\paramf=\paramf'\circ\left(P\mapsto P^Y\right)$), the goal is to estimate $\paramf(P)$, whereas when $\paramf(P,\omega)$ is free of $P$ (e.g. $\paramf=\paramf'\circ Y$), the goal is to predict $\paramf(\omega)$. The model for the observations is:
    $\left\{ P^\Obs \mid P\in\mathscr{P}\right\}$.

\subsection{Examples of models used in survey sampling}\label{sec:statisticalframework.2}

As explained by \citet{PfeffermannSverchkov2009}, the observations from a survey are the output of two random processes: the population generation and the sample selection, so survey sampling models account for those two processes. The framework we described in Subsection \ref{sec:statisticalframework.1} allows to distinguish those two processes: $Y$ is the population generation, and the observation on a random selection of units is the observation of a random transformation $\Transfo[Y]$ of $Y$. 
In this subsection, we show that the simple general framework of Subsection \ref{sec:statisticalframework.1} can cover all situations encountered in survey sampling.
In survey sampling, there exists a diversity of paradigms : model based inference, design based inference, bayesian inference, that imply different population models (fixed population model, superpopulation models, ...), that are associated with different population and observation models.
Populations models specificities come from:
\begin{enumerate*}
 \item the population size that can be observed or unobserved, a parameter or a random variable,
 \item the population that can be finite, infinite, discrete, or continuous (for example for spatial data).
 \item the population characteristics can be considered random or a parameter.
\end{enumerate*} The model will also reflect specificities about what is observed and how samples were obtained:
\begin{enumerate*}
 \item the sample can also be finite, discrete, or continuous,
 \item the sampling can be with or without replacement,
 \item one may or may not be able to link sample units to population units,
 \item sampling weights may or may not be observed,
 \item design variables may or may not be observed,
 \item duplicates may or may not be identifiable.
 \item measurement may be with or without errors errors.
\end{enumerate*}
We propose a general model that embraces all those cases, we show that Design, inclusion probabilities, that can be observed, are functions of a design variable $Z$, whereas the correspondance between the sample index and the population index is a function of $\Transfo$, so that in the general case, observations are a function of $(\Transfo,\Transfo[\Signal], \Designvar)$.

\comments{
JS 2019 02 14 - C. Definition 1: Before this Y is only an abstract concept . It would be good to add comments about its use in the finite population context 
In [..] I do not know what $\Extension{U}$ would be.

DB 2019 02 21 - I added an example.

JS 2019 02 23 - 1.P15, Definition 12. In the last sentence, the population itself is unknown, leading to the use of the superset.

DB 2019 02 24 - You are right. U can be a parameter or a random variable.
}

\begin{definition}[Population model]\label{def:popmodel}
 When $\range{Y}$ is of the form $\bigcup_{\pop \in \superpop} (\allfunc{\pop}{\range{Y}'})$, where $\range{Y}'$ is a set and $\superpop$ is a collection of sets, and when $\exists U:\Omega\to\superpop$ is such that $\forall \omega\in\Omega$, $Y(\omega)\in (U(\omega)\to \range{Y}')$, we will say that $Y$ follows a population model and the random variable $\Pop $ will be called the population, $N:\omega\mapsto\size(\Pop(\omega)) $ will be the population size. %For $y\in \range{Y}$, $\size(y)$ will denote the cardinal of the domain of $y$. 
\end{definition}

\comments{
JS 2019 02 23 -  I’m not sure that the last two lines are helpful. Also, my concept of superpopulation is a little different, i.e., selection of a finite population from a stochastic model. You are, I think, emphasizing that the composition of a finite population may be unknown, as, e.g., when only sampled clusters are enumerated.

DB 2019 02 24 - Hi Joseph, I see your point. I just want to say that the population can be unknown. In this case, we can model it (model the size for example, or make it a parameter, like in capture recapture models.
I think that the problem you are explaining is that using superpopulation there does not correspond to what people call super population. I just deleted the sentence {This superset can be conis called superpopulation.}. 
I agree it is confusing and people would try to make a link to another concept. And it is not usefull to name it. I am not even sure one should talk about superpopulation model here, although I mentioned that is,}

\delete{\begin{example}
 Consider a population $\pop$, that is a set of individuals, $\pop$ being unknwown. For an individual indexed by $k\in \Pop$, assume $Y[k]$ is the age in years of the individual $k$, which can vary in $\range{Y}'=\mathbb{N}$.
 So conditionnally to $\Pop$, $Y$ is a random variable with value in $(\allfunc{\Pop}{\range{Y}'})$, and unconditionnally on $\Pop$, $Y$ is a random variable with value in the union of all sets $(\allfunc{\Pop}{\range{Y}'})$ whith $\Pop$ varying in the set of possible values for $\Pop$. \delete{It is common in survey sampling to consider that the set of possible values for pop is just a sample (subset) of all the values of a superset that englobes all the possible sets. This superset can be conis called superpopulation.} In design based inference, when the population size is unknown, it becomes an unknown parameter, that can be the parameter of interest, as for example in capture-recapture problems. In a model-based approach, one could have a stochastic model on $U$, based on demographic projections for example. In this case $U$ is a random variable, that can be latent.
\end{example}}

\comments{
JS 2019 02 14
D. Definition 2
In most of the paper, $\theta$ refers to parameters [...] of a "superpop" process. Using it in a different way as here, is not helpful to resolve. But perhaps this cannot be avoided - in terms of generality.}

\comments{
DB 2019 02 20
Thank you for this comment. I think you are right, that some precautions are needed here. My goal is just to say that fixed population model fits in the usual parametric model. This is what I meant by "Definition \ref{def:popmodel} embraces the two cases". Do you have a suggestion to improve the writing ?.
}

In survey sampling, we usually consider selection under a finite population model. Frequentist inference consider two parametric types of population models refered to as fixed population model for design-based inference
\citep{CasselSarndalWretman1977}, \citet[p.~52]{Gourieroux1981theoriedessondages})
and superpopulation model for model-based inference.
% , the first one being a special case of the second,   
% when the second is not restricted to the iid case. A specificity of the fixed population model is that  $Y$ is considered non random. Definition \ref{def:popmodel} embraces the two cases.
There exists different interpretations of a superpopulation model. Nevertheless, in practice, it consists in using a model on $Y$ that is not the fixed population model. \cite{Nathan2011Superpopulation} gives more details about the origin of the term superpopulation and the design-based versus model-based inference debate.
For our purpose, which is only to propose general definitions, the only thing that matters is to consider a framework such as \ref{def:popmodel} general enough to fit any population model, a fixed population model or a superpopulation model.

\comments{
JS 2019 02 23 - Definition Fixed population model: I don’t see any way to make this more consistent with usual survey sampling notation. It is important, as you have done here, to explain the context, i.e., finite population model.

DB 2019 02 24 - Thank you for your comment. This definition is very close to \citet[p.~52]{Gourieroux1981theoriedessondages}, but it is in french, In this sense it is a classic. Gourieroux wrote a book about survey sampling although he is not a specialist, and the first thing he did was to write down the statistical model.  Unfortunately, I could not find this nice definition anywhere else.
}

 \begin{definition}[Fixed population model]\label{def:fixedpopmodel}
The fixed population model corresponds to the special case where $\left\{ P^Y\mid P\in\range{P}\right\}=\left\{\Dirac_{\{\signal\}}\mid \signal \in \range{\Signal}\right\}$ where $\range{Y}=\bigcup_{\pop\in\superpop}\allfunc{\pop}{\mathscr{Y'}}$. 
\end{definition}
\comments{JS 2019 03 20. Why not use $\bigcup_{\pop\in\superpop}\allfunc{\pop}{\mathscr{Y'}}$
You do not use caps/lower cases consistently. (see, e.g. def 10) and this may induce confusion.

DB 2019 04 18: I made the change: captital: random variable. lowercase: non random. bold greek: parameter function, greek: parameter or function.
}

Another example of fixed populations are fixed excnageable population models, (see \cite{Sugden1979}), that are models where the population index does not carry any information. 
Under an exchangeable population model, forall population size $\pop\in\range{\Pop}$, $P^{Y\mid \Pop=\pop}$ is invariant by any permutation of $\pop$, e.g. for all $P\in\mathscr{P}$, for all permutation $\sigma$ of $\pop$, one must have $P^{Y\mid \Pop=\pop}=P^{Y[\sigma]\mid \Pop=\pop}$.
\begin{example}[Fixed exchangeable population model ]
An exchangeable fixed finite population model is of the form $\left\{P^Y\setsep  P\in\range{P}\right\}=\left\{(\size(\domain(\signal))!)^{-1}\sum_{\sigma\in\somesubset_\signal}\Dirac_{\signal\circ\sigma}\setsep  \signal\in \range{Y}\right\}$, where $\range{Y}$ is defined as in Definition \ref{def:fixedpopmodel}, $\superpop$ contains only finite non empty sets, and here $\somesubset_\signal$ is the set of permutation of the codomain of $\signal$.
\end{example}

\comments{JS 2019 03 19: Not clear why this is helpful. Exchangeability  not usual in this context, seen frequently in Bayesian inference.

DB 2019 04 20: You are right. I want to give examples in last section with exchangeability, and I make more comment in the supplemented materials.}
\delete{ 
\begin{example}
Consider the income of US households at a certain time $t$: $\range{Y}'=\mathbb{R}$, $N$ is the number of households  at a time $t$. For a design-based analysis, the population size $N$ is a function of the parameter: $N=\paramf(P^Y)$ with $\paramf(\Dirac_\signal)=\size(\domain(\signal))$. For a model-based analysis, some analyst may want to put a stochastic model on $N$, for example, $N(t)=\beta_0 N(t-1)+\beta_1 N(t-2)+\varepsilon(t)$, and a parametric model on the distribution of household income. Then $N$ is not a parameter of the model to estimate, but a random variable to predict.
\end{example}}

%\comments{ Next definition is for partial non response, not fundamental, should certainly be replaced by a comment}

% \begin{definition}[Sample selection from the population with partial response]\label{def:selection} \delete{
% Assume that $Y$ follows a population model, and that $\range{Y}'$ is the direct product of $p$ spaces 
% $\range{Y}'=\range{Y}'_1\times \ldots\times\range{Y}'_p$, then we will say that $T$ is a (partial) selection from the population when 
% $ P^{(T,Y)}-a.s(t,y)$,  $\exists n\in\mathbb{N}, k_1,\ldots k_n\in\{1,\ldots,\size(y)\}$ such that  $t(y)=(y_{k_i})_{i\in\{1,\ldots,n\}}$ 
% (with or without replacement) of the coordinates of $y$.}
% \end{definition}

\comments{
 
JS 2019 02 20 Questions about definition 9 and example 12.
Transformation $\Samplepopmap(\omega)$ transforms sample to population indices.
Definition 9 seems to be clear. Using the simpler case: $\Transfo[Y] = (Y(R1), \ldots Y(Rn))$ assumes fixed $Y$ and transforms using $\Samplepopmap(\omega)$.

DB 2019 02 21 - I was also wrong there, I replaced $Y(R[1])$ by $Y[R[1]]$...

JS 2019 02 24 - I commented about these changes earlier today.

}

\begin{definition}[Sample selection from population]\label{def:selection} Assume that $Y$ follows a population model, then  
we will say that $T$ is a {\bf selection} or {\bf sample selection} (from the population) when 
$\forall \omega\in\Omega$,  
$\exists$ a set $\Sampleindexset(\omega)$ called sample index, a function $\Samplepopmap(\omega)$ from $\Sampleindexset(\omega)$ to $\Pop (\omega)$ called a sample to population indexes mapping, such that: 

$$
\Transfo:
\begin{array}[t]{lcl}
    \Omega&\to    &,\\
    \omega&\mapsto&T(\omega):
\begin{array}[t]{lcl}
    (\Pop(\omega)\to\range{Y}')&\to     & ,\\
    \quad\quad\quad\signal                  &\mapsto & (T(\omega))(\signal):
\begin{array}[t]{lcl}
    \Sampleindexset(\omega)&\to&\range{Y}',\\
    \ell&\mapsto& y\left(\left(R\left(\omega\right)\right)\left(\ell\right)\right).
                                       \end{array}\end{array}            
           \end{array}$$ Which, by using the $[.]$ notations, can be written when $\Pop$ is constant, $\Transfo[\signal]=\signal\circ R$. 
In fine, $\Transfo[\Signal]=\Signal[R]$, and when $P^\Sampleindexset$-a.s$(l)$, $l=\{1,\ldots,\samplesize\}$, for $\ell \in\{1,\ldots,\samplesize\}$, $(Y[R])[\ell]=Y[R[\ell]]$ can be interpretated as the value of $Y$ for the individual selected at the $\ell$th draw.
% Let $n=\size[\Sampleindexset]:\omega=\size(\Sampleindexset(\omega))$ denote the sample size. When $P^n$-a.s. $(\samplesize)$,  $\samplesize<+\infty$, then without loss of generality, $\Sampleindexset=\{1,\ldots, n\}$, and $\Transfo[Y]$ can be written:
% $\Transfo[Y]=(Y[R[1]], \ldots, Y[R[n]])$. 
When $ P^R-a.s(r)$, $r$ is injective, we say that the selection $R$ is \textbf{without replacement}.
% \delete{
% (in words, $t(y)$ is a vector whose coordinates are the output of a sample (with or without replacement)
%  of the coordinates of $y$. When $ P^{(T,Y)}-a.s.(t,y)$, $k_1,\ldots k_n$ are all distinct, we will say that the selection is without replacement (with replacement in the opposite case).}

\end{definition}

\comments{JS 2019 03 19: Too much details in $T:\Omega\to\ldots, \ell\mapsto y(\Samplepopmap(\omega))(\ell))$ I thought that the text is fine.

DB 2019 04 17: Thanks for the comment, we can discuss to make it less detailed.}

\comments{
 
JS 2019 02 20
Example 12 isn’t so clear.$T(\omega):\mathbb{R}^{\{1,\ldots,8\}}\to \mathbb{R}^5, y\mapsto (y(3),y(1),y(5),y(3), y(2))$.
The latter expression has two stages and seems to duplicate the function in definition 9, i.e., \Transfo[Y], in the simple case (except for your use of lower case letters in only one).
So, in your shorthand, what is T?
I feel that there are two processes, i.e., one selecting Y from a superpopulation, and, given Y, selecting the sample indices and associated values of the variable y. However, it’s not obvious how this is translated into these expressions.

DB 2019 02 21 - I rewrote this part. Is it better now?

JS 2019 02 24 -  I commented about these changes earlier today.
}

\begin{definition}[Selection count $J$ and selection indicator $I$]
For $\omega\in\Omega$, we define the (without replacement) sample $\Sampleset(\omega)=\Samplepopmap(\omega)(L(\omega))$ as the image of the function $\Samplepopmap(\omega)$. We define the selection indicator $I:\omega\mapsto I(\omega):\powerset(\Pop(\omega))\to \{0,1\},\somesubset\mapsto (I[\somesubset])(\omega)=1$ if $\somesubset\cup \Sampleset(\omega)\neq\emptyset$ and $0$ otherwise, 
and $J:\omega\mapsto J(\omega):\powerset(\Pop(\omega))\to \mathbb{N},\somesubset\mapsto J[\somesubset](\omega)=\size\left(\left(\Samplepopmap(\omega)\right)^{-1}(\somesubset)\right)$.
\end{definition}

\begin{example}\label{Ex:6}
 Assume a population model. For $\individual\in\Pop$, $Y[k]$ is the income of individual $k$. Assume a selection with replacement and random size.
 Let $\omega\in \Omega$ and denote $t=T(\omega)$, $r=\Samplepopmap(\omega)$ and assume $N(\omega)=8$ and $n(\omega)=5$, so $\Pop(\omega)=\pop=\{1,\ldots 8\}$, and $\sampleindex=\Sampleindexset(\omega)=\{1,\ldots,5\}$. Assume the units drawn for the draw 1 to 5 were $r(1)=3$, $r(2)=1$, $r(3)=5$, $r(4)=3$, $r(5)=2$ respectively, so 
 $r$ is the function $:\sampleindex\to \pop,~1\mapsto 3,~2\mapsto 1,~3\mapsto 5,~4\mapsto 3,~5\mapsto 2$ and 
  $t$ is the function $:\allfunc{\left(u\to\mathbb{R}\right)}{\left(\allfunc{\sampleindex}{\mathbb{R}}\right)}, y\mapsto \transfo(\signal)=\left(y\circ r:
 \sampleindex\to\mathbb{R},\ell\mapsto y(r(\ell))\right)$. Note that the selection is with replacement as $r(1)=r(4)=3$. The sample indicator vector is $(I[1](\omega),\ldots, I[8](\omega))=(1,1,1,0,1,0,0,0)$ and the sample count vector is $(J[1](\omega),\ldots J[8](\omega))=(1,1,2,0,1,0,0,0)$.
\end{example}

\comments{
JS 2019 03 19.Note typo on that line i.e. should be $r(4)=r(1)=3$.

DB 2019 04 17. Thank you , I corrected this line.
}

\begin{definition}[Design and design variables]\label{def:selection} 
A design is a random variable $D$ such that  $\forall \omega \in \Omega$, $\Design(\omega)$ is a probability measure
on $\bigcup_{l\in\mathscr{L}} \{l\to\Pop (\omega)\}$, where $\mathscr{L}$ is a set of possible sample indexes,  and $\forall P\in\range{P}$, 
\begin{equation}\label{prop:PRmidDdependsonlyD}
P^{\Samplepopmap\mid \Design,\Signal,\Designvar}\stackrel{ P-a.s.}{=}\Design.
\end{equation}
\end{definition}

\comments{JS 2019 03 19: Is all the notations needed ? Why $I:\Omega$ ?...

DB 2019 04 17: Thank you for the comment,
I simplified.
}

\permanentcomments{$ P^{R\mid D,Y,Z}=D$: we have a random distribution on one side, and a conditional probability in the other side. 
More rigourous notation would be $$ P^{D,Y}\text{-a.s}(d,y,z),\  P^{R\mid D=d,Y=y,Z=z}=d.$$ as distributions conditional to a r.v. $X$ are defined $ P^{X}\text{-a.s}$}

\begin{example}[Simple random sampling]
Let $\samplesize,\popsize\in \mathbb{N}$. Simple random sampling without replacement of size $\samplesize$ between $\popsize$ is the uniform distribution on the set (denoted here by $\somesubset$) of functions from $\{1,\ldots,n\}$ to $\{1,\ldots,N\}$: 
$\design=(\popsize)^{-\samplesize}\sum_{r\in \somesubsetB}\Dirac_r$, whereas  when $\samplesize\leq \popsize$, simple random sampling without replacement of size $\samplesize$ between $\popsize$ is the uniform distribution on the set (denoted here by $\somesubsetB$) of injective functions from $\{1,\ldots,\samplesize\}$ to $\{1,\ldots,\popsize\}$: 
$\design=\samplesize!(\popsize!)^{-1}\sum_{r\in \somesubsetB}\Dirac_r$.
\end{example}

\comments{
JS 2019 03 19. Example is clearer than before but still hard to read beacause of the notation. I suggest removing $\omega$ for simplification. This would be especially helpful for definition.

DB 2019 04 17. Thanks for the comment, I removed the $\omega$ in the expressions after defining all the values obtained for each random variable for a specific value of $\omega$.

JS 2019 03 19.I found $\sum_{r\text{ injective}...}$ to be very awkward, especially the world injective in a subscript. [...]

DB 2019 04 17.
THank you for your comments, I removed the injective from the subscript.
}

This series of examples illustrate that our general framework applies to any situation, and allow to discriminate between the different mathematical objects which ones are a function of $\Signal$, which ones are a function of $\Transfo$ and which ones are a function of $\Designvar$.

It is important to notice that the population is a function of $\Design$, e.g. $\Pop\dependsonly\Design$,  as $\forall\omega$, $\Design(\omega)$ is a probability distribution on a set of mappings to $\Pop(\omega)$. So when conditioning on $\Design=\design$, $\Pop$ is not random. When defining inclusion probabilities, it is then not necessary to condition on both the design and the population.

Because $\Pop $ can be deduced from $\Design$, the expression $\expected_P(I[\individual]\mid D)$ has a meaning. The expression  $\expected_P(I)$ may not have sense as we are taking the expected value of a random vector that take values in spaces of different dimensions when $\Pop$ is random, but $\expected_P(I\mid \Pop)$ does and so does $\expected_P(I\mid \Design)$ as $\Pop\dependsonly\Design$. 

\begin{definitions}[Inclusion  and selection probability densities] 
For each $\pop\in\range{\Pop}$, consider a $\sigma$-field $\tribu_{\pop}$ on $\pop$, and a measure $\countingmeasure_\pop$ on $\tribu_{\pop}$. The functions 
$:\tribu_{\pop}\to, A\mapsto \expected_P(I[A]\mid \Design=\design)$ and $:\tribu_{\pop}\to, A\mapsto \expected_P(J[A]\mid \Design=\design)$ are probability measure and measure respectively that may admit Radon Nykodym derivatives with respect to $\countingmeasure_{\Pop(\omega)}$. Equation \eqref{prop:PRmidDdependsonlyD} implies that these measures do not depend on $P$, but depend on $\design$ only: $\expected_P(I[A]\mid \Design=\design)=\design\left(\left\{\samplepopmap:\to u\mid A\cap\image(\samplepopmap)\neq\emptyset\right\}\right)$
and $\expected_P(J[\somesubset]\mid \Design=\design)=
\int(\size\left(\samplepopmap^{-1}(\somesubset)\right)\derive\left(\design(\samplepopmap)\right)$.

\begin{enumerate}\item
The inclusion probability density function is the Radon Nikodym derivative $\Inclusionprob(\omega)=\derive(A\mapsto \expected_P(I[A]\mid \Design=\Design(\omega)))/\derive\countingmeasure_{\Pop(\omega)}$ when defined.
\item
The selection density function is the Radon Nikodym derivative $\Inclusionexpec(\omega)=\derive(A\mapsto \expected_P(J[A]\mid \Design=\Design(\omega)))/\derive\countingmeasure_{\Pop(\omega)}$ when defined.
\end{enumerate}
\end{definitions}

If $\forall \pop\in\range{\Pop},$ $\pop$ is countable and $\countingmeasure_\pop$ is the counting measure, then $(\Inclusionprob[{k}])_{ k\in \Pop}$ is the inclusion probabilities of each unit being on the sample conditionnally on the design $\Design$ :
 $\Inclusionprob[k]=\expected(I[\individual]\mid \Design=\Design(\omega))$.

\comments{DB 
Later I want to use the previous definition in an example: informative core sampling with uniform $\Inclusionprob$ function.

}

\begin{property}[The sum of inclusion probabilities equal the sample size]
Let $\omega\in\Omega$, let $\pop=\Pop(\omega)$, $\design=\Design(\omega)$,$\inclusionexpec=\Inclusionexpec(\omega)$ and $\inclusionprob=\Inclusionprob(\omega)$, then $\forall P$, $\int_\pop \inclusionprob\derive\countingmeasure_\pop=\expected_P(n_{\text{wor}}\mid \Design=\design)$ and $\int_\pop \inclusionexpec\derive\countingmeasure_\pop=E[n \mid \Design=\design] $, where $n_{\text{wor}}=\sum_{\individual\in\Pop}\Sampleind[\individual]$ is the sample without duplicates.
\end{property}

\begin{definition}[Design variable]
Often the design is a function of design variables, denoted $Z$, that are variables defined on the population, and $D= P^{R\mid Z}$.
\end{definition}

% 
% \begin{definition}[Design variable, design]\label{def:selection} \delete{
% Assume that $Y$ follows a population model, then the function $D$ such that 
% $D(z)= P^{T\mid Z=z}$ is called the design function, $ P^{T\mid Z}$ the design (or sample scheme) and $Z$ the design variable.}
% \end{definition}

\begin{example}[Stratified sampling]
Consider a population model such that $\range{\Signal}'=\mathbb{R}$, let $H\in\mathbb{N}$, let $Z$ be a random variable such that $\forall\omega$, $Z(\omega): \Pop(\omega)\to\{1,\ldots,H\}$.
Define $N_h(\omega)=\#\{\individual \in\Pop(\omega)\mid Z[\individual](\omega) =h\}$, and let $\Samplesize_h$ be a random variable such that $\forall\omega. \Samplesize_h(\omega)\in\{0,\ldots,\Popsize_h(\omega)\}$ and such that $\sum_h\Samplesize_h=\Samplesize$. Define $D(\omega)$ as the uniform distribution on 
the set of injective functions $r$ from $\{1,\ldots, \Samplesize(\omega)\}$ to $\Pop(\omega)$ such that $\forall h\in \{1,\ldots,H\}$, 
$\#\{\ell\in\{1,\ldots,\Samplesize(\omega)\}\mid Z_{r(\ell)}=h\}=n_h(\omega)$.
Then $D$ is the stratified sampling based on the stratification $Z$ with allocation $n_1,\ldots,n_h$ and simple random sampling within strata.
\end{example}

%We show below how this model is adapted to different situations in sampling. 

\begin{remark}[Observations]
Observations,denoted $\Obs$, are a function of $(\Transfo[\Signal],Z,T)$. We show that this notation includes many possible scenarii.
In practice, different situations may happen, we may or may not observe the population size, we may or may not observe the 
index of the values $\Signal[\individual ]$, we may or may not observe the number of times each unit was selected, we may or may not observe the inclusion probabilities for the whole population or the sample units. 
We may or may not oberve the design that was used. We may or may not observe auxiliary information.
If $\Obs=\Transfo[\Signal]$ and the sample is with replacement, and unequal probabilities, then one is not able to identify duplicates.
In the case $\Obs=(\Transfo[\Signal], T$, one can deduce $R$ from $T$, so be able to identify the units and the duplicates.
If $\Obs=(\Transfo[\Signal],T, D)$ (D is a function of $Z$), then one can deduce $\Pop $, the inclusion probabilities, the double inclusion probabilities, all the units.
If $\Obs=(\Transfo[\Signal], \Inclusionprob(R[L])$, one cannot identify the sampled units, and knows the inclusion probabilities for the sampled units only.
When the sample is without replacement and the population index observed on the sample, $\Transfo[\Signal]=\Signal[\Samplepopmap]=\Signal\restrict_{[\Sampleset]}$.
\end{remark}

\begin{definition}[Models on the observations]
All different models for inference on survey data can be written with the same expression:
 $\left\{P^{\Obs}\mid P\in\range{P}\right\}$. Thanks to this simple expression, we can define the different concepts of interest without depending on a particular case.
\end{definition}
\begin{example}[\protect{Parametric model for design based inference, \citet[p.~52]{Gourieroux1981theoriedessondages}}]
In the classic case of finite sample selection from finite population, with fixed population $\pop$ and fixed design $\design$ model, is the combinaison of a parametric model for the population of the form $\{P^\Signal\mid P\in\range{P}\}=\left\{\Dirac_\signal\mid \signal\in\range{\Signal}=\bigcup_{\pop\in\range{\Pop}}\left\{\allfunc{\pop}{\range{\Signal}'}\right\}\right\}$, and of a 
   parametric model for the design $\{P^\Design\mid P\in\range{P}\}$ such that $\range{\Design}_\pop\subset\left\{\Dirac_\design\mid \design\in\{\text{probability measures on }\left(\bigcup_{n\in\mathbb{N}}\left(\{1,\ldots,n\}\to\pop\right)\right)\}\right\}$.
If for example the observation is $\Obs=\Transfo[\Signal]=\Signal[\Samplepopmap]$, then the parametric model for the observation is:
 $\left\{P^{\Obs}\mid P\in\range{P}\right\}=\left\{d^{g_\signal}\mid \signal\in\range{\Signal}, \design\in\range{\Design}_{\domain(\signal)}\right\}$, where 
 $g_\signal:\left(\samplepopmap:\{1,\ldots,n\}\to\domain(\signal)\right)\mapsto \signal\circ\samplepopmap$. When $\design$ and $\pop$ are known, the model is then just  
 $\left\{P^{\Obs}\mid P\in\range{P}\right\}=\left\{d^{g_\signal}\mid \signal\in(\pop\to\range{\Signal}')\right\}$.
\end{example}

\subsection{The ``$\Obs=\obsf(\Transfo,\Transfo[\Signal],Z)$'' model applies to other areas of statistics}\label{sec:statisticalframework.3}

\comments{JS 20190204: 2. The examples in section 2.3 probably shouldn’t be included, especially for our emphasis on survey sampling.

JS 2019 02 23 -  first highlighted box. It may be better to tighten the presentation by limiting to survey sample examples. This can be decided later

DB 2019 02 24 - I think we need to put examples that show that this model is not limited to survey sampling. I think the mistake here is to talk about informativeness or so. We should give a limited of examples showing that (T,Y, Z) is very common.
The same examples can be used after section 3. I agree on limiting the number, and deciding later which ones are pertinent, which ones are not.}

Selection is not the only possible transformation that one can think of. In Statistics, there is usually a true phenomenon of interest, and we can perceive some signal, that can be perturbated.
The definitions we propose in Section \ref{sec:statisticalframework.1} are general enough to embrace this situation.

\comments{JS 20190204 line 1b. $\theta$ is not defined.

DB 2019 02 24 I had to rewrite everything with the new notations.
}
\begin{example}[Measurement error]
  Assume $Y$ is  a dichotomous random variable that we observe for each units of a known population of size $N$ with measurement error.
  The random variable domain is $\range{Y}=\{1,\ldots,N\}\to\{0,1\}$, and define $Z$ as a measurement error indicator: $Z:(\Omega.\tribu_\Omega)\to\left(\{1,\ldots,N\}\to\{0,1\}\right)$.
%  Assume $\left\{ P^Y\mid P\in\range{P}\right\}=\left\{\setindicator_\theta\mid \theta:\allfunc{\{1,\ldots,N\}}{\{0,1\}\right\}}$, 
  Assume independent measurement error for each unit with probability $\xiparam$, and error independent on $Y$:
  $\forall P\in \range{P}, P^{(Y,Z)}=P^{Y}\otimes P^{Z}$ and $\forall P\in\range{P}$, $\exists \theta\in[0,1]$ such that $\forall z:(\{1,\ldots,N\}\to\{0,1\}$ $P\left(\left\{Z=z\right\}\right)=\Comp\theta^{\sum_{\individual =1}^N z(\individual )}(1-\Comp\theta)^{\sum_{\individual =1}^N (1-z(\individual ))}$. Then $\forall \omega\in\Omega$, $T(\omega):\range{Y}\to\range{Y},y\mapsto y\times Z(\omega)+(1-y)\times(1-Z(\omega)$, so that $\Transfo[Y](\omega)=Y(\omega)Z(\omega)+(1-Y((\omega)))\times(1-Z(\omega)$..
\end{example}

\comments{JS 2019 03 19: $Z$ is  usually used as related to design

DB 2019 04 20. I agree. $Z$ plays the role of an auxiliary variable, used for the design. In the case where we deal with another type of trasnformation (noise addition) $Z$ is auxiliary information related to the transformation. In this sense the notation is consistent.}
\begin{example}
Assume $Y$ is the trajectory of a rocket. We observe the signal $\Transfo[Y]=(Y_t+\varepsilon_t)_{t\in \{t_1,\ldots, t_{10}\}}$. In this case, we can see that we do not observe all the points of the trajectory, and some noise has been added in the signal. The determination of the points $t_1, \ldots, t_{10}$ can be informative, for example, we may not observe the rocket when it is deviated. The noise may also be more important when the rocket is deviated as the sensors may be perturbated by the origin of the deviation. Asking the question : Is the transformation informative ? is pertinent as well.
\end{example}

\begin{example}
\citet{Heitjan1991}, \citet{Jacobsen1995} define the notion of coarsening  at random. This is also a case where the observation is the result of a perturbated measure, where the perturbation in this case is made willingly to protect privacy.
 
\end{example}

\begin{example}
\citet{Wubailey2016}  describe right censoring in the case of clinical trials, and call this process informative. Our general framework also applies
\end{example}

\begin{example}
Consider a stochastic process (for example, electricity consumption of a group of households) measured on a grid, where the interval between two measuremnts will depend on the previous measurements (interval length will be reduced when volatility will be high, or when high values will be reached).
The random process we are interested in is the consumption over the time, but the signal we obtain is a transformation of this signal. Our general framework also applies. The censoring is a random transformation that transforms the uncensored data into the censored data.
\end{example}

\section{When observations are the outcome of two or more random processes}\label{sec:definitions}
%\section{Definition of ``At random'', ``R-ignorable'' and ``Informative''}\label{sec:definitions}

It is common in survey sampling to decompose the random process behind the observations into distinct random processes: the variable of interest generation, the design variable generation, the sample selection according to the design, the non response, the measurement error, ...etc. This has been explained by \cite{Pfeffermann1998a} as well as \cite{Skinner1994}. \cite{Rubin1976} explains what it means to ignore a particular random process (the one that causes missing data) for a particular type of inference (likelihood based inference).
In this section, we propose to generalise the definition of ignoring one of the random processes at the origin of the observations, as well as 
a to-go guide to define the concepts of informative, at random, or ignorable in different situations. The first step is to start from a model, and to transform this model by ignoring one of the generating processes. The second step is to specify the target of the inference, and to make sure that this target has a meaning in the transformed model. The third step is to compare the inference on the target in the original model with the inference on the transformed target in the transformed model. Based on the type of inference and the type of target, the comparison may be based on different criteria (same likelihood, same properties of an estimator, same posterior distribution...). 
We propose to define the concepts of ``At random'', ``R-ignorable'' and ``Informative'' in a statistical framework that is even more general than the ``$\Obs=\obsf(\Transfo,\Transfo[\Signal],Z)$`` model of Section \ref{sec:modelTYZ}.
The ``$\Obs=\obsf(\Transfo,\Transfo[\Signal],Z)$`` model will only be used to illustrate these definitions and to establish the connections with existing definitions that were given in restricted frameworks that fit in the ``$\Obs=\obsf(\Transfo,\Transfo[\Signal],Z)$`` model.

\subsection{General statistical framework}
Consider a statistical model $\left\{(\Omega,\tribu_\Omega,P)\right\}_{P\in\range{P}}$, where $\Omega$ is a  set, $\tribu_\Omega$ is a sigmafield on $\Omega$ and $\range{P}$ is a set of family distributions on the measured space 
$(\Omega,\tribu_\Omega)$.
Consider three  random variables $\Id_\Omega$, $\Obs:(\Omega,\tribu_\Omega)\to(\range{\Obs},\tribu_\range{\Obs})$, the observed statistic (e.g. the variable that the statistician will observe), and $\Rv:(\Omega,\tribu_\Omega)\to(\range{\Rv},\tribu_\range{\Rv})$, that can be latent, or not. A non latent variable being defined as any $\Rv$ variable such that $\Rv\dependsonly\Obs$.
Let $\tribu_\mathscr{P}$ be a sigmafield on $\mathscr{P}$ and let $\priordistset=\{\priordist\in\priordistset\}$ be a set of probatility measures on $(\mathscr{P},\tribu_\mathscr{P})$, called prior distributions. %Let $\paramf:\allfunc{\range{P}}{}$.
Let $\Comp\Rv$ be a (non necessarily distinct) complement of $\Rv$, that can be latent or not. Since $(\Rv,\Comp\Rv)\isogeq\Omega$, $\exists \obsf:(\Rv(\Omega)\times\Comp\Rv(\Omega))$ such that $\Obs=\obsf(\Rv,\Comp\Rv)$, $\obsf$ is given by the relationship: $\obsf=\Obs[\interclass_{\Rv,\Comp\Rv}]$.

\begin{examples}
In the ``$\Obs=\obsf(\Transfo,\Transfo[\Signal],\Designvar)$'' model:
\begin{enumerate*}
 \item $\Id_\Omega=(T,Y,Z)$.
 \item When $\Transfo$ is a selection, the random variable $\Comp\Rv$ can be an image of $Z$, (for example $D$, or $\Inclusionprob$), or an image of $\Transfo$ (for example $I$, $J$ $\Samplepopmap$ ...) a function of $\Signal$ ...etc. The question can we ignore $\Comp\Rv$ is then equivalent to can we ignore the transformation process $(\Comp\Rv=\Transfo)$ ? the inclusion probabilities ($\Comp\Rv=\Inclusionprob$) ? both the inclusion probabilities and the transformation ($\Comp\Rv=(\Inclusionprob,\Transfo$) ? anything but the values of the signal and the design variable on the sample ($\Rv=(\Signal_\individual,\Designvar_\individual)_{\individual\in\Pop}$)?
 %\item $\paramf$ can be $P\mapsto\int Y\d P^Y$, $P\mapsto \mathrm{quantiles}(P^Y)$, ...etc
\end{enumerate*}
\end{examples}

\begin{examples}[Complements of $\Rv$]
In the model $\Obs=\obsf(\Transfo, \Transfo[\Signal],\Designvar)$, 
\begin{enumerate*}
 \item if $\Rv=\Transfo$, then a complement is $\Comp\Rv=(Y,Z)$. 
 \item if $\Rv=Z$, then a  complement is $\Comp\Rv=(Y,T)$.
 \item if $\Rv=\Transfo[\Signal]=\Signal[\Samplepopmap]$, a natural complement may be $(\Pop,\Signal\restrict_{[\Pop\setminus\image[\Samplepopmap]]},\Samplepopmap,\Transfo)$, 
\end{enumerate*}
\end{examples}

In the following, $\Comp\Rv$, will be called the nuisance process and $\Rv$ will be the process of interest.
 Following \cite{Rubin1976}, we propose to define what it means to "ignore $\Comp\Rv$ in this section.

\subsection{Transformations induced by the ignoring of a random process}

%Assume a (non necessarily distinct) complement $\Comp\Rv:(\Omega,\tribu_\Omega)\to(\Comp{\range{\Rv}},\tribu_{\Comp{\range{\Rv}}})$  of a random variable $\Rv:(\Omega,\tribu_\Omega)\to(\range{\Rv},\tribu_{\range{\Rv}})$, 

We assume that  the sigmafield $\tribu_{\Rv(\Omega)}$ contains all the elements $\{\rv\}$ for $\rv\in \Rv(\Omega)$, and such that the sigmafield $\tribu_{\Comp{\range{\Rv}}}$ contains all the elements of the form $\Comp\Rv(A)$ for $A\in\tribu_\Omega$.

\subsubsection{Transformation of the model}
Ignoring $\Comp\Rv$ means that one also ignores its probability distribution, and there are different ways to do it. A first way would be to consider $\Comp\Rv$ as a fixed (non random) process. A second way would be to replace its distribution given by the model by any distribution. It does not really matter as the choice of the distribution of $\Comp\Rv$ should not have an impact after $\Comp\Rv$ is ignored in the inference. Let $\range{P}'$ be the set of possible distributions of $\Comp\Rv$ after ignoring the process.

\begin{examples}
\begin{enumerate*}
\item Considering that $\Comp\Rv$ is fixed is equivalent to chosing: $\range{P}'=\left\{\Dirac_{\Comp\rv}\mid\Comp\rv\in\Comp\Rv(\Omega)\right\}$
\item Considering that the distribution of $\Comp\Rv$ does not matter is equivalent of chosing $\range{P}'=\left\{P'\right\}$, with $P'$ an arbitrary distribution on $\Comp\Rv(\Omega)$.
\item It is also possible to define $\range{P}'=\left\{P^\Rv\mid P\in\range{P}\right\}$.
\end{enumerate*}
\end{examples}

More importantly, the second step when ignoring $\Comp\Rv$ is to cut all stochastic dependence between $\Rv$ and $\Comp\Rv$. It may not be possible to consider the two processes as totally independent, especially when $\image(\Rv,\Comp\Rv)\neq\image(\Rv)\times\image(\Comp\Rv)$. %For example when ignoring the missing data mechanism, one considers the missing data as a non random variable, but still, the observations are a function of this variable, as it consists only of the values for the non observed elements.

\begin{definition}[Transforming a model to ignore a process]

 For $\rv\in \Rv(\Omega)$, 
define 
 the measurable set $$\ignoreset_{\rv,\Rv,\Comp\Rv}=\Comp\Rv^{-1}\left(\Comp\Rv\left(\Rv^{-1}\left(\left\{\rv\right\}\right)\right)\right) .$$
% (e.g. $\ignoreset_{\rv,\Rv,\Comp\Rv}=\{\omega\in\Omega\mid \exists \omega'\in\Omega\text{ such that }\Rv(\omega')=\rv\text{ and }\Comp\Rv(\omega)=\Comp\Rv(\omega')\}=\{\Comp\Rv\in\Comp\Rv(\{V=v\})\}$) 
% and the function
%  $$\ignore_{\Comp\rv,\Comp\Rv,\Comp\Rv}:.$$
 
If $(\Rv,\Comp\Rv)\isoeq\Omega$, then $\forall \Comp\rv\in\Comp\Rv(\Omega)$,  $\ignoreset_{\Comp\rv,\Comp\Rv,\Rv}=\Omega$, and 
$:\omega\mapsto \left(\interclass_{\Comp\Rv,\Rv}\right)(\Comp\rv,\Rv(\omega))$ is measurable $:(\Omega,\tribu_\Omega)\to(\Omega,\tribu_\Omega)$.

%For $A\in \tribu_\Omega$, such that $P(A)>0$, $P\restrict_A$ is the measure $\tribu_\Omega\to [0,1], B\mapsto P\large|_A(B)=P(A)^{-1} P(A\cap B)$.

If $\forall P\in\range{P}$, $\forall \Comp\rv\in\Comp\Rv(\Omega)$,  $P\left(\ignoreset_{\Comp\rv,\Comp\Rv,\Rv}\right)>0$, define:

 $$\range{P}^\star=\left\{\left.
 \left(\left<\left.\Comp\rv\mapsto P^{(\Rv,\Comp\rv)\mid \ignoreset_{\Comp\rv,\Comp\Rv,\Rv}}\right|P'\right>_{(\Rv,\Comp\Rv)(\Omega),\Comp\Rv(\Omega)}\right)^{\interclass_{\Rv,\Comp\Rv}}\right|(P,P')\in\range{P}\times\range{P}'\right\}.$$

%and
% $$\mathrm{Ignore}^{\Obs=\obs}_{\Comp\Rv,\Rv}:\range{P}\to,P\mapsto 
% \left\{\left(P^{\Rv\mid \ignoreset_{\Comp\rv,\Comp\Rv,\Rv}}\otimes \Dirac_{\{\Comp\rv\}}\right)^{\interclass_{\Rv,\Comp\Rv}}\right\}_{\Comp\rv\in \Comp\Rv(\Obs^{-1}(\obs))},$$

The model obtained after ignoring $\Comp\Rv$ is the model $\left(\codomain(\Obs),\tribu_{\codomain(\Obs)},\{(P^\star)^\Obs\mid P^\star\in\range{P}^\star\}\right)$.
\end{definition}

\begin{example}\label{ex:rubin1976}
When $\range{P}'=\left\{\Dirac_{\Comp\rv}\mid\Comp\rv\in\Comp\Rv(\Omega)\right\}$, this definition is the generalisation of the concept of "ignoring the process that causes missing data" from \cite[Sec.~6, p.~585]{Rubin1976}. Indeed, Rubin explains that when the nuisance process, equivalent to the missing data indicator, it can be considered as a fixed variable.
\end{example}

\begin{example}\label{ex:ignoringthedesignandthepopulation}
Fixed population model for design based inference is the model obtained after ignoring the population generation and the design variable generation ($\Comp\Rv=(\Signal,\Designvar)$) and chosing $\range{P'}=\{\Dirac_{}\mid (\signal,\designvar)\in(\Signal,\Designvar)(\Omega)\}$.
\end{example}

%\begin{remark}
Ignoring a process $\Comp\Rv$ is in general different from marginilising on $\Comp\Rv$ (which consists in only considering $P^{\Rv}$) or conditioning on $\Comp\Rv$ (which consists in only considering $P^{\Rv\mid\Comp\Rv}$, ignoring a process is about changing the model to suppress stochastic relations between $\Rv$ and $\Comp\Rv$ as much as possible.
When conditioning, one assumes that $\Rv$ is observed.
\cite{Rubin1976} gives conditions on the missing data process so that it can be ignored, e.g. conditions under which conditioning and ignoring lead to the same inference.
%\end{remark}

There is another way to ignore a nuisance process, without having to consider it as a parameter, which consists in considering that $\forall\rv$, 
$\rv\mapsto P^{\Comp\Rv\mid\Rv=\rv}=P^{\Comp\Rv\mid\ignoreset_{\rv,\Rv,\Comp\Rv}}$. This last approach is equivalent to chosing $\range{P}'=\{P^{\Comp\Rv}\mid P\in\range{P}\}$, and if in addition, $(\Rv,\Comp\Rv)\isoeq\Omega)$, this is equivalent of chosing 
$\range{P}^\star=\left\{\left(P^\Rv\otimes (P')^{\Comp\Rv}\right)^{\interclass_{\Rv,\Comp\Rv}}\mid P,P'\in\range{P}\right\}$.

\begin{property}
If $(\Rv,\Comp\Rv)\isoeq\Omega$, and $\range{P}'=\{P^\Comp\Rv\mid P\in\range{P}\}$, then
 $\range{P}^\star=\left\{\left(P^\Rv\otimes P^{\Comp\Rv}\right)^{\interclass_{\Rv,\Comp\Rv}}\setsep P\in\range{P}\right\}$.
\end{property}

\subsubsection{Transformation of the target of the inference}
We have seen that ignoring a random process consists in changing the family of probability distributions of the model. The object of the inference can be a function defined on the initial model. To completely describe what ignoring the model is, one needs to specify what is the new object of inference after ignoring the model.
The target of the inference is any function $\paramf:(\range{P}\times\Omega)\to$. For example, when $\paramf\dependsonly((P,\omega)\mapsto P)$, $\paramf$ is a parameter and the object of the inference is estimation, whereas when $\paramf\dependsonly((P,\omega)\mapsto \Signal(\omega))$, $\paramf$ is a function of $\Signal$ to predict and the object of the inference is prediction.
Given a target $\paramf$,  what is the corresponding target $\paramf^\star$ in the model $\range{P}^\star$ ?
There is a  natural answer in the following cases.
\begin{enumerate}
\item  If there exists  $\Signal:(\Omega,\tribu_\Omega)\to$ such that $\paramf:(\Omega\times\range{P})\to,(\omega,P)\mapsto \Signal(\omega)$, a natural choice is $\paramf^\star:\range{P}^\star\times\Omega\to,(P^\star,\omega)\mapsto (P^\star,\omega)\mapsto \Signal(\omega)$.
\item  If there exists  $\Signal:(\Omega,\tribu_\Omega)\to$ such that $\paramf\dependsonly(P,\omega\mapsto P^\Signal)$: $\paramf=\paramf'\circ(P\mapsto P^\Signal)$ and if $\{(P^\star)^\Signal\mid P^\star\in\range{P}^\star\}\subset\codomain(\paramf')$,  then one may define: $\paramf^\star:P^\star\mapsto\paramf'((P^\star)^\Signal)$.
\item When $\range{P}^\star\subset\range{P}$, a natural choice is $\paramf^\star=\paramf\restrict_{\range{P}^\star\times\Omega}$.
\end{enumerate}

%\cite{Rubin1976} is interested in , $\paramf:\P\mapsto P^\Signal$, so one naturaly defines $\paramf^\star:P^\star\mapsto (P^\star)^\Signal$.

\subsubsection{Transformation of the set of prior distributions}
Switching from $\range{P}$ to $\range{P}^\star$ after ignoring one process, makes it difficult to define a new set of prior distributions $\range{\priordist}^\star$ because $\priordist$ and $\priordist^\star$ are not distributions defined on the same set. When each prior distribution is the product of a prior on the process of interest and a prior on the nuisance process, a natural choice for the new set of prior distributions on the new model can be made.

The natural way to define a new set of prior distributions is to proceed this way:
\begin{enumerate*}
\item Define a set $\tilde{\range{\priordist}}$  of prior (not necessarily proper) distributions on $\range{P}'$.
\item Define $\range{\priordist}^\star=\left\{\left(\priordist^{P\mapsto \left(\Comp\rv\mapsto P^{(\Rv,\Comp\rv)\mid \ignoreset_{\Comp\rv,\Comp\Rv,\Rv}}\right)}  \otimes \priordist' \right)^{(P'',P')\mapsto\left(\left(\left<\left.P''\right|P'\right>_{(\Rv,\Comp\Rv)(\Omega),\Comp\Rv(\Omega)}\right)^{\interclass_{\Rv,\Comp\Rv}}\right)}\setsep \priordist\in\range{\priordist},\priordist'\in\range{\priordist}'\right\}$.

\end{enumerate*}
\begin{example}
When $\range{P}'=\{\Dirac_{\Comp\rv}\setsep \Comp\rv\in \Comp\Rv(\Omega)\}$, define $\left(\int P^{\Comp\Rv}\derive \priordist(P)\right)^\Dirac$ as the measure on $\range{P}'$: $\left(\int P\derive \priordist(P)\right)=\left(\Comp\Rv(\tribu_\Omega)\to,\somesubset\mapsto\int_\range{P} P^{\Comp\Rv}(A)\derive Q(P)\right)^{\Comp\rv\mapsto\Dirac_{\Comp\rv}}$. A possible choice for the set of prior distbution on $\range{P}'$ is 

$$\range{\priordist}'=\left\{\priordist''\setsep \priordist'<<\left(\int P^{\Comp\Rv}\derive \priordist(P)\right)^\Dirac \text{ and }\priordist'>>\left(\int P^{\Comp\Rv}\derive \priordist(P)\right)^\Dirac ,\priordist\in \range{\priordist} \right\}.$$
\end{example}

\subsubsection{Transformation of the dominant measure}

The model $\range{P}$ may be dominated by a $\sigma$-finite measure $\dominant$ used to define the likelihood for example. When $\range{P}'=\{\Dirac_{\Comp\rv}\mid\Comp\rv\in\Comp\Rv(\Omega)\}$, the model $\{P^\Obs\mid P\in \range{P}^\star\}$ is dominated by the measure $\left(\dominant^\Rv\otimes\countingmeasure_{\Comp\Rv(\Omega)}\right)^{\Obs[\interclass_{\Rv,\Comp\Rv}]}$, where $\countingmeasure_{\Comp\Rv(\Omega)}$ is the counting measure on $\Comp\Rv(\Omega)$. This dominant measure is not necessarily $\sigma$-finite. %The choice of $\range{P}'$% unless ${\Comp\Rv(\Omega)}$ is countable. 

\subsection{Equivalent inferences}

Comparing the inference before and after ignoring the model is relevant when the two inferences are of the same nature.
%For example, in the $\Obs=\obsf(\Signal,\Transfo[\Signal],\Designvar)$ model, it does not seem relevant to compare the properties of make sense to compare the prediction of a function of $\Signal$ in a frequentist framework with the prediction of the same quantity after ignoring the selection $\Transfo$ in a frequentist framework.
The only thing that the model after ignoring a selection mechanism and the model before have in common is the probability space $(\Omega,\tribu_\Omega)$ and the random variables defined on it.
No matter what type of inference (Bayesian inference, maximum likelihood estimation in a frequentist framework, testing, model choice), inference always comes down to studying the properties of some random quantities for specific sets of probability distributions and prior probability distributions. 

\subsubsection{Likelihood based inference.}

The likelihood based inference on $\paramf:\range{P}\to$ is equivalent to the likelihood based inference on $\paramf^\star:\range{P}^\star\to$ when $\Obs=\obs$ if $\codomain(\paramf)=\codomain(\paramf^\star)$ and $\exists~ \alpha\in\mathbb{R}^+$ such that
$\left(\param\mapsto\likelihood^\star_{\paramf^\star;\Obs}(\param;\obs)\right)=\left(\param\mapsto\alpha~\likelihood_{\paramf,\Obs}\left(\param;\obs\right)\right)$.

\subsubsection{Estimation in the frequentist framework}

The inference based on the estimation of $\paramf:\range{P}\to$ by $\hat{\paramf}:\Omega\to\codomain(\paramf)$ is equivalent to the estimation of $\paramf^\star:\range{P}\to$ by $\hat{\paramf}$ if $\codomain{\paramf}=\codomain{\paramf}^\star$ and 
$\forall\param\in\codomain(\paramf)$, $\{P^{\hat\paramf}\mid P\in\paramf^{-1}(\param)\}=\{(P^\star)^{\hat\paramf}\mid (P^\star)\in(\paramf^\star)^{-1}(\param)\}$.

\subsubsection{Bayesian inference.}

The Bayesian inference on $\paramf:(\range{P}\times\Omega)\to$ with the set of prior distribution $\range{\priordist}$ is equivalent to the likelihood based inference on $\theta^\star:(\range{P}^\star\times\Omega)\to$  with the set of prior distribution $\range{\priordist}$ when $\Obs=\obs$ if 
$\{(\priordist^\star)^{\paramf^\star\mid\Obs=\obs}\mid\priordist^\star\in\range{\priordist}^\star\}=
\{\priordist^{\paramf\mid\Obs=\obs}\mid\priordist\in\range{\priordist}\}$

\section{Ignorable vs. Informative Process: Definition and Characterisation}

\subsection{Definition}

The process $\Comp\Rv$ will be called ignorable for a specific type of inference if for $\range{P}'=\left\{\Dirac_{\Comp\rv}\setsep \Comp\rv\in\Comp\Rv(\Omega))\right\}$, and $\range{\priordist}'=\left\{\priordist''\setsep \priordist'<<\left(\int P^{\Comp\Rv}\derive \priordist(P)\right)^\Dirac \text{ and }\priordist'>>\left(\int P^{\Comp\Rv}\derive \priordist(P)\right)^\Dirac ,\priordist\in \range{\priordist} \right\}$ for which this type of inference based on $(\range{P},\paramf,\likelihood,\range{\priordist},...)$ is equivalent to the same type of inference based on 
$(\range{P}^\star,\paramf^\star,\likelihood^\star,\range{\priordist}^\star,...)$.
In the case they are not equivalent, the process $\Comp\Rv$ will be called informative.

\begin{examples}
\begin{enumerate*}
\item In the $\Obs=\obsf(\Transfo[\Signal],\Transfo,\Designvar)$ Assume that $\Obs=\Transfo[\Signal]$, that $\Signal$ follows a population model, and that $\Transfo$ is a selection without replacement. 
The population model is that $(Y[k])_{k\in\Pop}$ are iid realisations of some distribution, 
and the goal is to estimate $\paramf=E[Y[1]]$. If $\forall P\in\range{P}$, 
$P^{\Designvar,\Transfo,\Signal}=P^{\Designvar,\Transfo}\otimes P^{\Signal}$, then 
the selection and the design variables are not informative.
\item If we assume that  selection can be with replacement, then 
the selection and the design variables are informative.
\end{enumerate*}
\end{examples}

\section{Discussion}
This section is a draft that gives an idea of the topics we would like to discuss in the final version.

\subsection{Ignorability and Ancillarity, Informativeness and information}

Consider the following example:
Assume $Y\mid n\sim\mathrm{Bernoulli}(\param,1)^{\otimes n}$,
$(n-1)\sim\mathrm{Binomial}(\theta,N)$. 

The inference consists in estimating $\param$. The theory tells us that the best estimator is 
a affine combinaison of $(n-1)/N$ and the  average of $Y$. And this statistic contains all the information about $\param$.
Should $n$ be ignored ? no.
Can $n$ be ignored ? yes if the statistician has chosen to use the average of $Y$ only. The average of $Y$ has the same distribution in the model where $n$ is not ignored, and in the model where $n$ is ignored.
Is $n$ informative ? yes, if we consider that informative means: contains additional information on the parameter of interest, with respect to the information contained in $Y/n$. No if we consider that informative means: contains information that will condition the distribution of $Y/n$.

The choice of the term informative selection seems to have appeared for the first time in \cite{Scott1975}, it was a reference to design for which any prior would have no effect on the posterior distribution of the parameters of distribution of the process of interest: it was clearly stated for a joint model on the nuisance and interest processes was separated in to a nuisance process model and a study process model, which excludes the case considered above. 

At first sight, it seems that an easy answer would be to consider a sufficient statistic in the original model with respect to $\param$, and require that the distribution of the sufficient statistic is the same in the orignal model and in the ignored model. This raises the problem of defining sufficient and ancillary statistics in presence of nuisance parameters. So it seemed wise to separate the question about what information could be used 
from the question about doing a correct inference. Another aspect about sufficient statistic and ancillary statistic is that those notions only apply to non latent processes, whereas  the nuisance process is also latent. 

%In the case where the nuisance and interest processes are separated, the likelihood of the observations can be expressed 
%
%$\mathscr{L}(

\subsection{Sufficient conditions for ignorable transformation}

In this section we will illustrate the sufficient conditions for ignorable transformation on real life examples.
As mentioned, independence of the transformation process and the process of interest are not necessarily sufficient to ensure ignorability.

\subsection{Testing ignorability vs informativeness of a process}

In this subsection we propose to give sufficient conditions for informativeness.
%In certain cases, ignorability of $\Comp\Rv$ can be seen as a property of $P$. In this case, one can define testing ignorability. 
%Sufficient conditions for ignorability $\{\range{P}^\star\subset\range{P}$ are
%\begin{enumerate*}
%\item $\range{P}'=\{P^{\Comp\Rv}\setsep P\in\range{P}\}$, and 
%\item $\left(P\mapsto\left(\Comp\rv\mapsto P^{(\Rv,\Comp\rv)\mid \ignoreset_{\Comp\rv,\Comp\Rv,\Rv}}\right),P\mapsto P^{\Comp\Rv}\right)\isoeq\range{P}$.
%\item $\forall P\in\range{P}$, $P=\left<\left.\Comp\rv\mapsto P^{(\Rv,\Comp\rv)\mid \ignoreset_{\Comp\rv,\Comp\Rv,\Rv}}\right|P'\right>_{(\Rv,\Comp\Rv)(\Omega),\Comp\Rv(\Omega)}\right)^{\interclass_{\Rv,\Comp\Rv}}$
%\end{enumerate*}
%Under these conditions, testing that $\Rv$ is ignorable consists in testing that:
%$$.
%For example, assume an experiment with missing data, with a logistic model linking the probability for a missing sampled unit to the value of the study variable for the same unit.
%When chosing $\range{P}'=\{P^\Comp\Rv\mid P\in\range{P}\}$, the transformed model $\range{P}^\star$ corresponds to a model under which $\Rv=\Signal$ and $\Comp\Rv=\Transfo$ are independent. If, for $P\in\range{P}$, the regression coefficient associated to the study variable  is null, it corresponds to the case where $P\in\range{P}^\star$. Testing ignorability of selection is equivalent to testing the nullity of this regression coefficient.
 
\subsection{On the question of weighting or not weighting}

The question of using the samplign weights has been treated by numerous authors. 
Without narrowing too much the debate, we consider the case where $\Transfo$ is a selection and where 
$\Obs=\obsf(\Transfo[\Signal],\Transfo,Z)=(\Transfo[\Signal],\Transfo[\inclusionprob],\func(Z))$: we observe the values of $\Signal$ on the sample, the values of the inclusion probabilities on the sample, as well as some information derived from the design variable (which can be the strata for the whole population, the list of units in the population, ...etc).

There are different ways of using or not using the weights and different environments: 
\begin{enumerate*}
\item ignore the weights, and use a model where selection and weights are non informative, \item ignore the weights, because auxiliary information is available (for example the design and design variables on the population) and contains already the information brought by the weights, 
\item use the weights and $\Transfo[\Signal]$ and discard the information brought by all other auxiliary information.
\end{enumerate*}.

We can remark that the questions of using or not using the weights is linked to the two distinct following questions: Can we ignore the selection ? Do the weights bear any information  that should be used ?

\bibliographystyle{authordate1}
%\bibliography{Art_Def_Informative}

\appendix
\newpage
{\huge Supplemented Materials}

We took some notes on different papers that deal with sufficiency, partial sufficiency, 
ignorability, and informative selection, that we reproduce below. These notes can be useful for the reader, but are not necessary for the comprehension of the paper.

\section{Remarks and notes on set theory}

\subsection{On the ambiguity of ``dependence''.}

The statement ``$\func:(\elt_1,\elt_2,\elt_3,\elt_4)\mapsto (\elt_1+\elt_3)$ depends only of $\func':(\elt_1,\elt_2,\elt_3,\elt_4)\to (\elt_1,\elt_3)$'' means that there exists a function $\funcc$ such that $\func=\funcc\circ\func'$. The statement  ``$\func:(\elt_1,\elt_2,\elt_3)\mapsto (\elt_1+\elt_3)$ does depends on $\func':(\elt_1,\elt_2,\elt_3,\elt_4)\to (\elt_3)$'' means that the value of $\func$ is not necessarily the same for all $\elt_3$, all other things being equal, e.g. $(\elt_1,\elt_2,\elt_4)=(\elt'_1,\elt'_2,\elt'_4)\nLeftrightarrow\func(\elt)=\func(\elt')$. The statement ``$\func:(\elt_1,\elt_2,\elt_3)\mapsto (\elt_1+\elt_3)$ does not depend on $\func':(\elt_1,\elt_2,\elt_3,\elt_4)\to (\elt_2)$'' means that the value of $\func$ is the same for all $\elt_2$, all other things being equal, e.g. $(\elt_1,\elt_3,\elt_4)=(\elt'_1,\elt'_3,\elt'_4)\Leftrightarrow\func(\elt)=\func(\elt')$. There is an ambiguity in the definition of "depends on" that comes from the difficulty to give a mathematical meaning to the expression "all other things being equal".
Consider the statement  $\func:(\elt_1,\elt_2,\elt_3)\mapsto (\elt_1+\elt_3)$ does not depend on $(\elt_1,\elt_2,\elt_3)\to \elt_2$. This statement is ambiguous because 
the codomain of the two functions are not specified. If the codomain is $\mathbb{R}^3$, the statement is not ambiguous: the function $\func$ depends only on $\func':(\elt_1,\elt_2,\elt_3)\mapsto (\elt_1,\elt_3)$, which is a variation independent complement of $\func$. If the domain is 
$\spaceE:\{(\elt_1,\elt_2,\elt_3)\in\mathbb{R}^3\mid \elt_1=\elt_2$, then in this case, $\func:\spaceE\to,(\elt_1,\elt_2,\elt_3)\mapsto (\elt_1+\elt_3)=:\spaceE\to,(\elt_1,\elt_2,\elt_3)\mapsto (\elt_2+\elt_3)$. 
So even if there is a formal independence in $(\elt_1+\elt_3)$, when applied to functions, it is ambiguous, as the domain of $\func$ may be the outcome of some constraining that link $x_1$ and $x_2$.
The conclusion is that independence is with respect to a complete reparametrisation of the parameter, $\func$ does not depend of $\func'$ has a meaning when it is understood as $\func$ depends only of $\Comp\func'$ where $\Comp\func'$ is a complement of $\func'$ that has been specified. In the case where $\spaceE:\{(\elt_1,\elt_2,\elt_3)\in\mathbb{R}^3\mid \elt_1=\elt_2$, $Id_\spaceE$, $\elt\mapsto (\elt_1,\elt_3)$ and  $\elt\mapsto (\elt_3)$ are in 
As we have seen, there always exist a complement, although it is not a distinct complement. When we say  $\func:(\elt_1,\elt_2,\elt_3)\mapsto (\elt_1+\elt_3)$ does not depend on $(\elt_1,\elt_2,\elt_3)\to x_2$ all other parameters being equal, we implicitely fix $\elt_1$ and $\elt_3$ and check that for every possible $\elt_2,\elt_2'$ such that $\{(\elt_1,\elt_2,\elt_3),(\elt_1,\elt_2',\elt_3)\}\in\spaceE$, $\func(\elt_1,\elt_2,\elt_3)=\func(\elt_1,\elt_2',\elt_3)$. But it should be avoided and the use of ``depends only'' should be preferred, the natural definition of $\func$ is independent of $\func'$ when $\Comp\func$ is fixed is the following.

\begin{definition}
  $\func:\spaceE\to$ is deterministically independent of $\func':\spaceE\to$ when $\Comp\func'$ is fixed if and only if any of the following statement is true:
  \begin{enumerate*}
  \item $\forall \elt,\elt'\in\spaceE, \Comp\func'(\elt)=\Comp\func'(\elt')\Rightarrow \func(\elt)=\func(\elt')$,
  \item $\func$ depends only on $\Comp\func'$,
  \item $\exists \funcc:\to$ such that $\func=\funcc\circ\Comp\func'$.
  \end{enumerate*}
\end{definition}

When we deal with nuisance parameters, nothing prevents to link the nuisance process to the process of interest by saying that the same parameter governs them, as for example in cutoff sampling.
The parameter space is not the product of the two parameters spaces and still we may want to be able to apply the definition of informative or non informative process to this case. It is then necessary to be clear about this statement ``the likelihood does not depend on the nuisance parameter".

\begin{property}[Variation independence]
Given two functions $\func,\func':\spaceE\to$, the following statements are equivalent:
\begin{enumerate*}
\item $\func$ and $\func'$ are variation independent,
\item $\image(\func,\func')=\image(\func)\times\image(\func')$.
\item $\forall \eltt\in\image(\func), \eltt'\in\image(\func')$, $\exists \elt\in\spaceE$ such that $\func(\elt)=\eltt$ and $\func'(\elt)=\eltt'$.
\item $\forall \eltt\in\image(\func')$, $\image\left(\func\restrict_{\func'^{-1}(\{\eltt\})}\right)=\image(\func)$.
\end{enumerate*}
\end{property}

%We may also use the following notations for complements:
% Let $\func: \allfunc{\spaceE}{}$, and let $\Comp\func:\allfunc{\spaceE}{}$ be a complement of $\func$.
% If $(\func,\Comp\func)\isogeq\spaceE$, we will note $\spaceE\subseteq(\spaceE/\func)\rtimes(\spaceE/{\Comp{\func}})$, if $(\func,\Comp\func)\isoeq\spaceE$, we will note $\spaceE=(\spaceE/\func)\rtimes(\spaceE/{\Comp{\func}})$  if not, and if $(\func,\Comp\func)\isog\spaceE$, we will write $\spaceE\subsetneq(\spaceE/\func)\rtimes(\spaceE/{\Comp{\func}})$.

\section{Notes on selected work on the definition of information, sufficiency, in particular in presence of nuisance parameters}\label{appendix:nuisanceparameters}

%To generalise definitions of informative or ignorable processes, we use 

\subsection{Conditional distributions}

%In this section, we 

\begin{propertydefinition}[Transition probabilities](see \citet[p.~310]{Sverdrup1966Present})
Consider a measured space $(\Omega,\tribu,P)$, and a subsigmafield $\tribu'$ of $\tribu$, then define the conditional probability of $A$ given ''the most accurate description of $y$ by means of statements from $\tribu'$ as the value un $(A,\omega)$ of the function: $\tribu\times \Omega\to[0.1], (A,\omega)\mapsto P(A\mid \tribu_Y, \omega)$, which $P(A\mid \tribu_Y, .):\Omega\to[0,1],\omega\mapsto P(A\mid \tribu_Y, \omega)$ is the class of almost unique $\tribu'$-measurable functions of $\omega$ which satisfies $[\forall B\in \tribu'],\ P(A\cap B)=\int_B P(A\mid \tribu',\omega) \derive P(\omega)$. The function $P(A\mid\tribu',.)$ does always exist and is almost uniquely defined.
\end{propertydefinition}

\subsubsection{Conditions for existence of conditional distributions.}

One could have assumed that the measurable space $(\Omega,\tribu_\Omega)$ satisfies the regular conditional probability property (this assumed property conveniently states that all conditional probabilities are defined). As \cite{DieudonnesurLebesgueNikodym} mentions, it is only available for certain distributions. So this assumption is not convenient there. For more reading on regular conditional probability, one can look at \url{https://en.wikipedia.org/wiki/Regular_conditional_probability}, but for academic references, some are listed here: \url{https://scielo.conicyt.cl/pdf/proy/v23n1/art02.pdf}, which includes Kolmogorov (1933, in german) and \cite{DieudonnesurLebesgueNikodym}, the paper is in french but the notion is not the main topic, it is lost in chapter 7. So it seemed reasonable here to just say: Assume that all given conditional distributions are defined, without adding unnecessary refinment.

\subsection{Densities}

\fancyquote{\citet[p.~229, ]{HalmosSavage1949}}{
{\bf 3. Conditional probabilities and expectations.} \textsc{Lemma 4.} {\it If $\mu$ and $\nu$ are measures on $\mathbf{S}$ such that $\nu\dominatedby\mu$, then $\nu T^{-1}\dominatedby\mu T^{-1}$}.

\textsc{Proof.} If $F\in T$ and $0=\mu T^{-1}(F)=\mu(T^{-1}(F))$, then 
$$0=\nu(T^{-1}(F)=\nu T^{-1}(F).$$

Lemma 4. is the basis of the definition of a concept of great importance in probability theory. If $\mu$ is a measure on $\mathbf{S}$ and $f$ is a non negative integrable function on $X$, then the measure $\nu$ defined by $d\nu=fd\mu$ is absolutely continuous with respect to $\mu$. It follows from Lemma 4 that $\nu T^{-1}$ is absolutely continuous with respect to $\mu T^{-1}$; we write $d\nu T^{-1}=gd\mu T^{-1} $. The function value $g(y)$ is known as the {\it conditional expectation} of $f$ given $y$ (or given $T(x)=y$)
}

Let $(\Rv,\Comp\Rv)$ be a separation of $\Omega$,
then $(P^\Rv\otimes P^{\Comp\Rv})^{\interclass_{\Rv,\Comp\Rv}}$ is a measure on $\Omega$ that may differ from
$P=\left(P^{(\Rv,\Comp\Rv)}\right)^{\interclass_{\Rv,\Comp\Rv}}$, even when $(\Rv,\Comp\Rv)$ is a distinct separation of $\Omega$. 
For example, if $P=\mathrm{Normal}\left(0,\begin{pmatrix}1&\rho\\\rho&1\end{pmatrix}\right)$, with $\rho\neq 0$, $\Rv=X_1:(x_1,x_2)\mapsto x_1$, $\Comp\Rv=X_2:(x_1,x_2)\mapsto x_1$, $\interclass=\mathrm{Id}_{\mathbb{R}^2}$.
%Assume now $P'=P^\X_1\otimes \delta_{0}$

Let $f$ be a density of $P$ with respect to $\dominant$ and $f_\Rv$ a density of 
$P^\Rv$ with respect to $\dominant^\Rv$.
If $P\dominatedby \dominant$, then $P^\Rv\dominatedby \dominant^\Rv$ (\citep[See][p.~229, Sec.~3]{HalmosSavage1949}). The relationship between $\derive P/\derive\dominant$ and $\derive P^\Rv/\derive\dominant^\Rv$ is as follows:
$\forall \func$ measurable 
$\int_{\Rv(\Omega)}            \func         ~\derive P^\Rv
=\int_{\Omega}                 \func\circ\Rv ~\derive P  
=\int_{\Omega}                 \func\circ\Rv ~\derive (f.\dominant)
=\int_{\Omega}         f      .\func\circ\Rv ~\derive \dominant
=\int_{\Rv(\Omega)}    f_{\Rv}.\func         ~\derive \dominant^\Rv$.
In the case where $\Rv$ is a $C^{(1)}$ diffeomorphism and that the measurable spaces involved are subsets of $\mathbb{R}^n$. The change of variable formula gives a relationship between $f$ and $f_\Rv$.

% \subsubsection{Sufficiency, ancillarity, ignorability}
% \begin{definition}[Sufficient statistic]
% A sufficient statistic (\cite{Sverdrup1966Present}[p.~])\todomargin{Give page} is a random variable $\Sufstat:(\Omega,\tribu_\Omega)\to$ such that 
% \begin{equation}\label{eq:sufficiency}
% \forall A\in\tribu_\Omega,\ \exists\Phi_A:(\Omega,\tribu(S))\to([0,1],\Borel_{[0,1]})\text{ such that }\forall P\in\range{P},\ P\text{a.s.}(\omega),\ \Phi_A(\obs)=P(A\mid \tribu(S),\omega).
% \end{equation}
% \end{definition}

%
%
%We propose to extend the definition of \cite{Fraser1956} for partially sufficient statistic with respect to $\thetaf$. Generanlisation comes from 
%\begin{definition}[$\Comp\paramf$-sufficient statistic with respect to $\paramf$]
%Let $\paramf:\{P^\Obs\mid P\in \range{P}\}\to$ and let $\Comp\paramf$ be a complement of $\paramf$, then the statistic $\Sufstat:\allfunc{(\Omega,\tribu_\Omega)}{}$ is a $\Comp\paramf$-partially sufficient statistic with respect to $\paramf$  if and only if 
%\begin{multline*}\label{eq:partialsufficiency}
%(:\allfunc{\range{P}}{},P\mapsto P^\Sufstat)\perp \Comp\paramf\text{ and }\\
%\forall A\in\tribu_{\Omega}, \forall \range{P}\in\left(\range{P}/\Comp\paramf\right) \exists\Phi_A:(\Omega,\tribu(\Sufstat))\to([0,1],\Borel_{[0,1]})\text{ such that }\\
%\forall P\in\range{P},\ P\text{-a.s.}(\omega),\ \Phi_A(\omega)=P(A\mid \tribu(\Sufstat),\omega)
%\end{multline*} 
%\end{definition}

\subsection{Sufficient statistic}
For sufficient statistics with respect to $\theta$, without nuisance parameters details can be found in the following books: \cite{LehmannCasellaTheoryofPointEstimation},\cite{LoeveProbabilityTheory} and \cite{Sverdrup1967Laws}.
Three things are important to be able to define sufficient statistics: 1. The definition consist on a property satisfied by some transition probabilities, which means that they must be defined, and there exists a list of conditions for transition probabilities to be defined and to be probabilities. 2. The definition of a sufficient subsigmafield. 3. The definition of a sufficient statistic.

\begin{enumerate}
 \item Definition of transition probabilities for a probability space $(\Omega,\tribu_\Omega,P)$ wrt a subsigmafield $\mathscr{B}$ of $\tribu_\Omega$.
 \citet[Vol.~2, Sec.~ D.5, p.~306]{Sverdrup1967Laws}, \citet{LehmannCasellaTheoryofPointEstimation}.%, \cite{Rudin}. 
 Below are some essential elements about transition probabilities
    \begin{enumerate}
     \item Existence of transition probabilities is a consequence of the Radon Nikodym theorem.
     \item The result can be extended to measures.
     \item Assumptions may be added to get 3 nice conditions satisfied by the transition probabilities, one needs additional assumption: the spaces must be euclidean. This is explained in \citet[Vol.~2, Sec.~ D.5, p.~306]{Sverdrup1967Laws}, \citet[]{LehmannCasellaTheoryofPointEstimation}. The three nice properties are as follows:
     \item transition probabilities are ${B}$ measurables (into $[0,1],\mathscr{B}_{[0,1]}$, defined by  $$\forall A\in\tribu_\Omega,\ P(A\mid \mathscr{B},.):\omega\mapsto P(A\mid\mathscr{B},\omega),$$
     such that $\forall B\in\mathscr{B}$, $$P(B\cap A)=\int_B P(A\mid \mathscr{B},\omega)\d P(\omega).$$
     \item transition probabilities are all defined $P$-a.s. and ${B}$ measurables. which means that if one takes two representants of the class, the difference is also $A_0$ measurable. The set where they differ is the reciprocal image of $0$ by the difference so it is an element of $\mathscr{B}$ of measure $0$.
     \item when working on euclidean spaces with sigmafields being Borel fields, with densities, and when $\mathscr{B}=\tribu(Y)$, conditional density can be derived.
     \item Although it has no meaning at this stage(when conditional distributions have not been defined, Define the equivalence relation $\omega\sim_\mathscr{B} \omega'\leftrightarrow[\forall B\in \mathscr{B},[t\in B\leftrightarrow t'\in B]]$. Then 
     $$P(A\mid\mathscr{B},\omega)=\int_{\class_{\sim_\mathscr{B}}(\omega)} \setindicator_A \d P^{\Id\mid \class_{\sim_\mathscr{B}}=\class_{\sim_\mathscr{B}}(\omega)}.$$
     One can see that $P(A\mid\mathscr{B},.)$ is necessarily constant on each $\class_{\sim_\mathscr{B}}(\omega)$.
     The problem is that the measures $P^{\Id\mid \class_{\sim_\mathscr{B}}=\class_{\sim_\mathscr{B}}(\omega)}$ are not necessarily defined unless $P(\class_{\sim_\mathscr{B}}(\omega))>0$, and at this stage, have no meaning as conditional probabilities were not defined.
     \item once the existence of the transition probabilities up P-a.s. equivalence is proven, the conditional expectation of a random variable X is defined as a $\mathscr{B}$ measurable random variable, defined up to a equivalence relationship, such that $\forall f$, integrable[...] $E[f(X)]=\int f(E[X\mid\mathscr{B}]) \d P$.
     (although one understands it better when written like this:)
     $$E[X\mid B]=\class\left(:\Omega\to,\omega\mapsto\int_{\class_{\sim_\mathscr{B}}(t)} X \d P^{\Id\mid \class_{\sim_\mathscr{B}}=\class_{\sim_\mathscr{B}}(\omega)}\right)$$
    \end{enumerate}
 
%  \item definition of conditional probability for a probability space $(\Omega,\tribu_\Omega,P)$ wrt ${B}$ a subsigmafield of $\tribu_\Omega$. 
%  \begin{itemize}
%   \item Here again, no need for density, although it is easier to define with a joint density.
%  \end{itemize}
 
 \item Definition of a sufficient subsigmafield for a set of probability distributions.
   \begin{itemize}
    \item Definition of sufficiency does not require domination of the model. 
    \item Many books use a parametric dominated model. \citet[p.~344]{LoeveProbabilityTheory} does not use a parametric model, but a dominated one. He is defining a sufficient sigmafield wrt a set of measures (not necessarily probability measures.
    \item more generally, sufficient statistic consist in saying that forall $A\in \tribu_\Omega$, there exists $\Psi_A$  such that classes of transition probabilities  $(P(A\mid \mathscr{B},.)$ all contain $\Psi_A$,There is a difficulty there because the class of transition probabilities for $P,A$ is the class wrt $P-a.s$ equality, whereas  the class of transition probabilities for $P',A$ is the class wrt $P'-a.s$ equality.
    Let take a representent $\Psi_{A,P}\in (P(A\mid \mathscr{B},.)$ for each $P$. Then define $B_P=(\Psi_A-\Psi_{A,P})^{-1}(\{0\})$. So $B_P\in \mathscr{A_0}$ for each $P$, and for each $P$ $P(B_P)=0$.
    We do not necessarily have $P'(B_P)=0$. For example: $\Omega=\{0,1\}$, $\mathscr{B}=\mathscr{A}=\mathrm{Subsets}(\Omega)$, $\mathscr{P}=\{\delta_0,\delta_1\}$, $\Psi_A=\delta_A$. 
    \item Note that  \citet[p.~344]{LoeveProbabilityTheory} definition is the definition of a sufficient sub sigma field for a family of random variables defined on the same measured space.
For each of these random variables $Z_t:(\Omega,\tribu_\Omega,\mu)\to(\mathscr{R}^+,\mathbb{B}_{R^+})$?, Loeves denotes by $\mu_t$ the indefinite integrals of $X_t$.
(see p. 129 for defintiion of indefinite integrals),$\mu_t$ is the measure $:A\mapsto \int_AX_t\d\mu$. So $\mu_t$ in that case are all continuous wrt $\mu$, and $X_t$ is a density of $\mu_t$ wrt $\mu$. 
$X_t$ are supposed to be positive and integrable. 
    \end{itemize}

\end{enumerate}

\begin{definition}[Sufficient statistic]
A sufficient statistic (\cite{Sverdrup1966Present}[p.~310]) is a random variable $\Sufstat:(\Omega,\tribu_\Omega)\to$ such that 
\begin{equation}\label{eq:sufficiency}
\forall A\in\tribu_\Omega,\ \exists\Phi_A:(\Omega,\tribu(S))\to\left([0,1],\Borel_{[0,1]}\right)\text{ such that }\forall P\in\range{P},\ P\text{-a.s.}(\omega),\ \Phi_A(\omega)=P(A\mid \tribu(S),\omega).
\end{equation}
\end{definition}

\begin{remark}
 Condition \eqref{eq:sufficiency} is almost equivalent to 
\begin{equation*}
\forall A\in\tribu_{\range\Obs},\ \left(:\allfunc{\range{P}}{},P\mapsto P(A\mid \tribu(S),.)\right)\dependsonly (\mathds{1}_{\range{P}}).\label{eq:sufficiency2}
\end{equation*} 

It is almost equivalent and not equivalent because of the almost sure equality between $\Phi_A$ and the conditional distributions derived from the distributions of the model.
If we consider that $P(A\mid \tribu(\Sufstat),.)$ is one representative of the class and not the whole class, in a dominated context, the equivalence holds. In a non dominated context, the equivalence does not hold.

So Condition \eqref{eq:sufficiency} is weaker than Condition \eqref{eq:sufficiency2}.

\end{remark}

\subsection{Inference in presence of nuisance parameters. Specific and partial sufficiency.}

In the current section, we will give some quotes from some papers that deal with the definition of specific and  partial sufficiency or ancillarity in presence of nuisance parameters. We start with giving all the definitions in our notations.

\subsubsection{Definitions and comments in the current paper notations}
In our notation, the \citet[Eq. (1)]{Fraser1956} Definition is equivalent to the following:
\begin{definition}[Specific sufficiency]
Consider the statistical model $(\Omega,\tribu_\Omega,\range{P})$. Let $\paramf:\range{P}\to$, and let $\Comp\paramf$ be a complement of $\paramf$ . The random variable $\Sufstat:(\Omega,\tribu_\Omega)\to$. is a $(\paramf,\Comp\paramf)$-specific sufficient statistic if and only if:
\begin{eqnarray}
\label{eq:specificsuff1} \left(:\range{P}\to,P\mapsto P^\Sufstat \right) \dependsonly\paramf,\\
\label{eq:specificsuff2}\left(
            :\range{P}\to,P\mapsto \left(:\tribu_\Omega\times \Omega\to, (A,\omega)\mapsto P(A\mid\tribu(\Sufstat),\omega)\right)\right) \dependsonly\Comp\paramf,\\ 
\label{eq:specificsuff3}\text{ and }(\paramf,\Comp\paramf)\isoeq\mathscr{P}.
\end{eqnarray}
\end{definition}

There are three  conditions: (Eq. \eqref{eq:specificsuff1}) $P^\Sufstat$ depends only on $\paramf$,(Eq. \eqref{eq:specificsuff2}) the conditional distribution of $\Id_\Omega$ given $\Sufstat$ depends only on $\Comp\paramf$, and (Eq. \eqref{eq:specificsuff3}) $\paramf$ is a distinct complement of $\Comp\paramf)\isoeq\mathscr{P}$.
 We can see it is a generalisation of the definition of a sufficient statistic by choosing $\paramf=\Id_{\range{P}}$, $\Comp\paramf=\mathds{1}_{\range{P}}$. Heuristically them A $(\paramf,\Comp\paramf)$-specific sufficient statistic is a statistic that (Eq. \eqref{eq:specificsuff1}) contains no information about $\Comp\paramf(P)$ but the one already included in $\paramf(P)$ , such that (Eq. \eqref{eq:specificsuff2}) $P^{\Id_\Omega\mid\Sufstat}$ contains no information about $\paramf(P)$ but the one already included in $\Comp\paramf(P)$  , and such that (Eq. \eqref{eq:specificsuff3}) $\paramf$ contains no information about $\Comp\paramf$ and vice versa..
In the case where 
$\paramf=\Id_{\range{P}}$, $\Comp\paramf=\mathds{1}_{\range{P}}$, the two first and third conditions are always true. Other ways to generalize the definition of sufficient statistic would have been not to include first or third condition or both, which would have insured the existence of a sufficient statistic.
If we drop the condition that $(\paramf,\Comp\paramf)\isoeq\mathscr{P}$, the previous definition still makes sense, but it does not correspond to the definition as proposed by \cite{Fraser1956}. It is important to note that in \cite{Fraser1956} definition though, to define a $\theta$-sufficient statistic, one needs to specify $\eta$, the distinct complement of $\theta$. It is important to say that if $\paramf,\Comp\paramf\isoeq\mathscr{P}$, for any $\func:\Comp\paramf(\range{P})\to$ bijective, necessarily, $\paramf,\func\circ\Comp\paramf\isoeq\mathscr{P}$. So when $(\paramf,\Comp\paramf)\isoeq\mathscr{P}$, in that sense, specifying $\Comp\paramf$ does not matter.

Having to impose the condition that $(\paramf,\Comp\paramf)\isoeq\mathscr{P}$ has drawbacks: Consider for example the model $\left(\mathbb{N}^2,\Borel_{\mathbb{N}^2},\left\{\mathrm{Poisson}(\lambda_1)\otimes\mathrm{Poisson}(\lambda_2)\mid \lambda_1<\lambda_2\right\}\right)$. With $\Id_\Omega=(\Rv_1,\Rv_2)$, just because the parameter space is not separated, one just cannot say that $\Rv_1$ is specific sufficient for $\lambda_2$, one needs to take the precaution of plundging the parameter space in a separated one first, then apply the definition, or test the different conditions one by one and report on each of them.  
If one drops the condition $(\paramf,\Comp\paramf)\isoeq\mathscr{P}$, then one can always take $\left(\paramf:\range{P}\to,P\mapsto P^\Sufstat \right)$ and $\Comp\paramf:
            \range{P}\to,P\mapsto \left(:\tribu_\Omega\times \Omega\to, (A,\omega)\mapsto P(A\mid\tribu(\Sufstat),\omega)\right)$ and then $\paramf{}_{\paramf}\!\!\pperp_{\Comp\paramf} \Comp\paramf$ and say that $\Sufstat$ is
            a $(\paramf,\Comp\paramf)$-specific sufficient statistic.

\subsubsection{Fisher 1922}

\cite{Fisher1922Mathematical} gives definitions of sufficiency in two places
First he gives the following heuristic definition:
\fancyquote{\citet[p.~310]{Fisher1922Mathematical}}{
A statistic satisfies the criterion of sufficiency when no other statistic which can be calculated from the same sample provides any additional information as to the value of the parameter to be estimated.}

Then he adds a more mathematical definition.

\fancyquote{\citet[p.~316]{Fisher1922Mathematical}}{
[...] In mathematical language we may interpret [the] statement[: ``the statistic chosen should summrise the whole of the relevant information supplied by the sample''] by saying that if $\theta$ be the parameter to be estimated, [$X_1$[ a statistic which contains the whole of the information as to the value of $\theta$, which the sample supplies, and [$X_2$] any other statistic, then the surface of the distribution of pairs for values of [$X_1$] and [$X_2$] for a given value of $\theta$, is such that for a given value of $\theta$, the distribution of [$X_2$] does not involve $\theta$. In other words, when [$X_1$] is known, knwoledge of the value of [$X_2$] throws no further light upon the value of $\theta$.}

Ambiguity is that Fisher did not talk about the complement of $\theta$, and more discussion needs to be given about the parameter space, but right in the beginning, the definition about a sufficient statistic seemed to be the definition of a specific sufficient statistic.

\subsubsection{Fisher (1934) Two new properties of mathematical likelihood.}
\cite{Fisher1934}

\subsubsection{Neyman Pearson 1936}
\cite{NeymanPearson1936}
 
This paper acknowledges for the need ot theory of functions:
\fancyquote{\citet[p.~240]{NeymanPearson1936}}{it is inevitable, [...] that a paper dealing with [the] problem [of what conclusions regarding sufficient statistics may be drawn from the existence of uniformly most powerful tests, or vice versa] should bear some mark of the theory of functions, in spite of its concern with statistical questions}

In this paper, the model is parametrized by $n$ parameters, but there is no discussion about the joint space of those parameters. It is implicit in the paper that all conditional probabilities, densities are defined.
The author say they are not concerned about sufficiency for estimation, but for decision. They recall the results of Fisher linking sufficient statistics and uniformly most powerful tests.
Then they say this:

\fancyquote{\citet[p.~243]{NeymanPearson1936}}{In treating [the] problem of [the use of sufficient statistics with heir bearing on the theory of testing statistical hypotheses], we have found it necessary to use not only the conception of sufficient statistics as introduced by R.A. Fisher, but to introduce also some new conceptions which, as far as we are aware, have not been considered before, namely the conceptions of a sufficient set of statistics and of a shared sufficient statistic. [...]. We believe that our definition of a ``specific sufficient statistic'' corresponds to Fisher's conception of a sufficient statistic, but thought he has written on sufficient statistics in several places the definitions he has given appearm in our opinion, to leave some room for misunderstanding.} 
The authors refer to the paper [\cite{Fisher1922Mathematical}] and copy his definitions of sufficiency.
%They remark that the heuristic definition contains the word ``information'' that is not defined.
Then \citet{NeymanPearson1936} define a statistic.
\fancyquote{\citet[p.~243]{NeymanPearson1936}}{\color{black}
4. {\it Definitions and properties of sufficient statistics.}

[...]

{\it Definition I.} If a function $T$, or random variables $x_1,\ldots,x_n$ possesses the following properties:
\begin{enumerate}[label={\alph*)}]
 \item $T$ is defined and single valued at almost every point of the sample space $W$,
 \item whatever be a number $T'$, the locus of points in the sample space un which $T<T'$ is such that the probbility law of the $x's$ may be integrated over it, giving the probability $P(\{T<T'\})$,*
 \item there exist such values, $T'$, that the locus of points, $W(T'$, in which $T=T'$ is of at least $(n-1)$ dimensions, i.e. one less than the number of dimensions of the sample space $W$, \item $T$ does not depend upon any unknown parameters which may be involved in the probability law,

\end{enumerate}
then it will be called a statistic
}

They comment on that definition, and then define a sufficient statistic:
\fancyquote{\citet[p.~246]{NeymanPearson1936}}{
{\it Definition II.} The statistic $T$ is called a specific sufficient statistic with regard to the parameter $\theta_1$ if, whatever other statistic $T_2$ be taken, the relative probability law $P(T_2\mid T_1)$ of $T_2$, given $T_1$ is independent of $\theta_1$. This we believe to correspond with Fisher's original definition of a sufficient statistic. We have added the adjective ``specific'' for convenience in comparison.}

It is interesting to note that \cite{NeymanPearson1936} only retain one condition for a statistic to be sufficient, contrary to the later definition of \cite{Fraser1956} for example. Note that the definition is still imprecise, because of the ambiguity on the term independence  and the lack of discussion on the parameter space. The rest of the parameter defines a share sufficient statistic (when $T_1,...,T_n$ is a specific sufficient statistic. 
In their setup, they have a multidimensional parameter, but they have no assumption on the parameter space. Their definition is thus imprecise as they use the term ''depend on`` and we have seen that this is ambiguous.

\begin{definition}[Shared sufficient statistic]\cite{}
The statistic $T_1$ is called a shared sufficient statistic of the parameter $\theta_1,\ldots,\theta_q$ if, whatever other statistic $T_2$ be taken, the relative probability law of $T_2$ given $T_1$ is independent of these $q$ pqrqmeters, while it depends on the remaining $l-q$ parameters $\theta_{q+1},\ldots,\theta_l$.
\end{definition}
\citet[p.~247, Sec. 4, Proposition I]{NeymanPearson1936} characterize specific sufficiency by givin a necessary and sufficient condition on the likelihood, which looks like the factorisation theorem.

\fancyquote{\citet[p.~246]{NeymanPearson1936}}{
{\it Proposition I.} The necessary and sufficient condition for a statistic $T$ to be specifically sufficient with regard to a parameter $\theta$ (in the case of Definition II) is that in any point of the sample space (as defined on p.117), except perhaps fo a set of measure zero, it should be possible to present the probability law of the $x's$ in the form of the product $$p(x_1,\ldots x_n\mid\theta)=p(T\mid\theta)\phi(x_1,\ldots x_n)\mid_{T=T(x_1,\ldots,x_n)},$$, where $p(T\mid\theta)$ denotes the probability law of $T$, and $\phi$ is a function of the $x's$ independent on $\theta$.}
By writing $p(T\mid\theta)$ they seem to implicitely require that the distribution of $T$ only depends on $\theta$, which is in contradiction with their Definition II.

\subsubsection{Blackwell (1947), Conditional expectation and unbiased sequential estimation}
This early paper is important, as one of the importance of identifying the sufficient statistics, is to be able to only base the information on the sufficient statistic,

\fancyquote{\citet[p.~105]{Blackwell1947}}{It is shown that $E[f(x) E[y\mid x]]=E(fy)$ whenever $E(fy)$ is finite, and that $\sigma^2E(y\mid x)\leq\sigma^2 y$, where $E(y\mid x)$ denotes the conditional expectation of $y$ with respect to $x$. These results imply that whenever there is a sufficient statitic $u$ and an unbiased estimate $t$, not a function of $u$ only, for a parameter $\theta$, the function $E(t\mid u)$, which is a function of $u$ only, is an unbiased estimate for $\theta$ with a variance smaller than that of $t$. }

\subsubsection{Halmos \& Savage (1949). Application of the Radon-Nikodym theorem to the theory of sufficient statistic}
\cite{HalmosSavage1949} show how in a dominated model, after applying the Radon Nikodym theorem, one can deduce sufficiency from the densities (factorisation).
The result was already stated in \citet[p.~247, Sec. 4, Proposition I]{NeymanPearson1936}, \cite{HalmosSavage1949} provide a very general result, by using Radon-Nikodym theorem. Nevertheless, they 
do not mention specific sufficiency, only sufficiency. to say that but This paper also discusses pairwise sufficiency and likelihood ratio.

Section 1. starts with definitions from the theory of measure. Section 2. is about Measures and their derivatives. The Radon-Nikodym theorem. Section 3. is about probabilities and expectations. It contains the following lemma:

\fancyquote{\citet[p.~229]{HalmosSavage1949}}{
\textsc{Lemma} 4.If $\mu$ and $\nu$ are measures on $\mathbf{S}$ such that $\nu<<\mu$ then $\nu T^{-1}<<\mu T^{-1}.$}. Other classical lemma are given. Section 4 is about dominated sets of measures. Lemma 7 is the one that states that every dominated set of measures has an equivalent countable set.
Section 5 is titled Sufficient statistics for dominated sets. It contains the definition of a sufficient statistic for a set of measures on $S$. 
It does not contain a defintion of a specific statistic.
The authors develop the concept of pairwise sufficiency, which is sufficiency for any two measures from the model. Pairwise sufficiency is weaker than sufficiency, but equivalent in the dominated case.

\subsubsection{Dinkyn (1951), Necessary and sufficient statistics for a family of probability distributions}
This paper is cited by \cite{Andersen1967}
DYNKIN,
Necessary
sufficient statistics
family
. Vspechi Matem. Nauk (N.S.),
probability Uspechi
vol. 6, pp. (1-41), 68-90.

This paper contains results about the existence of minimal sufficient a·field.

\subsubsection{Fraser (1952), Sufficient statistics and selection depending on the parameter}
\cite{Fraser1952SufficientStatisticsandSelection} defines what is a functional sufficient (or f-sufficient) function with respect to a family of functions, as well as what is meant by parameter of selection. An index on the probability measures of the model is a parameter of selection if all the measures of the model are inductions of the same measure on subsets of the observation space.

Fraser establishes some following result, as for example:

\fancyquote{
\citet[p.~419]{Fraser1952SufficientStatisticsandSelection}}{\textsc{Lemma 3.}\textit{Any sufficient statistic for a dominated set of measures $\mathfrak{M}<<\lambda$ is an f-sufficient statistic for an equivalent set of densities (relative to $\lambda$)}}

I was mislead by the title. Here selection does not refer to the process of selecting units in a population.

\subsubsection{Rao (1952). Minimum variance estimation in distributions admitting ancillary statistics.}
This paper explains how to use an ancillary statistic to improve estimates.

\subsubsection{Bahadur, (1954). Sufficiency and statistical decision functions}
This paper is cited by \cite{Andersen1967}

\subsubsection{Basu (1955), On statistics independent of a complete sufficient statistic.}
I was hoping to find in \cite{Basu1955a} the definition of an ignorable statistic, that I only found in \cite{Schervish1995}, an ignorable statistic being a statistic independent from a sufficient statistic.
I have had a hard time to find the definition elsewhere, although I had this definition in my lecture notes as a student.
In this paper, it is not question of nuisance parameters.
It contains a theorem saying that the distribution of a statistic independent from a sufficient statistic $\Rv$ for $\paramf=\Id_\range{P}$ is necessarily the same for any $P$, e.g 
($(P\mapsto P^{\Id\mid \Signal})\dependsonly\mathds{1}_{\range{P}}$ and $\Rv\perp\Signal$ $\Leftrightarrow$ $(P\mapsto P^\Rv\dependsonly\mathds{1}_\range{P}$, e.g. being independent from a sufficient statistic implies being an ancillary statistic. If in addition $\Signal$ is boundedly complete, the reciprocal is true. This paper is important though at its generalisation in the presence of nuisance parameters was discussed later. Note that Basu does not use the term ancillary, but independent of $\param$, although ancillary was introduced by \cite{Fisher1934}, and was reused by \cite{Rao1952Minimumvariance}.

\subsubsection{Fraser 1956}
\cite{Fraser1956} Proposes a definition of sufficient statistic in presence of nuisance parameters.
This paper does not have any assumption on the measured space (it is not necessarily euclidean with Borel sigmafield). It is limited to the case where the parameter space is of the form $\Theta\times\bar\Theta$, and does not consider a subset of $\Theta\times\Comp\Theta$. The paper also suggest that they may not exist a specific sufficient statistic.This is mentioned in the paper. The paper is decision theory orientated and cites nice results, which are about testing based on a statistic sufficient wrt $\theta$. And there is also a theorem about unbiased estimator with minimum variance. The definition coincides with the one I give.
\fancyquote{\citet[p.~229]{Fraser1956}}{
For some of these problems a generalized definition of sufficiency can be applied. Let $X$ be a random variable over the measurable space $\mathcal{X}(\mathcal{A})$ and let $\{P_{\theta,\eta}\mid (\theta,\eta)\in\Theta\times H\}$ be the class of possible probability measures for $\mathcal{X}$. Also, let $t(x)$ be a statistic mapping $\mathcal{X}(\mathcal{A})$  into the measurable space $\mathcal{I}(\mathcal{B})$ and let $P^T_{\theta\eta}$ designate the measure on $\mathcal{I}(\mathcal{B})$ induced by $t(x)$ from the measure $P^T_{\theta\eta}$, over $\mathcal{X}(\mathcal{A})$. Then we propose the following extension of the concept of sufficiency: {\it $t(x)$ is a sufficient statistic $(\theta)$ for the class of measures $\{P^T_{\theta\eta}\mid (\theta,\eta)\in\Theta \times H\}$ if there exists a function $P_\eta(A \mid t)$ such that
$$(1)\quad\quad\quad\quad\quad\quad P_{\theta,\eta}(A\cap t^{-1}(B)) \equiv \int_B P_\eta(A \mid t)dP_\theta^T(t)$$
for all $A \in \mathcal{A}$, $B \in \mathcal{B}$ where the induced measure of $t(x)$, $P^T_\theta$, is independent of $\eta$.}

The conditional probability that $X$ falls in the set $A$ given $t(X) = t$ is given by a function which will serve as the integrand in the integral of (1). The definition says that this conditional probability must depend only on the nuisance parameter $\eta$, and that the marginal distribution of the statistic $t(x)$ should depend only on the parameter of interest $\theta$. Thus it can be seen intuitively that the statistic $t(x)$ is in a general sense sufficient for problems concerning the parameter $\theta$.
For the particular case in which there are no nuisance parameters, this defini- tion reduces to the ordinary definition of sufficient statistic. However, there need not exist a sufficient statistic (0), whereas there always exists a sufficient statistic by the usual definition. Another drawback to the formulation above is the requirement that the parameter space be a Cartesian product.}

{\textbf{$(\Rv,\Comp\Rv)$-Extension and $(\paramf,\Comp\paramf)$ extensions of a model}}

Here, we discuss why it makes sense to impose $(\paramf,\Comp\paramf)\isoeq\range{P}$ as a necessary condition for partial sufficiency.

First, we show how one for any random variable $\Rv$, one may consider that any model $\range{P}$ can be ''plunged`` in another model $\Extension{\range{P}}$ that will lead to the same inference and such that $(P\mapsto P^\Rv,P\mapsto P^{\Id\mid\Rv})\isoeq\Extension{\range{P}}$, we will call this extension a $(\Rv,\Comp\Rv)$ extension of $\range{P}$, and we explain how to move from the inference on $\range{P}$ to the inference on $\range{P}$ as long as the object of the inference as a function of $(P\mapsto P^\Rv,P\mapsto P^{\Id\mid\Rv})$ can be extended to $\image(P\mapsto P^\Rv)\times\image(P\mapsto P^{\Id\mid\Rv})$ .

Given a random variable $\Rv$, and a non necessarily distinct complement $\Comp\Rv$.Then 
a (non necessarily distinct) separation of $\range{P}$ is $((\range{P}\to,P\mapsto P^\Rv),\range{P}\to,P\mapsto P^{\Comp\Rv\mid\Rv})$.
We define the extension of $\range{P}$ based on $\Rv$ and $\Comp\Rv$ as the model
$\left(\Extension{\Omega},
\tribu_{\Extension{\Omega}},
\Extension{\range{P}}\right),$ where 
$\Extension{\Omega}=\Omega\times\{0,1\}$,
$\tribu_{\Extension{\Omega}}$ is the sigmafield $\tribu_\Omega\times\{\emptyset,\{0\},\{1\},\{0,1\}\}$ ,
$\Extension{\range{P}}=\left\{\left.P\otimes\Dirac_{\{1\}}\right| P\in\range{P}\right\}\cup \left\{\left.\left<P^\Rv\mid P'^{\Comp\Rv\mid\Rv}\right>\otimes\Dirac_{\{0\}} \right| (P,P')\in\range{P}, 
\left<P^\Rv\mid P'^{\Comp\Rv\mid\Rv}\right>\notin\range{P}\right\}$
$\Extension{P}_{\param,\Comp\param}$.  
for $\somesubset\in \Rv(\tribu_\Omega)$, $i\in\{0,1\}$, 
% \begin{equation}\label{eq:extension}
% \Extension{P}_{\param,\Comp\param}(\somesubset\times\{i\})=\left<\param\mid\Comp\param\right>(A)\times\left|
%                             \begin{array}{ll}
%                                  1 &\text{if } i=1 \text{ and }(\param,\Comp\param)\in (\paramf,\Comp\paramf)(\range{P})\\
%                                  1 &\text{if } i=0 \text{ and }(\param,\Comp\param)\notin (\paramf,\Comp\paramf)(\range{P})\\
%                                  0 &\text{otherwise}
%                             \end{array}\right.\end{equation}
For any statistic $\Obs:(\Omega,\tribu_{\Omega})\to$, 
define the statistic $\Extension{\Obs}:(\Extension{\Omega},\tribu_{\Extension{\Omega}})\to, (\omega,i)\mapsto (X(\omega),i)$.
For any parameter $\paramf:\range{P}\to$, such that $\paramf\dependsonly (\range{P}\to, P\mapsto (P^\Rv,P^{\Comp\Rv\mid\Rv}))$ $\exists \paramf':\to$ such that $\paramf=\paramf\circ(P\mapsto (P^\Rv,P^{\Comp\Rv\mid\Rv})$,  and so it is possible to define the extension of $\paramf$ as $\Extension\paramf=\paramf'\circ(\Extension{\range{P}}\to P\mapsto (P^\Rv,P^{\Comp\Rv\mid\Rv})$.
The likelihood-based inference on $\paramf$ given $\Obs=\obs$ with the original model will be equivalent to the likelihood based inference on $\Extension\paramf$ given $\Extension\Obs=(\obs,1)$.
In a Bayesian framework, we also extend the sigmafield on $\range{P}$ to $\mathring{\range{P}}$ by taking any sigmafield $\tribu_{\mathring{\range{P}}}$ on $\mathring{\range{P}}$ containing all elements from $\tribu_{\range{P}}\times\{1\}$, and we define the set of prior distributions on $\Extension{\range{P}},\tribu_{\Extension{\range{P}}}$ as the set $\Extension\priordistset=\{\Extension\priordist\mid\priordist\in\priordistset\}$, where $\Extension\priordist$ is any measure on $(\mathring{\range{P}},\tribu_{\mathring{\range{P}}})$ such that $\forall \somesubset \in\tribu_{\range{P}}$, $\Extension{\priordist}(\somesubset\times\{1\})=\priordist(\somesubset)$.

% 
% This definition, requires forcing $P_\eta(A\mid t)$ to be independent on $\theta$ and $P_\theta^T$ to be independent on $\eta$, which is the reason why there may not exist a specific sufficient statistic. 
% One could have introduced a weaker condition, so that any sufficient statistic is also a specific sufficient statistic.
% Another remark with this definition is that it requires to specify what the complement $\Comp\theta$ (here $\eta$) is. 
% 
% For our purpose, which is to find a way to say that a part of the observation can be ignored, a weaker condition, that  would consist in only imposing that $P_\eta(A\mid t)$ to be independent on $\theta$, would be enough.
% Consider for example the case where
% 
% $Y_1\mid Y_2\sim\mathrm{Normal}(\theta Y_2,\sigma^2)$ and $Y_1=Z_1$, $Z_1\perp Z_2$, design is Poisson sampling (without replacement), function of $Z_1,Z_2$ two Bernoulli variables and  parameter $\Comp\theta$: $P(\{I_\individual=1\mid \Design\})=\Comp\param_0+\Comp\param_1\Design_{1,k}+\Comp\param_2\Design_{2,k}$, selection is independent of $(Y,Z)$ conditionnally to the design. 
% Then the likelihood of 
% $\likelihood(\expected[Y_1\mid Y_2],(Y[R])_{1,\ldots,[L]})(\theta,(y_1,\ldots,y_l))$ is, when the population size is a parameter, 

\cite{LehmannBookTesting} General theory book Definition of sufficiency in the non dominated case.

\subsubsection{Rasch (1960), On general laws and the meaning of measurement
meaning
in psychology}
and also Rasch, 
-- (1960). Probabilistic Models for some Intelligence and Attainment
Probabilistic
Intelligence
Attainment

This paper is cited by \cite{Andersen1967}

\subsubsection{On the foundations of statistical inference - Birnbaum (1962)}
This paper does not contain many formulae, there is a lot of text. It is a deep reflexion about certain principles in statistics. Sufficiency and Conditionality Principles are presented there. There is a theorem that says that Sufficiency and conditionality principle implies likelihood principles. This text is necessary to understand discussion of \cite{Dawid1975Ontheconcept}. Note that {\it informative inference} is used here, but not defined by opposition to another statistical inference approach. 

\fancyquote{\citet[p.~275]{Birnbaum1962}}{Since the problem-area of informative inference has not received a generally accepted delineation or terminology, it will be useful to note here some of the terms and concepts used by writers representing several different approach}.
Notions like informative and uniformative experiments are explained 
\fancyquote{\citet[p.~291]{Birnbaum1962}}{Since this experiment gives the same distribution on the two-point sample space under each hypothesis, it is completely uninformative, as is any outcome of this experiment. According to the likelihood principle, we can therefore conclude that the given likelihood function has a simple evidential interpretation, regardless of the structure of the experiment from which it arises, namely, that it repre- sents a completely uninformative outcome. (The same interpretation applies to a constant likelihood function on a parameter space of any form, as an essentially similar argument shows.)}

For the current paper objective, what we may retain from this paper is that it uses the term informative and uninformative for experiments, not for statistics. An experiment is uninformative when the likelihood function is the same for different values of the parameters basically. This does not help to define an uniformative process, but it is worth mentioning that in ''informative inference``, definition of informative experiment is based on the likelihood.
\cite{Basu1977} will comment on this paper. 
\subsubsection{Barnard (1963), Some logical aspects of the fiducial argument}

I read \cite{Barnard1963} because it was cited by \cite{Dawid1975Ontheconcept}.
As its name indicates, it contains comments about the fiducial argument (Fiducial argument is about considering that parameters as  random variables). The paper is not helpful with respect to defining specific-informative latent process in presence of a nuisance parameter.

\subsubsection{Fraser (1964), Local Conditional Sufficiency}
This paper is about local conditional sufficiency. It is not mention of nuisance parameter, and is not helpful with respect to defining specific-informative latent process in presence of a nuisance parameter.

\subsubsection{Le Cam (1964), Sufficiency and approximate sufficiency.}
Ouf.

\subsubsection{Sverdrup (1966), The present state of the decision theory and the Neyman-Pearson theory}
This paper contains a characterisation of the sufficiency not based on likelihood, and gives an overview of basic mathematical definitions from probability theory.

\fancyquote{\cite{Sverdrup1966Present}}{
The starting point is a sigmafield $\mathscr{A}$ of subsets [,,,] on the sample space $\range{X}$ of sample points $x$. A subsigmafield $\mathscr{A}_0$ of $\mathscr{A}$ may be of special interest. Let $P(A)$ be a probability measure over $\mathscr{A}$. Furthermore $P(A\mid \mathscr{A}_0, x)$ is the conditional probability of $A$ given ''the most accurate description of $x$ by means of statements from $\mathscr{A}_0$``. It is defined as the almost unique $\mathscr{A}_0$-measurable function of $x$ which satisfies
$[\forall B\in \mathscr{A}_0],\ P(A\cap B)=\int_B P(A\mid \mathscr{A}_0,x) \derive P$. $P(A\mid\mathscr{A}_0,x)$ does always exist and is almost uniquely defined[...].

Let $\range{P}$ be a family of probability measures $P$ for a random variable $X$ over $(\range{X},\mathscr{A})$. A subsigmafield $\mathscr{A}_0$ of $\mathscr{A}$ is sufficient for the family $\range{P}$ if 
$\forall A\in\mathscr{A}$, [$\exists\Phi_A(x)$ a $\mathscr{A}_0$ measurable function which for all $P$ is the conditional probability relatively to $\mathscr{A}_0$, i.e. for which 
 $$\forall P\in\range{P},\ \Phi_A(x)=P(A\mid\mathscr{A}_0,x)\text{ a.e. }(\mathscr{A}_0, P).$$

We can rewrite the definition more consisely

A subsigmafield $\mathscr{A}_0$ of $\mathscr{A}$ is sufficient for the family $\range{P}$ if 
$$\forall A\in\mathscr{A},\ \exists\Phi_A:(\Omega,\mathscr{A}_0)\to(\mathbb{R},\mathscr{B}_\mathbb{R}),\ \forall P\in\range{P},\ \Phi_A(x)=P(A\mid\mathscr{A}_0,x)\text{ a.e. }(\mathscr{A}_0, P).$$
}

The paper explains that there are different definitions of sufficient statistics (p. 312). Classical, for decision, for Bayesian. Also discuss ``specific sufficiency'' in presence of ``nuisance parameter''.  The paper Gives the definition of an ancillary statistic by Fisher, cites papers by Basu and Rao about nuisance parameters, and explains that \cite{Fraser1956} do things differently by requiring that the distribution of a sufficient statistic should be independent of the nuisance parameter, which is a very restrictive property. Gives the definition of a sufficient statistic in the general case (Eq 1. p. 313). Note that the paper redefines from the begining conditional probabilities.

This is an important summary quote from their paper (p. 312),

\fancyquote{\citet[p.~312]{Sverdrup1966Present}}{\color{black}
[\cite{NeymanPearson1936}] and Fisher were talking about specific sufficiency relatively
to a certain parameter. Let $\theta = (\rho, \tau)$, where $\rho$ is the decision parameter, i.e. the parameter which the decision situation is concerned with, whereas $\tau$ is the nuisance parameter. {\bf If $\Samplepopmap(x)$ for any given $\tau$  is minimal sufficient, then $R (x)$ is specifically sufficient for $\rho$}. Suppose e.g. that $X_1, \ldots , X_n$, are independent normal $(\xi, \sigma)$. Then $\bar{X}$ is specifically sufficient for $\xi$. Whereas minimal sufficient statistics exist under very general conditions, this is far from being true of specific sufficient statistics. Thus in the example just given no specific sufficient statistic exists for $\sigma^2$. Because for any given $(\xi)$ it is 
$\sum_(X_j - \xi)^2$ which is the minimal sufficient statistic, but this is not a statistic if $\xi$ is unknown. 

It seems doubtful whether specific sufficiency in the sense taken above is an important concept in the decision theory. It is difficult to find any direct connection between this concept and decision functions. Consider e.g. Student's situation with 
testing of $\theta=0$ (or constructing confidence interval for $\xi$). Then $\bar{X}$ is specifically
sufficient, but in order to perform the testing we have to consider $\sum (X_j - \bar{X})^2$. It is
of course $(\bar{X}, \sum(X_j - \bar{X})^2)$ which is the minimal sufficient statistic for the model. 
Fisher [9] was aware of the difficulty and introduced the concept of ancillary statistic. $T(x)$ is ancillary if it jointly with the specific sufficient statistic $\Samplepopmap(x)$ is minimal sufficient and the probability distribution of $T (x)$ only depends on the nuisance parameter $\tau$. Rao [25] and Basu [3] have proved some interesting mathematical properties about ancillary statistics.
A very interesting approach from a statistical point of view, is due to [\cite{Fraser1956}]. He does not need the concept of ancillary statistic, instead he adds to the above definition of specific sufficiency the property that the distribution of $R (x)$ shall be independent of the nuisance parameter. This is a rather restrictive property ($\bar{X}$ in the example above is then not specifically sufficient). On the other hand he is then able to establish links with decision problems. He shows that by testings and point estimations concerning $\rho$, the statistician may limit himself to procedures depending on $R (x)$ without losing power or efficiency.
Below we shall expand upon some, but not all, of the ideas which we have sketched above.
}

The paper continues with important comments, including Section II.E (Tests with optimal power, Justification of Conditioning)

\subsubsection{Andersen 1967, On partial sufficiency and partial ancillarity}
\cite{Andersen1967} 
This paper gives the definition of partial sufficiencty and partial ancillarity.
In the abstract, it explains that the target is not to discuss the 

\subsubsection{Durbin 1969. Inferential aspects of the randomness of sample size in survey sampling}
\cite{durbin1969inferential} is cited by \cite{Dawid1975Ontheconcept}.

\subsubsection{Durbin(1970), On Birnbaum's theorem on the relation between sufficiency, conditionality and likelihood}
\cite{Durbin1970} This note contains comments about Birnbaun's theorem (see notes on \cite{Birnbaum1962}

% \subsubsection{Birnbaum 1970}
% \cite{Birnbaum1970}

\subsubsection{Sandved, A principle for conditioning on an ancillary statistic. 1967}.

\cite{Sandved1966Principle} gives a definition of ancillarity for a decision problem.
\fancyquote{\citet[p.~1]{Sandved1966Principle}}{
A statistical decision problem can usually be formulated in the following way: $X$ is a stochastic variable with probability measure $P$. (In this paper the stochastic variables and the parameters may be multidimensional.) A priori $P\in\mathscr{P}$, where $\mathscr{P}$ is a given class of prbability measures. Based upon an observation of $X$ we shall choose a subclass of $\mathscr{P}$ which we believe contains $P$, or we shall give an estimate for $P$.}
\fancyquote{\citet[p.~6.]{Sandved1966Principle}}{\color{black}
Let $\Theta$ be a set of indices and $\{\mathscr{P}_\theta,\theta\in\Theta\}$ be a family of non-empty subclasses of $\mathscr{P}$ such that if $P\in\mathscr{P}$, then $P$ belongs to one and only one $\mathscr{P}_\theta$. We want to make inferences on $\theta$ on the basis of an observation of $X$.

For instance, if we want to test a hypothesis, we may let $\Theta$ consist of two elements, $\theta_0$ and 
$\theta_1$, i.e. to accept or to reject the hypothesis.

Let $\mathbf{a}(X)$ be $\mathscr{A}$-measurable. $a(X)$ induces a sub-$\sigma$- algebra $\mathscr{A}^{\mathbf{a}}$ of 
$\mathscr{A}$. Let $P^\mathbf{a}$ be the measure $P$ restricted to $\mathscr{A}^{\mathbf{a}}$, and let $\mathscr{P}^{\mathbf{a}}_\theta$ be the class of all $P^{\mathbf{a}}$ with $P\in \mathscr{P}_\theta$.
Let $\mathbf{a}(X)$ be a statistic such that 
\begin{enumerate}[label={\roman*)}]
 \item The classes $\mathscr{P}_\theta^{\mathbf{a}}$, $\theta\in\Theta$ are identical.
 \item The class of conditional probability distributions of $X$, given $\mathbf{a}(X)$, $\theta$ and $P^{\mathbf{a}}$, is independent of $P^{\mathbf{a}}$.
\end{enumerate}
We then define $\mathbf{a}(X)$ to be an ancillary statistic for the decision problem, and we propose the following principle: 
In the decision problem at hand start with the conditional distribution of $X$ given $\mathbf{a}(X)$.}

\subsubsection{Barndorff-Nielsen (1973)}
\cite{BarndorffNielsen1973} definitions are discussed in two papers later \cite{Dawid1975Ontheconcept}, and \cite{Basu1977}.
This is the citation from \citet[p.~249]{Dawid1975Ontheconcept} (notations were changed to match ours, and for \cite{Dawid1975Ontheconcept}, one can consider that $\Omega=\range\Obs$ and $X=\Id_\Omega$):
{\em ``The data of a statistical experiment may be expressed as a random quantity [$\Obs:(\Omega,\tribu_\Omega)\to\range\Obs$] taking values in a space [$\range\Obs$] and with a [set of distributions $\{P^X\mid P\in\range{P}\}$] over [$\range\Obs$] [...]. Suppose that [$\paramf:\{P^X\mid P\in\range{P}\}\to$]  is a parameter-function. Let [$U=\mathbf{u}(X)$,$V=\mathbf{v}(X)$] be statistics. Then we can ask what useful meaning can be assigned to the statements:
\begin{enumerate}
 \item $U$ is sufficient for $\paramf(P^\Obs)$;
 \item $V$ is ancillary for $\paramf(P^\Obs)$. 
\end{enumerate}
We shall be primarily concerned with the following one of many attemps to give precise meaining to these ideas. \cite{BarndorffNielsen1973} has called a statistic $U$ a cut with corresponding parameters [$\paramf,\Comp\paramf$] if 
\begin{enumerate}
 \item [$\paramf(P)$] and [$\Comp\paramf(P)$] vary independently as [$P$]  varies in $\range{P}$, and are together equivalent to $\{P^\Obs\mid P\in\range{P}\}$ [e.g. $(\paramf,\Comp\paramf)\isoeq\{P^X\mid P\in\range{P}\}$,ed]
 \item The distribution [$P^{U}$] depends on $\theta$ alone [e.g, $\exists \paramf':\allfunc{\paramf(\{P^X\mid P\in\range{P}\})}{}$ such that $(P\mapsto P^{U})=\paramf'\circ\paramf$,ed] and 
 \item The conditional distribution [$P^X\mid U$] depends only on [$\Comp\paramf(P)$] [ed, e.g. $\exists \Comp\paramf':\allfunc{\Comp\paramf(\{P^X\mid P\in\range{P}\})}{}$ such that $(P^X\mapsto (P^X)^\mathbf{u})=\Comp\paramf'\circ\Comp\paramf$] 
\end{enumerate}
}

\subsubsection{Dawid, 1975}
\cite{Dawid1975Ontheconcept} shows ''ambiguities [that] may arise from attempts to define and apply analogues of sufficiency and ancillarity in the presence of nuisance parameters``.
\cite{Dawid1975Ontheconcept} refers to \cite{BarndorffNielsen1973}. 
The general definition of S sufficiency is recalled.

\fancyquote{\citet[p.~ 248]{Dawid1975Ontheconcept}}{ \color{black}

The data of  statistical experiment may be expressed as a random quantity $x$ taking values in a space $\mathscr{X}$, and with a distribution $P_\omega$ ober $\mathscr{X}$ depending on a parameter $\omega$ taking values in $\Omega$. Suppose that $\theta=\theta(\omega)$ is a parameter-function (henceforth called simply ''parameter``). Let $u$, $v$ be statistics. Then we can ask what useful meaning can be assigned to the statements:
\begin{enumerate}[label={\roman*)}]
 \item $u$ is sufficient for $\theta$;
 \item $u$ is ancillary for $\theta$.
\end{enumerate}
We shall be primarily concerned with the following one of many attemps to give precise meaining to these ideas. \cite{BarndorffNielsen1973} has called a statistic $u$ a cut with corresponding parameters $(\theta,\Phi)$ if 
(i) $\theta$ and $\Phi$  vary independently as $\omega$  varies in $\Omega$, and are together equivalent to $\omega$ (ii) The distribution $P_{u,\omega}$ of $u$ depends on $\theta$ alone  and 
 (iii) the conditional distribution $P_\omega^u$ of $x$ given $u$ involves only $\Phi$.
 In the dominated case, 
 $(ii)$ and $(iii)$ above are equivalent to 
 \begin{equation*} f(x\mid \omega)=f(u\mid \theta)f(x\mid u,\Phi),\tag{(1.1)}\end{equation*}

 where $f$ represents a density with respect to an appropriate underlying measure.
 (If $Q$ is a probability distribution which dominates $\{P_\omega\mid\omega\in\Omega\}$, the respective measures in $(1.1.)$ may be taken as $Q$, $Q_u$ and $Q^u$.)
 
If (1.1) holds, there seems to be an intuitively appealing case for supposing that all relevant information about $\theta$ may be extracted by referring $u$ to its marginal distribution; while for inference about $\Phi$, we can refer $x$ to its distribution conditioned on $u$. This prompts the following definitions (Barndorff-Nielsen, 1973-after Fraser, 1956; Sverdrup, 1966; Sandved, 1967).

If $t$ is a statistic and $\lambda = \lambda(u)$ a parameter, we call $t$ {\it S-sufficient} for $\lambda$ if there exists some parametrization $\omega = (\theta, \Phi)$ as above, and a corresponding cut $u$, such that $u$ is a function of $t$ and $\lambda$ is a function of $\theta$. Similarly, $t$ is {\it S-ancillary} for $\mu$ if $t$ is a function of $u$ and $\mu$ is a function of $\Phi$. }

\cite{Dawid1975Ontheconcept}[Sec. 2-3] gives a historical overview of the issue: How to define specific sufficiency, and where the debate comes from. 
Basu cite \cite{Barnard1963},\cite{Andersen1967},\cite{durbin1969inferential},
\cite{Durbin1970},\cite{Birnbaum1962},\cite{BarndorffNielsen1973}.

\cite{Dawid1975Ontheconcept}[Sec. 3] mentions boundedly completeness for $\theta$ in presence of a nuisance parameter without defining it. It is implicit that boundedly complete for $\theta$ means boundedly complete when the complement $\phi$ of $\theta$ is fixed.

\subsubsection{Basu (1975), Statistical Information and Likelihood}
\subsubsection{Basu, 1977 - On the Elimination of Nuisance Parameters}
\cite{Basu1977} summarizes and critics all that precedes and that is related to definition of partial and specific sufficiency, although the definitions are given in a dominated context. If one paper has to be given as a reference, it should be this one.

\citet[Sec.2, Marginalization and Conditioning]{Basu1977} explains in what consists the marginalization procedure (when  a sufficient statistic is identified) and the conditioning procedure (when an ancillary statistic is identified).
\citet[Sec.3, Partial sufficiency and partial ancillarity]{Basu1977} contains the definitions of 
\begin{enumerate}
\item an ancillary statistic $\Rv$: $(P\mapsto P^\Rv)\dependsonly\mathds{1}$, 
\item a $\paramf$-oriented statistic $\Rv$ : $(P\mapsto P^\Rv)\dependsonly\paramf$,
\item variation independence, 
\item a $(\paramf,\Comp\paramf)$-specific sufficient statistic $\Rv$ : $[(\paramf,\Comp\paramf)\isoeq\range{P}$ and $(P\mapsto P^{\Id_\Omega\mid \Rv})\dependsonly\Comp\paramf]$, 
\item a $(\paramf,\Comp\paramf)$-specific ancillary statistic $\Rv$ : $[(\paramf,\Comp\paramf)\isoeq\range{P}$ and $(P\mapsto P^{\Rv})\dependsonly\Comp\paramf]$, 
\item a $(\paramf,\Comp\paramf)$-partially sufficient  (or $(\paramf,\Comp\paramf)$-Sandved-sufficient or $(\paramf,\Comp\paramf)$-S-sufficient) statistic   is a $(\paramf,\Comp\paramf)$-specific sufficient and $\paramf$-orientated statistic : $[(\paramf,\Comp\paramf)\isoeq\range{P}$ and $(P\mapsto P^{\Id_\Omega\mid \Rv})\dependsonly\Comp\paramf$ and $(P\mapsto P^{\Rv})\dependsonly\paramf]$
\item a $(\paramf,\Comp\paramf)$-partially ancillary (or $(\paramf,\Comp\paramf)$-Sandved-ancillary or $(\paramf,\Comp\paramf)$-S-ancillary) statistic  is a $(\paramf,\Comp\paramf)$-specific ancillary and $(\Comp\paramf,\paramf)$-specific sufficient statistic : $[(\paramf,\Comp\paramf)\isoeq\range{P}$ and $(P\mapsto P^{\Id_\Omega\mid \Rv})\dependsonly\paramf$ and $(P\mapsto P^{\Rv})\dependsonly\Comp\paramf]$.
\item a Barndorff cut is a statistic $\Rv$ such that $(P\mapsto P^\Rv,P\mapsto P^{\Id_\Omega\mid\Rv})\isoeq\mathscr{P}$.
\end{enumerate}

Basu then explains that given $\paramf$, there does not necessarily exist a Barndorff cut $\Rv$ such that $\paramf:P\mapsto P^\Rv$, and gives examples.

Then Basu gives equivalent definitions in the case of a dominated model.

\fancyquote{\citet[p.~356 ]{Basu1977}}{{\it Definition 3 ($\theta$- Oriented Statistic):} The statistic $T$ is $\theta$ oriented if the marginal distribution of $T$ depnds on $\omega$ only through $\theta=\theta(\omega)$. That is $\theta(\omega_1)=\theta(\omega_2)$ implies $P_{\omega_1}(T^{-1}B)=P_{\omega_2}(T^{-1}B)$ for all $B\in\mathcal{B}$.

Observe that every ancillary statistic is $\theta$-oriented irrespectrive of what $\theta$ is.}

\fancyquote{\citet[p.~356 ]{Basu1977}}{{\it Definition 4 (Variation Independence):} The two functions $\omega\to a(\omega)$ and $\omega \to b(\omega)$ on the space $\Omega$ with respective ranges $A$ and $B$ are said to be variation independent if the range of the function $\omega\to (a(\omega), b(\omega))$ is the Cartesian product $A\times B$.}

\fancyquote{\citet[p.~357 ]{Basu1977}}{If the universal parameter $\omega$ can be represented as $\omega=(\theta,\phi)$, where $\theta$ and $\Phi$ are variation independent in the preceding sense - that is $\Omega=\Theta\times\Phi$ where $\Theta$ and $\Phi$ are the respective ranges of $\theta$ and $\phi$ - then we cal $\phi$ a variation independent complement of $\theta$. With $\theta$ as the paremter of interest, we may then call $\phi$ the nuisance parameter. 
We have not come across a satisfactory definition of the notion of a nuisance parameter. It is only hoped that the above working definition will meet with little resistance. (See \cite{BarndorffNielsen1973} for durther details on the notion of variation independence.}

\subsubsection{Dawid (1979) - Conditional Independence in Statistical Theory}

\subsubsection{Basu, Pereira (1983) Conditional Independence in Statistics}

\subsubsection{Godambe, 1984, On Ancillarity and Fisher Information in the Presence of a Nuisance Parameter}

\subsubsection{Yamada Morimoto - 1992 - Sufficiency}

\subsubsection{Severini, 1993, Local Ancillarity in the Presence of a Nuisance Parameter}

\subsubsection{Zhu, Reid, 1994 - Information, Ancillarity, and Sufficiency in the Presence of Nuisance Parameters}

\subsubsection{Fraser (2004) Ancillaries and Conditional Inference}
\cite{Fraser2004}

\subsubsection{Ghosh Reid and Fraser, 2010, Ancillary Statistics: a review}
\cite{GhoshReidFraser2010Ancillary}
This paper cites Basu1964, and a lot of others. It deals with nuisance parameters. It refers to Bayesian ancillarity by Severini. 

\section{Notes on selected works related to the definition of ignorable design of missing data mechanism, and informative selection}
In this section, we tracked the appearance of the term informative and ignorable in the litterature.
We provide some notes on a selection of articles, and present them in an alphabetical order. We also provide notes on some fundamental articles that precede the appearance of those terms, as the need to take into account the selection process was identified before some names were put on situations where it had  or had not to.

\subsection{Early papers by Godambe and Ericson}
\citet{Godambe1966a} explains why the likelihood function is ``uninformative'' in design based inference, without giving a definition of the term. In this paper, a particular definition of indepent linear estimators, linear sufficiency for linear 
 combinaisons of observations, are given. Those concepts did not prove popular afterwards and there is not much in this paper to help us define properly what is informative selection.
 
\subsubsection{Godambe 1969, ``A fiducial argument with application to survey sampling''} 
 
\citet{Ericson1969} A Bayesian framework, (the first for Bayesian ?). No mention of informative selection.
``No new principles of inference are necessary''. Develops results on the posterior mean in special cases (special prior, special scheme).

\subsection{Two papers by Scott, 1975, 1977}
% \citet{Holt1980Regression}  Refers to Scott 1975
% A basic assumption of model (2.1) is that $E(\varepsilon\mid X_2) = 0$. The property must hold after selection of the sample if least squares methods are to be used, and in this case the design can be said to be uninformative. If, after selection, $E(\varepsilon\mid X_2)\neq 0$ then the design is informative. (See, Scott, 1975.) 

The paper by \citet{Scott1977} that follows \cite{Scott1975} is particularly important for two reasons:
\begin{enumerate}
 \item it is likely the first one that actually provides a definition for a non-informative design.
 \item it undelines the importance of taking into acccount what is observed in the definition of an informative design.
\end{enumerate}

The paper also states results such as ``simple random sampling'' is the only non informative design, and underlines the importance of the design variables and the importance of their use if they are available to the analyst during the analysis phase.

Nevertheless, two main problems remain:
First, the authors did not consider sampling with replacement, or taking into account the order in which units were selected, or other possible observed or latent information related to the selected sample and the design, second they did not tackle the issue of inference on a function of the parameter or prediction of a particular function of the study variable. Although it acknowledge the fact that sometimes units of the sample can be non identified and refers to \cite{ScottSmith1973}

The paper is done in a Bayesian framework.

In the introduction, Scott gives a definition of non-informative design:
\fancyquote{\citet[Sec. 1 Introduction, p.~ 1.]{Scott1977}}{It is well known (see \cite{Godambe1966a} or Ericson \cite{1969} that this posterior distribution depends only on the sample actually drawn and not on the sampling desing used to draw it, provided that [the sample of the population elements] and [the study variable on the population] are independent random variables. (If this condition of independence is satisfied the design is said to be non-informative.) \cite{Godambe1969} has called this conflict about the way inferences should depend on underlying design the problem of randomization and most of the comments in the paper are concerned with this small facet of the general problem of randomization in surveys.}

Scott gives a result that is announced in the introduction:
\fancyquote{\citet[Sec. 1 Introduction, p.~ 1-2.]{Scott1977}}{
In the next section we show that simple random sampling is the only design that is always non-informative. Other desingns are non-informative only if the information used at the design stage is also available for the analysis. This may not be the case in a secondary analysis carried out well after the collection of the original data, for example. If supplementary information used in the design is not available for the analysis, the posterior distribution depends on the design.}

\citet{Sec. 2, p. 502}, Scott considers the case where the design variables on the population are known and the sample and the study variables are independent conditionnally on the design variables, and argues that the distribution of the study variable on the population conditional to the design variable and all the information provided by the survey is proportional to the distribution of the study variable on the population conditional to the design variable and only the values of the study variable on the sample, and concludes that the design is then non-informative.

\fancyquote{\citet[Sec.2, p.~2]{Scott1977}}{Consider a situation in which information is available on some related characteristic, say $\mathbf{x}=(x_1,\ldots,x_N)$, which can be used in the design. We write $p(s\mid \mathbf{x}$ for the probability of drawing sample $s$ and $\xi(\mathbf{Y}\mid \mathbf{x})$ for the conditional prior density of $\mathbf{y}$ given $\mathbf{x}$. If $s$ and $\mathbf{y}$ are conditionally independent for the given value $\mathbf{x}$, then the posterior distribution of $\mathbf{y}$ is 
$$p(\mathbf{y}\mid (i,y_i;i\in s),\mathbf{x})\begin{array}[t]{l}=\frac{p(s\mid\mathbf{x}\xi(\mathbf{y}\mid \mathbf{x})}{\int_{Y_s} p(s\mid \mathbf{x})\xi(\mathbf{y}\mid \mathbf{x}\derive \mathbf{y}}\text{  for }\mathbf{y}\in Y_s\\=0\text{ otherwise,}\end{array}$$
where $Y_s=\{y^{\ast}\in Y:y^{\ast}_i=y_i\text{ for }i\in s\}$.
This reduces to $$p(\mathbf{y}\mid (i,y_i;i\in s),\mathbf{x})\begin{array}[t]{l}\propto\xi(\mathbf{y}\mid \mathbf{x})\text{  for }\mathbf{y}\in Y_s\\=0\text{ otherwise,}\end{array}$$

so that the design is non-informative to anyone who knows $\mathbf{x}$ and plays no part in the analysis.
}

\fancyquote{\citet[p.~3]{Scott1977}}{
Even if all supplematary information is available the design may still be informative if the sampled elements cannot be identified. The situation was discussed in detail in \cite{ScottSmith1973}}

The remaining does not contain more that could be helpful for the definition of informativeness of the paper contains a discussion about ``robust sampling strategies'', that take into account both the design and the estimation.  The paper discusses robustness follwing Blackwell and Girshick, and participates in the debate about using unequal probability sampling and refers to Basu, Royall and Ericson again.

\subsection{Rubin 1976}
\cite{Rubin1976} 
\begin{enumerate}
 \item defines what it means to ignore the missing data mechanism 
 \item define three conditions on the density of the distribution of the missing data mechanism conditional 
\item For different type of inference (a,. sampling distribution inference about $\theta$, b. direct likelihood inference, c. Bayesian inference), \cite{Rubin1976}: 
\begin{enumerate}
    \item defines what it means to say that it is appropriate to ignore the missing data mechanism, e.g the definition of what is an ignarable missing data mechanism.
    \item Give conditions for ignoring appropriately the missing data mechanism
\end{enumerate}
\end{enumerate}

Note that in this paper, the term ignorable is not used. It has been introcuced later and it is now popular in the expressions ignorable or non ignorable missing data mechanism.

In this section, we give the definitions and theorems of \cite{Rubin1976}.

\fancyquote{\citet[p.~584]{Rubin1976}}{\color{black}
 {\begin{center} 3. \textsc{Notation for the random variables}\end{center}}
 Let $U=(U_1,\ldots U_n)$ be a vector random variable with probability density function $f_\theta$. The objective is to make inferences about $\theta$, the vector parameter of this density.[...] Let $M=(M_1,\ldots M_N)$ be the associated ``missing data indicator'' vector random variable [...]. The probability tat $M$ takes the value $m=(m_1,\ldots,m_n)$ given that $U$ takes the value $u=(u_1,\ldots,u_n)$ is $g_\phi(m\mid u)$, where $\phi$ is the nuisance vector parameter of the distribution.
The conditional distribution $g_\phi$ corresponds to ``the process that causes missing data''[...].More precisely, define the extended vector random variable $V=(V_1,\ldots, V_n)$ with range extended to include the special value $*$ for missing data.``
}

The transition from \cite{Rubin1976} notations to ours is the following:
$(U\translated\Signal)$, $(n\translated N)$, $(f_\theta\translated \density_{Y})$, $(\theta\translated\param)$, $(M\translated I)$, $(\phi\translated\Comp\param)$, $(m\translated i)$, $(g_\phi(m\mid u)\translated \density_{I\mid \Signal=\signal}(i)$. In our notations, the framework chosen by \cite{Rubin1976} is the following: $\Id_{\Omega}=(\Transfo,\Signal)$, $\Signal$ follows an i.i.d population model and $\Transfo$ is a without replacement selection. The conditioning variable $\designvar$ is not needed here. We have naturally $(\paramf:P\mapsto P^\Signal,\Comp\paramf:P\mapsto P^{I\mid \Signal}$. The observation is $\Obs=(\Transfo[\Signal],I)$, which is equivalent to $\Obs=(\Transfo[\Signal],\Transfo)$ in this setup.

Implicitely, \cite{Rubin1976} assumes a parametric dominated model both for $P\mapsto P^\Signal$ and $P\mapsto P^{I\mid \Signal}$. Note that what rubins calls a ''missing data indicator`` is indeed a ''non-missing data indicator`` as $m_i=1\Leftrightarrow Y_i$ is observed.

\fancyquote{\citet[p.~584]{Rubin1976}}{
Hence the observed value of $M$, namely $\tilde m$, effects a partition of each of the vectors of random variables and the vectors of observed values into two vectors corresponding to $\tilde m_i=0$ for missing data and $\tilde m_i=1$ for observed data.
}

Rubin defines new variables, the link from Rubin to our notation is as follows:
$\left(U_{(0)}\translated (Y_k)_{k\in\{1,\ldots,N\mid I_k=0\}}\right)$, $\left(U_{(1)}\translated (Y_k)_{k\in\{1,\ldots,N\mid I_k=1\}}\right)$.

\fancyquote{\citet[p.~584]{Rubin1976}}{\color{black}
 {\begin{center} 5. \textsc{Ignoring the process that causes missing data.}\end{center}}

[...]

{\it Definition} 1. The missing data are missing at random if for each value of $\phi$, $g_\phi(\tilde m\mid \tilde u)$ takes the same value for all $u_{(0)}$,
}

In Rubin notation, $\tilde u$ is the concatenation of the unobserved $u_{(0)}$ and the observed $u_{(1)}$.\citet[Definition~1]{Rubin1976} is ambiguous for the following reasons:
\begin{enumerate}
 \item The notation $g_\phi(\tilde m\mid (u_{(0)},\tilde u_{(1)})$ can denote the value taken by $g_\phi$ for specific $\tilde m$, $u_{(0)}$ and $u_{(1)}$, but it sometimes denotes the function $g_\phi(.\mid (.,\tilde u_{(1)}):(\tilde m,u_{(0)})\mapsto g_\phi(\tilde m\mid (u_{(0)},\tilde u_{(1)})$ for a specific $(\tilde m,\tilde u_{(1)})$ and it could as well denote  $g_\phi(.\mid (u_{(0)},\tilde u_{(1)}):(\tilde m)\mapsto g_\phi(\tilde m\mid (u_{(0)},\tilde u_{(1)})$ for a specific $\tilde u_{(1)}$. Which leads to 3 possibilities.
 \item In which order should we read ''for each value of $\phi$``, ''takes the same value`` and ''for all $u_{(0)}$``? 2 consecutive forall are interchangeable, and ''for each value of $\phi$`` was stated before ''for all $u_{(0)}$`` so out of the 6 permutations, there are 4 possible interpretations.
\end{enumerate}
Combining these ambiguities lead to 12 possible interpretations.
Nevertheless, when Rubin considers whether data is missing at random for each of his examples, he clearly distinguishes cases for different values of $\tilde m$ (More precisely, for example 2, Table 1 reads: MAR only if all $\tilde m_i=1$, which means that depending on the observed value for $\tilde m$, one same missing data mechanism may be considered at random or not). Rubin uses the notation ''$\tilde~$`` to notice fixed specific values of the different parameters. The use of $\tilde m$, $\tilde u$ and $\tilde v$ shows the definition apply for a specific value of $\tilde m$ and $\tilde u_{(1)}$.
Confronting the different possible interpretations of the definition with \citet[Table 1]{Rubin1976}, one can conclude that the only right interpretation is:
[For the specific observed $(\tilde m,\tilde u_{(1)}),~\forall \phi,~ \exists C \text{ such that }\forall u_{(0)},~ g_\phi(\tilde m\mid (u_{(0)},\tilde u_{(1)})=C$],
and one can discard the other interpretations as for example:
[For the specific observed $(\tilde m,\tilde u_{(1)}),~\exists C$ such that $\forall \phi, \forall u_{(0)},~ \forall \tilde m,~g_\phi(\tilde m\mid (u_{(0)},\tilde u_{(1)})=C$], 
or
[$\forall \tilde m,~ \exists C \text{ such that }~\forall u_{(0)},~ g_\phi(\tilde m\mid (u_{(0)},\tilde u_{(1)})=C$].

%  are clearly too strong: it would mean that the distribution of the missing parameter should not depend on both the nuisance parameter $\phi$ and the non observed values for the missing mechanism to be considered at random. It becomes clear after reading \citet[Example 1]{Rubin1976} that the right interpretation is 
% Equation \eqref{eq:rubinarint3}.

This definition by Rubin has also the following drawabacks: 
\begin{enumerate*}
\item when the conditional distribution of $P^{T\mid Y=y}$ is uniquely defined, the conditional density of $T\mid Y=y$, when it is defined, is a class of density distribution that are equal $P^{T\mid Y=y}$-almost surely$(t)$, so this definition depends on a choice of representants of the density that must be made prior to its application.  
                                                             \item the domain of variation of $\{u_{(0)}$ is not specified. It is implicitely $U_{(0)}((M,U_{(1)})^{-1}(\{(\tilde m,\tilde u_{(1)})\})$.
                                                            \end{enumerate*}

 So the rigth interpretation of \citet[Definition 1]{Rubin1976} is the folling: 
 \begin{definition}[Missing at random, Rubin's definition, Rubin's notations, detailed]
The missing data are missing at random contitionnally to $M=\tilde m$ and $U_{(1)}=\tilde u_{(1)}$ if 
\begin{equation}\label{eq:rubinarint1}
\forall \phi,~ \exists C \text{ such that }\forall u_{(0)}\in U_{(0)}((M,U_{(1)})^{-1}(\{(\tilde m, \tilde  u_{(1)})\}),~ g_\phi(\tilde m\mid (u_{(0)},\tilde u_{(1)}))=C, 
\end{equation}
 \end{definition}

 Below is the translation of \citet[Definitions 1]{Rubin1976} in our notations. 
It is valid in the following framework:
Assume that $\Signal$ follows a population model of known size $N$, that $\Transfo$ is a selection, that $\Obs=\left(\Transfo[\Signal]=\left(Y_k\right)_{k\in\{1,\ldots,N\mid I_k=1\}},I\right)$, and that $\Id_\Omega=(\Transfo,\Signal)$. Consider $\paramf:P\to P^\Signal$ and $\Comp\paramf:P\mapsto P^{I\mid\Signal}$. %Note that $(\Transfo,\Signal)\isoleq(\Omega)$ so 

\begin{definition}[Missing at Random, Rubin's definition, our notations]\label{eq:MARRubinDBnotations2}
Let $\{\density_{\Transfo\mid \Signal;P}(.)\mid P\in\range{P}\}$ be a set of unique representants for each $P\in\range{P}$ of the classes of the conditional distribution densities of $\Transfo\mid\Signal$.
The missing data is missing at random with respect to $\{\density_{\Transfo\mid \Signal;P}(.)\mid P\in\range{P}\}$ and the (not necessarily distinct) separation $(\paramf:P\mapsto P^\Signal,\Comp\paramf:P\mapsto P^{\Transfo\mid\Signal})$ conditional to $\Obs=\obs$, with $\Obs=(\Transfo[\Signal],\Transfo)=(,\signal^\star,\transfo)$, if and only if any of the following equivalent conditions apply: 

\begin{equation}
\forall \Comp\param\in\Comp\paramf(\range{P}),\ \exists C\in\mathbb{R}\text{ such that }\forall \signal\in \Signal\left(\left(\Transfo[\Signal],\Transfo\right)^{-1}\left(\left(\signal^\star,\transfo\right)\right)\right), \forall P\in\Comp\paramf^{-1}(\Comp\param),\ \density_{\Transfo\mid \Signal=\signal;P}(t)=C.
\end{equation}

\begin{equation}
\left(
\begin{array}{lll}
:~\range{P}\times\left((\Signal,\Transfo)\left(\Obs^{-1}(\{\obs\})\right)\right)&\to    &,\\
~~(P,(\signal,\transfo))                                                          &\mapsto&\density_{\Transfo\mid \Signal=\signal;P}(\transfo)\end{array}\right)
\dependsonly\left(\begin{array}{lll}
                                                  :~\range{P}\times\left((\Signal,\Transfo)\left(\Obs^{-1}(\{\obs\})\right)\right)&\to    &,\\
                                                  ~~(P,(\signal,\transfo))                                                          &\mapsto& (\Comp\paramf(P),\transfo(\signal))
                                                 \end{array}
\right)\label{MAR3}\end{equation}
\end{definition}

There are many choices of generalisation that would match Rubin's definition, because Rubin is observing $\Transfo=\transfo$ and $\Transfo[\Signal]=\transfo(\signal)$. In another setup, where only $\Transfo$ is observerd, we could replace $(\Comp\paramf(P),\transfo(\signal)$ in the second part of Equation \eqref{MAR3} by $(\Comp\paramf(P),\transfo)$, or by $(\Comp\paramf(P),\transfo,\transfo(\signal))$. We choose the weakest condition that allows to get the equivalent of Rubin theorem and it that  still generalize Rubin's definition.

Which leads to the following general definition

\begin{definition}[At Random, 1st generalisation]\label{eq:MARRubinDBnotations2}
Let $\Rv$, $\Rvv$ be two random variables. Let $\paramf,\Comp\paramf$ a (non necessarily distinct) separation of $\mathscr{P}$ and let $\{\density_{\Rvv\mid \Rv;P}(.)\mid P\in\range{P}\}$ be a set of unique representants for each $P\in\range{P}$ of the classes of the conditional distribution densities of $\Rvv\mid\Rv$. Let $(\paramf,\Comp\paramf)$ be a (non necessarily distinct) separation of  $\mathscr{P}$. 
Then $\Rvv$  is at random with respect to the separation $(\paramf,\Comp\paramf)$ conditional to the observation is $\Obs=\obs$,  if and only if: 
$(P\mapsto P^\Rv)\dependsonly \paramf$, $(P\mapsto P^{\Rvv\mid\Rv})\dependsonly \Comp\paramf$, and 
\begin{equation}
\left(
\begin{array}{lll}
:~\range{P}\times\left((\Rv,\Rvv)\left(\Obs^{-1}(\{\obs\})\right)\right)&\to    &,\\
~~(P,(\rv,\rvv))                                                          &\mapsto&\density_{\Rvv\mid \Rv=\rv;P}(\rvv)\end{array}\right)
\dependsonly\left(:\begin{array}[t]{lll}    (P,(\rv,\rvv))                                                          &\mapsto& (\Comp\paramf(P),\rvv)
                                                 \end{array}
\right)\label{MAR4}\end{equation}
\end{definition}

A sufficient condition for the condition \eqref{MAR4} to hold is:
\begin{equation}
\left(
\begin{array}{lll}
:~\range{P}\times\left((\Rv,\Rvv)\left(\Omega)\right)\right)&\to    &,\\
~~(P,\rv,\rvv)                                                          &\mapsto&\density_{\Rvv\mid \Rv=\rv;P}(\rvv)\end{array}\right)
\dependsonly\left(:\begin{array}[t]{lll}    (P,\rv,\rvv)                                                          &\mapsto& (\Comp\paramf(P),\rvv)
                                                 \end{array}
\right)\label{MAR4}\end{equation},

A sufficient condition for the condition \eqref{MAR4} to hold is:
\begin{equation}\left(:(\range{P}\times\range{\Rv})\to,(P,\rv)\mapsto P^{\Rvv\mid \Rv=\rv}\right)\dependsonly\left(:(\range{P}\times\range{\Rv})\to,(P,\rv)\mapsto\Comp\paramf(P)\right)\label{MAR4}\end{equation}, 
which itself implies that 
\begin{equation}\forall P, \Rv\indep_P\Rvv\label{MAR5}\end{equation},

%A bigger problem is that the choice of the domain for $y$ is an issue: we have taken $\{y\in Y(\Omega)\mid t(y)=y^\star\}$, 
Rubin follows with the definition of observed at random:
\fancyquote{\citet[p.~584]{Rubin1976}}{
{\it Definition} 2. The observed data are observed at random if for each value of $\phi$ and $u_{(0)}$, $g_\phi(\tilde m\mid \tilde u)$ takes the same value for all $u_{(1)}$.}

The right interpretation is the following, in Rubin's notations:
\begin{equation}\label{eq:rubinarint3}\forall \phi\in\Omega_\phi,
\forall u_{(0)}\in U_{(0)}(M^{-1}(\{\tilde m\})),~ \exists C\in\mathbb{R} \text{ such that }\forall u_{(1)}\in U_{(1)}((M,U_{(0)})^{-1}(\{(\tilde m,\tilde u_{(0)})\})),~ g_\phi(\tilde m\mid (u_{(0)},u_{(1)}))=C, 
\end{equation}

Rubin definition of ''observed at random`` is less local than its definition of 
missing at random as it implies the variation of both $u_{(0)}$ and $u_{(1)}$. 

Rubin definition of missing at random is the definition of a local property for a specific $u_{(1)}$ and $\tilde{m}$. From there two versions of the definition can be proposed: a local version, for a specific value of $u_{(1)}$ and $\tilde{m}$, and a uniform version, where the local condition must hold for all possible $u_{(1)}$ and $\tilde{m}$.

The definition of observed as random is more difficult to translate in a general framework as it requires to define what we do not observe about $\Signal$, wich was not the case for the definition of observed at random.
Define $\Comp\Transfo$ as the transformation such that for all $\omega$, $\Comp\Transfo[\Signal](\omega)$ is the function $\Pop(\omega)\setminus \Sampleindexset(\omega)\to, k\mapsto Y[k](\omega)$.

\begin{definition}[Observed at Random, Rubin's definition, our notations]
The observed data are observed at random  if $\forall \Comp\param\in\Comp\paramf(\range{P})$, 
$\forall \Comp \signal\in\Comp\Transfo[\Signal](\Transfo^{-1}(\{\transfo\})$, $\exists C\in\mathbb{R}$ such that
$\forall \signal\in \Transfo[\Signal]\left((\Transfo,\Comp\Transfo[\Signal])^{-1}(\{(\transfo,\Comp\signal)\})\right)$, $\forall P\in\Comp\paramf^{-1}(\Comp\param)$,~$\density_{\Transfo\mid \Signal=\signal}(\transfo)=C.$
\end{definition}

This definition requires to define $\Comp\Transfo$. We opt for defining the complementation on $\range(\Transfo)$. 

\begin{definition}[Observed at Random, Rubin's definition, our notations, 2d version]\label{eq:MARRubinDBnotations2}
Let $\{\density_{\Transfo\mid \Signal;P}(.)\mid P\in\range{P}\}$ be a set of unique representants for each $P\in\range{P}$ of the classes of the conditional distribution densities of $\Transfo\mid\Signal$.
Let $\Comp\transfo$: $(\range(\Transfo)\to$, such that $\forall \transfo$, $\Comp\transfo$ is a (non necessarily distinct) complement of $\transfo$.
The transformed data is observed at random with respect to $\{\density_{\Transfo\mid \Signal;P}(.)\mid P\in\range{P}\}$ and the the (not necessarily distinct) separation $(\paramf:P\mapsto P^\Signal,\Comp\paramf:P\mapsto P^{\Transfo\mid\Signal})$ conditional to $\Obs=\obs$, with $\Obs=(T[Y],T)$, if and only if: 

\begin{equation}
\left(
\begin{array}{lll}
:~\range{P}\times\left((\Signal,\Transfo)\left(\Transfo^{-1}\left(\Transfo\left(\Obs^{-1}(\{\obs\})\right)\right)\right)\right)
&\to&\\
~~(P,(\signal,\transfo))
&\mapsto&
\density_{\Transfo\mid \Signal=\signal;P}(\transfo)
\end{array}\right)
\dependsonly\left(:(P,(\signal,\transfo))\mapsto (\Comp\paramf(P),\Comp\transfo(\signal))\right)\label{MAR3}\end{equation}
\end{definition}

\begin{definition}[Observed at Random, our notations, General]\label{eq:MARRubinDBnotations2}
Let $\Rv$, $\Rvv$ be two random variables. Let $\paramf,\Comp\paramf$ a (non necessarily distinct) separation of $\mathscr{P}$ and let $\{\density_{\Rvv\mid \Rv;P}(.)\mid P\in\range{P}\}$ be a set of unique representants for each $P\in\range{P}$ of the classes of the conditional distribution densities of $\Rvv\mid\Rv$. Let $(\paramf,\Comp\paramf)$ be a (non necessarily distinct) separation of  $\mathscr{P}$. And for $\rvv\in\range{\Rvv}$, let $\Comp\rvv$ denote a (non necessarily distinct) complement of $\rvv$. 
Then $\Rvv$  is observed at random with respect to the separation $(\paramf,\Comp\paramf)$ and to the complementation on $\range{\Rvv}$ conditionally on the observation is $\Obs=\obs$,  if and only if: 
$(P\mapsto P^\Rv)\dependsonly \paramf$, $(P\mapsto P^{\Rvv\mid\Rv})\dependsonly \Comp\paramf$, and 
\begin{equation}
\left(
\begin{array}{lll}
:~\range{P}\times\left((\Rv,\Rvv)\left(\Rvv^{-1}(\Rvv(\Obs^{-1}(\{\obs\})))\right)\right)&\to    &,\\
~~(P,(\rv,\rvv))                                                          &\mapsto&\density_{\Rvv\mid \Rv=\rv;P}(\rvv)\end{array}\right)
\dependsonly\left(:\begin{array}[t]{lll}    (P,(\rv,\rvv))                                                          &\mapsto& (\Comp\paramf(P),\Comp\rvv(\rv))
                                                 \end{array}
\right)\label{MAR4}\end{equation}
\end{definition}

\fancyquote{\citet[p.~585]{Rubin1976}}{\color{black}
{\it Definition} 3. The parameter $\phi$ is distinct from $\theta$ if their joint parmeter space factorises into a $\phi$-space and a $\theta$-space, and when prior distributions are specified for $\theta$ and $\phi$, if these are independent.}

\fancyquote{\citet[p.~585]{Rubin1976}}{\color{black}
 {\begin{center} 6. \textsc{Missing data and sampling distribution inference}\end{center}}
A sampling distribution inference is an inference that results solely from comparing the observed value of a statistic, e.g an estimator, test criterion or confidence interval, with the sampling distribution of that statistic under various hypothesized underlying distributions. Whithin the context of sampling distribution inference, the parameters $\theta$ and $\phi$ have fixed hypothesized values.
Ignoring the process that causes missing data when making a sampling distribution inference about the true value of $\theta$ means comparing the observed value of some vector statistic $S(\tilde v)$, equivalently 
$S(\tilde m,\tilde u_{(1)})$, to the distribution of $S(v)$ found from $f_\theta$. More precisely, the sampling distribution of $S(\tilde{v})$ ignoring the process that causes missing data is found by fixing $M$ at the observed $\tilde m$ and assuming that the sampling distribution of the observed data follows from density (5$\cdot$1) which is the marginal density of $U_{(1)}$ but from the conditional density of $U_{(1)}$ given that the random variable $\tilde m$:

\begin{equation}\int\{f_\theta(u)g_\phi(\tilde m\mid u)/k_{\theta,\phi}(\tilde m)\} du_{(0)},\tag*{(6$\cdot$1)}\end{equation}

where $k_{\theta,\phi}(\tilde m)=\int f_\theta(u)g_\phi(\tilde m \mid u) du$, which is the marginal probability that $M$ takes the value $\tilde m$. Hence the correct sampling distribution of $S(\tilde v)$ depends in general not only on the fixed hypothesized $f_\theta$ but also on the fixed hypothesized $g_\phi$.}

Ignoring a (non necessarily) observed latent random transformation $\Transfo$ when making inference on $P^{\Signal}$ consists in chosing a transformation $\transfo$ compatible with $\Obs=\obs$ and considering that the distribution of 
$\Signal$ given $\Obs=\obs$ is the distribution of $\Signal$ given $\obsf(\transfo(\Signal),\transfo)=\obs$.

\fancyquote{\citet[p.~585]{Rubin1976}}{\color{black}

\textsc{Theorem} 6$\cdot$1. {\it Suppose that (a) the missing data are missing at random and (b) the observed data are observed at random. Then the sampling distribution of $S(\tilde v)$ under $f_\theta$ ignoring the process that causes missing data, i.e. calculated from density $(5\cdot 1)$, equals the correct conditional sampling distribution of $S(\tilde v)$ given $\tilde m$ under $f_\theta g_\phi$, that is calculated from density $(6\cdot 1)$ assuming $k_{\theta,\phi}(\tilde m)>0$.}}

%\begin{theorem}
% 
% Suppose that \begin{enumerate}\item\end{enumerate}
%\end{theorem}

\fancyquote{\citet[p.~585]{Rubin1976}}{

\textsc{Theorem} 6$\cdot$2. {\it The sampling distribution of $S(\tilde v)$ under $f_\theta$ calculated by ignoring the process that causes missing data equals the correct condition sampling distribution of $S(\tilde v)$ given $\tilde m$ under $f_\theta g_\phi$ for every $S(\tilde v)$, if and only if
\begin{equation}E_{u_{(0)}}\left\{g_\phi(\tilde m \mid u)\mid \tilde m,u_{(1)},\theta,\phi\right\}=k_{\theta,\phi}(\tilde m)>0.\tag*{(6$\cdot$2)}
\end{equation}.}
[...]
The phrase 'ignoring the process that causes missing data when making sampling distribution inferences' may suggest not only calculating sampling distributions with respect to denstity (6$\cdot$1) but also interpretin g the resulting sampling distribtions as unconditional rather than conditional on $\tilde m$.

\textsc{Theorem} 6$\cdot$3. {\it The sampling distribution of $S(\tilde v)$ under $f_\theta$ calculated ignoring the process that causes missing data equals the correct unconditional sampling distribution of $S(\tilde v)$ under $f_\theta g_\phi$ for all $S(\tilde v)$ if and only if $g_\phi(\tilde m \mid u)=1$.}

}

\fancyquote{\citet[p.~586]{Rubin1976}}{\color{black}
 {\begin{center} 7. \textsc{Missing data and direct-likelihood inference}\end{center}}
 A direct-likelihood inference is an inference that results solely from ratios of the likelihood function for various values of the parameter (Edwards, 1972). Within the context of direct likelihood inference, $\theta$ and $\phi$ take values in a joint paramter space $\Omega_{\theta,\phi}$.
 Ignoring the process that causes missind data when making a direct-likelihood inference for $\theta$ means defining a parameter space for $\theta$, $\Omega_\theta$ and taking radios, for various $\theta\in \Omega_\theta$, of the marginial likelihood function based on density (5$\cdot$ 1):
 
 \begin{equation}\mathscr{L}(\theta\mid\tilde v)=\delta(\theta,\Omega_\theta)\int f_\theta(\tilde u) du_{(0)},\tag*{(7$\cdot$1)}
\end{equation}sw21  
where $\delta(a,\Omega)$ is the indicator function of $\Omega$. Likelihood  (7$\cdot$1) is regarded as a function of $\theta$ given the observed $\tilde m$ and $\tilde u_{(1)}$.

The problem with this approach is tat $M$ is a random variable whose value is also observed, so that the actual likelihood is the joint likelihood of the observed data $\tilde u_{(1)}$ and $\tilde m$:

 \begin{equation}\mathscr{L}(\theta,\phi\mid\tilde v)=\delta\{(\theta,\phi),\Omega_{\theta,\Phi}\}\int f_\theta(\tilde u) g_{\Phi}(\tilde m\mid \tilde u)du_{(0)},\tag*{(7$\cdot$2)}
\end{equation}
regarded as a function of $\theta$, $\phi$ given the observed $\tilde u_{(1)}$ and $\tilde m$.

\textsc{Theorem} 7$\cdot$1. {\it Suppose (a) that the missing data are missing at random and (b) that $\phi$ is distinct from $\theta$. Then the likelihodd ratio ignoring the process that causes missing data, that is $\mathscr{L}(\theta_1\mid\tilde v)/\mathscr{L}(\theta_2\mid\tilde v)$, equals the correct likelihood ration, that is $\mathscr{L}(\theta_1,\phi\mid \tilde v)/\mathscr{L}(\theta_2,\phi\mid \tilde v)$, for all $\phi\in\Omega_{\phi}$ such that $g_{\phi}(\tilde m\mid \tilde u)>0$.}
[...]

\textsc{Theorem} 7$\cdot$2. {\it Suppose $\mathscr{L}(\theta\mid\tilde v)>0$ for all $\theta\in\Omega_\theta$. All likelihood ratios for $\theta\in\Omega_\theta$ ignoring the causes missing data are correct for all $\phi\in \Omega_\phi$, if (a) $\Omega_{\theta,\phi}=\Omega_{\theta}\times\\Omega_{\phi}$ and (b) for each $\phi\in\Omega_\phi$, $E_{u_{(1)}}\{g_\phi(\tilde m\mid \tilde u_{(1)},\theta,\phi\}$ takes the same positive value for all $\theta\in \Omega_\theta$.}

}

In the lineage of \cite{Rubin1976}, a series of paper targeted specific issues in survey sampling that are missing data issues and proposed solutions to solve them. Those papers do not bring more concepts that are of interest for the current paper, although they were fundamental for survey statistics in general: \cite{Rubin1978} introduced multiple imputation in sample surveys, and \cite{Little1982Models} extends Rubin's results to nonresponse: Non response mechanism is seen as a missing data process, and contains more details about ways to deal with ignorable or non ignorable response, as for example the EM algorithm. It also stresses on the fact that there exist cases where the observation is not only the values of the study variables for all units in the sample.
\fancyquote{\citet{Little1982Models}}{Following
Rubin (1976, 1978, 1980), we say that the sample design is ignorable
if inferences based on the distribution $f(x_s \mid z; \theta)$ are equivalent
to inferences
based on the full distribution $f(x_s,\delta \mid z; \theta,\phi)$}

\todof{ Harel Schafer - Partial and latent ignorability in missing-data problems}

\fancyquote{\citet[Abstract]{Mealli2015Clarifying}}{
We clarify the key concept of missingness at random in incomplete data analysis. We first distinguish
between data being missing at random and the missingness mechanism being a missing-at-random one, which we call missing always at random and which is more restrictive. We further discuss how, in general, neither ofthese conditions is a statement about conditional independence.We then consider the implication ofthe more restrictive missing-always-at-random assumption when coupled with full unit-exchangeability for the matrix ofthe variables ofinterest and the missingness indicators: the conditional distribution of the missingness indicators for any variable that can have a missing value can depend only on variables that are always fully observed. We discuss implications of this for modelling missingness mechanisms.}

\subsection{The ambiguous definition of Cassel, Sarndal and Wretman, 1977, for design-based inference under the fixed population model} \citet{CasselSarndalWretman1977} propose a definition of informative selection in a framework that proves to be of little interest, because they define it in a fixed population framework and refer to the dependency between two variables that are not stochastic in their framework (the design variable and the study variable), thus necessarily independent. According to this definition, in their framework, all design is non informative. The term ''depend`` is ambiguous as it may refer to stochastic or deterministic dependence, but in none of the case it is precised in the book. 

This is the comment made in \cite{Bonnerytheseen}:

        \fancyquote{\cite{Bonnerytheseen}}{\color{black}
			Definition 1.6. of non-informative selection applied to the fixed population model
			is not consistent with the definition proposed 
			in \citet[p. 12]{CasselSarndalWretman1977}, which reads, after replacing their $\operatorname{p}$ by $\operatorname{q}$
			for consistency with our notation:			
			{\it{``A sampling design [...] $\operatorname{q}(.)$ is called a non-informative design if and only if 
			$\operatorname{q}(.)$ is a function that does not depend on 
			the $y$-values associated with the labels in [the sample] $s$[...].''}}
			This definition is ambiguous.
			The ambiguity comes from the fact that in \citet{CasselSarndalWretman1977}, the model for inference is not complete because
			the set of parameters is not given, and what is really meant by ``does not depend on'' is not explained.
			
			In \citet{CasselSarndalWretman1977},  the design measure $\operatorname{p}$ is fixed, but may 
			depend (non stochastically) on a function of the parameter $y$.
			So we consider that in \citet{CasselSarndalWretman1977}, the statistical model which is implicitly referred to is 
			parametrized by the couple $(\operatorname{p},y)$ that belongs to some subset $\Thetaset'$ of $\mathscr{S}\times  Y(\Omega)^N$.
			The statement that $\operatorname{p}$ does not depend on $(y_k)_{k\in \Pop , I_k\geq 1}$ is not ambiguous in some particular cases:
			\begin{itemize}
			\item In the case where $\Thetaset'=\left\{\operatorname{p}_0\right\}\times B$ where $\operatorname{p}_0\in\mathscr{S}$,
						$B\subseteq  Y(\Omega)^N$, then $i\mapsto\operatorname{p}(\{i\})$ is a constant function of $(y_k)_{k\in \Pop , i_k=1}$.
			\item In the case where $\Thetaset'=\left\{\left.\left(\mathrm{Poisson}_y,y\right)\right|y\in Y(\Omega)^N\right\}$, $ Y(\Omega)=[0,1]$ 
			then the definition indicates that the sample is informative, as the probability to draw a sample depends in part  on
			the values of $y$ on the sample.
			\end{itemize}
			But the definition can be  ambiguous in cases where the dependence on $y$ is not direct, 
			but is indirectly imposed by some non trivial correspondence between $\operatorname{p}$ and $y$. For example, in the case where 
			$\Thetaset'=\left\{\left.\left(\mathrm{Poisson}_z,y\right)\right|y\in[0,1]^N,z\in[0,1]^N,\|y-z\|^2\leq \frac{1}{2}\right\}$, 
			it is very difficult to apply the \citet{CasselSarndalWretman1977} definition to determine whether the sample design is informative or not.		
			
			Nevertheless, according to a certain interpretation of ``does not depend on'', the \citet{CasselSarndalWretman1977}  definition 
			can be understood as:
% 			\begin{equation}\text{there does not exist a non constant function }\varphi:\text{ such that }
% 				\forall i \in \mathbb{N}^N,\ \Design(z)(\{i\})=\varphi(\obsf_{y}(i))\label{eq:1.7}
% 			\end{equation}
% 			or 
			\begin{equation*}\text{The selection is non-informative if }
			\exists A\subseteq\mathscr{S},\  B\subseteq  Y(\Omega)^N,\text{such that }\Thetaset'=A\times B.
					\label{eq:1.8}
			\end{equation*}}
					
% 			

% This definition of \cite{Scott1977} is not exactly the same than a definition that will be used later by \cite{Pfeffermann1998a} for example, that caracterize informativeness as the case where the population distribution and the sample distribution are different.

%\subsection{} \citet{Rubin1978} Paper on multiple imputation.
%\subsection{} \citet{little1980superpopulation}  I do not have the paper.

\subsection{Sugden Smith 1984} \cite{SugdenSmith1984}, establish ``Conditions under which partially known designs can be ignored''. This paper is in the lineage of \cite{Scott1977} and \cite{Rubin1976}. The paper is particularly important for the current paper for the following reason:
following \cite{Scott1977}, it considers different cases based on what is observed: the values on the sample only, or with the design variables values on the sample, or with the design variable values on the population. This stresses the importance of having a definition that depends on what is observed, and that should not be restricted to the particular case of observing the study variables on the design only as in \cite{Rubin1976}.
The paper does not contain a definition on informative or uninformative, but it refers to \cite{Scott1977}. 

The introduction is of particular interest as it gives a short but detailed review of the litterature.
\fancyquote{\citet[Introduction,~p.~495]{SugdenSmith1984}}{\color{black}
In a model-based approach to survey sampling inference the role played by the survey design is not completely clear. Sone ahtors such as \cite{Godambe1966a}, \cite{Basu1971}, \cite{Ericson1969}, \cite{RoyallPfeffermann1982}, \cite{Little1982} and \cite{Smith1983} recognize that random sampling schems may have desirable robustness properties in a model-based approach but that other designs, such as balanced samples may be vetter for some purposes.
\cite{Scott1977} and \cite{ScottSmith1973} examine the conditions under which any survey design can be ignored for Bayesian inference. If these conditions are not satisfied the averages over subsets of the randomization distribution may be necessary for valid Bayesian inference. \cite{Rubin1976} is a fubndamental paper on missing values, interprets sampling as a special case of missing values and establishes conditions under which the selection method can be ignored for model-based inferences from the Bayesian, likelihood or sampling theory viewpoints. \cite{Little1982Models} extends Rubnin's results to nonresponse and \cite{Smith1983} to non random designs such as quota sampling which depend on response variables. The key to understanding the role of survey design is to follow \cite{Scott1977} and introduce the idea of design variables, known to the sampler before the sample is drawn, in addition to the response variables measured in the survey.}

In section 4, Sugden and Smith give examples of informative and non informative designs. The term informative and noninformative appear in the paper three times together, once in the title, and in the two quotes below:

\fancyquote{\citet[Sec.~4. Ingorable Designs, p.~501]{SugdenSmith1984}}{\color{black}
[Simple random sampling] is independent of any variable, design or response, so certainly satisfies .[\citet[Condition 1]{Scott (1977)}] regards this as the only uniformly noninformative design.
}

\fancyquote{\citet[Sec.~4. Ingorable Designs, p.~504]{SugdenSmith1984}}{\color{black}
Before embarking on an analysis of survey data collected by others an analyst must
examine his data set and his knowledge of the selection mechanism carefully to see if Condition 1 is satisfied or not. If not, then the design forms an explicit part of the model- based inference, and is 'informative' \cite{Scott1977}

Also it is important to notice that \cite{SugdenSmith1984} do not differentiate non-informative design and ignorable design.
}

\subsection{Baird 1983}
\cite{Baird1983} Fisher Pearson controversy. This is a controversy about using degrees of freedom, and informativeness of the hypothesis. MB: uses the word informative hypothesis. This paper is of little interest for us.

%\subsection{Sugden 1993} \citet{Sugden1993} This paper is about exchangeability.

\subsection{Pfeffermann et al., 1998, ''Parametric distributions of Complex Survey Data Under Informative Probability Sampling`` and subsequent in line work.} \citet{Pfeffermann1998a} as well as other papers on the same topic by Pfeffermann give a heuristic definition of informative selection: Selection is informative when it has to be taken into account. Which has the merit to be the most consensual and general definition of informative selection given. They emphasize a consequence of informative selection:
the sample and population distributions are different. They also develop tools for the analyst to make inference in presence of informative selection. The tool propose is derived from the approximation of the sample probability density function by a weighted version of the population probability distribution. The inference proposed is made in two stages: 1. the modeling of the design variable conditional to the study variable, from which can be derived an expression of the weighted distribution that will depend on this nodel parameters. 2. the estimation of both the nuisance and interest parameters. 

\fancyquote{\citet[Abstract]{Pfeffermann1998a}}{
The sample distribution is defined as the distribution of the sample mea-
surements given the selected sample. Under informative sampling, this distribution
is different from the corresponding population distribution, although for several
examples the two distributions are shown to be in the same family and only differ
in some or all the parameters. A general approach of approximating the marginal
sample distribution for a given population distribution and first order sample se-
lection probabilities is discussed and illustrated.}
% Theoretical and simulation results
% indicate that under common sampling methods of selection with unequal probabilities, when the population measurements are independently drawn from some
% distribution (superpopulation), the sample measurements are asymptotically independent as the population size increases. This asymptotic independence combined
% with the approximation of the marginal sample distribution permits the use of stan-
% dard methods such as direct likelihood inference or residual analysis for inference
% on the population distribution.
% Key words and phrases: Design variables, independence, likelihood, mixtures, PPS
% sampling, weighted distribution.
% 1. Introduction

\fancyquote{\citet[Sec. 1, Introduction]{Pfeffermann1998a}}{
Survey data may be viewed as the outcome of two random processes: The
process generating the values in the finite population, often referred to as the
‘superpopulation model’, and the process selecting the sample data from the
finite population values, known as the ‘sample selection mechanism’. Analytic
inference from survey data relates to the superpopulation model, but when the
sample selection probabilities are correlated with the values of the model response
variables even after conditioning on auxiliary variables, the sampling mechanism
becomes informative and the selection effects need to be accounted for in the
inference process.
In this article, we propose a general method of inference on the population
distribution (model) under informative sampling that consists of approximating
the parametric distribution of the sample measurements. The sample distribu-
tion is defined as the distribution of measurements corresponding to the units in
the sample.``}

Still, a general mathematical definition of informative selection is needed. And with respect to the new concepts of sample and population pdfs, as well as limit sample pdfs, more general definitions were also needed.

Following \cite{Pfeffermann1998a}, a series of paper focus on the consequence of informative selection of inducing a difference between the sample and population distribution. Below are some quotes from a selection of these papers related to the definition of informative selection.

\fancyquote{\citet{Qin2002} }{Pfeffermann, Krieger, and Rinott (1998) pointed out that
survey data can often be viewed as outcomes of two random
processes. In the first process, data are generated from a (finite
or infinite) population. However, we do not fully observe the
data in the first process. Given the data from the first process,
the second process generates a subset of data that are fully
model on the
distribution of the data and the covariates.
Often in a study, apart from the variable of interest, auxiliary information is available that can be used to increase the
precision of estimators. For example, in many surveys, the
population means (or totals) of several covariates are known,
which leads to the well-known ratio estimators and regression estimators (Cochran 1977; Rao, Kovar, and Mantel 1990;
Sarndal, Swensson, and Wretman 1992). Various other methin two-phase sampling. Otherwise, the nonresponse is 
nonignorable and the second-phase sampling is informative}

\fancyquote{\citet{Pfeffermann2002SmallAreaEstimation}}{
``If, however, the selection probabilities are related to the response variable values even after conditioning on the values of the explanatory variables included in the model, the sampling design becomes informative and the model holding for the sample data differs from the model holding in the population. Ignoring the sampling process in such cases may yield biased predictors for the target characteristics of interest. ''
and mathematic:
``Equation (2) defines the relationship between the population and sample
distributions, so that if $\pi_i$ depends on $y_i$, then $E(\pi_i|y_i) = E(\pi_i) and fp(y_i) = f_s(y_i)$. In this case the population distribution differs from the sample distribution and the sample design is informative.}

\fancyquote{\citet{Eideh2006a}}{ [...] that the selection probabilities are correlated with the variables of interest, even when conditioned on explanatory variables.}

\fancyquote{\citet{Eideh2006b}}{ When the sample selection probabilities depend on the values of the model response variable, even after conditioning on auxiliary variables, the sampling mechanism becomes informative and the selection effects need to be accounted for in the inference process.}

\fancyquote{\citet{Pfeffermann2006}}{ When the sampling probabilities are related to the values of the outcome variable after conditioning on the model covariates, the sampling process becomes informative and the model holding for the sample data is
then different from the corresponding population model before sampling}

\fancyquote{\citet{Eideh2009a}}{When the selection probabilities are related to the values of the response variable, even after conditioning on concomitant variables included in the population model, the sample design is defined as informative.}

\fancyquote{\cite{Pfeffermann2011Modelling}}{The sample selection probabilities in at least some stages of the sample selection are often unequal; when these probabilities are related to the model outcome variable, the sampling process becomes informative and the model holding for the sample is then different from the target population model.}

\subsection{In the lineage of Fuller, 2011, Sampling Statistics} 
In the book by Fuller, the definition of informative design is given  like follows:
\fancyquote{\citet[p.~ 355 sec. 6.3.2]{Fuller2009}}{If $E[x_i\pi_ie_i]\neq 0$, it is sometimes said that the design is ``informative'' for the model.``}
The quantities $x_i$, $\pi_i$ and $e_i$ are the values of the study variable, the inclusion probability and an error term of a linear model between $x_i$ and the design variable for the individual $i$ of the sample.
This definition only gives a vague idea in a very restricted context of what informative selection refers to.
Still, some authors refer to \cite{Fuller2009}, as for example \citet{Kim2013} to define informtive selection.

A series of article are in line with \cite{Fuller2009}, and define informative selection in the same very restrictive model. Quotes from a selection of such papers are provided  below:

\fancyquote{\citet{verret2010small}}{However, in many practical situations the inclusion probabilities
$\pi_{j,i}$ may be related to associated $y_{i,j}$ even when conditioning on the model covariates, $x_{ij}$. In
such cases, we have “informative sampling” and the model (1), holding for the population, no longer holds for the sample.
As a result, the estimators obtained by assuming that the 
model (1) holds for the sample may be heavily biased and their mean squared errors (MSE) significantly increased. It is therefore necessary to use methods that can account for informative sampling.}

\cite{kim2013weighting}are in line with \cite{Fuller2009} when defining informative selection. Below are two quotes from this paper:

\fancyquote{\citet[Abstract p.~386]{kim2013weighting}}{
Sampling related to the outcome variable of a regression analysis conditional on covariates is
called informative sampling.}

\citet{kim2013weighting} refer to \cite{Fuller2009}, although, how we have seen, this definition is very restrictive.

\fancyquote{\citet[Abstract p.~386]{kim2013weighting}}{
Surveyweights are often used in regression analysis of survey data to ensure consistent estimation of parameters when sampling may be informative, that is when sample inclusion may be related to the outcome variable conditional on covariates (Fuller, 2009, \S 6.3)}

The authors also give the following mathematical definition of informative sampling:

\fancyquote{\citet[Sec. 2. Basic Setup]{kim2013weighting}}{
We assume a probability sampling design, where inclusion in the sample is represented by the indicator variables $I_i(i=1,\ldots,N)$, where $I_i=1$ if unit $i$ is included in the sample and $I_i=0$ otherwise and $\pi_i=\mathrm{pr}(I_i=1\mid i)$ is the first-order inclusion probability. THen the ordinary least squares estimator of $\beta_0$ solves 
\begin{equation}\sum_{i=1}^N I_i(y_i-x'_i\beta)x_i=0 \tag{(2)}\end{equation}
for $\beta$, and this estimator will generally be biased unless sampling is noninformative, that is $I_i$ and $y_i$ are pairwise independent conditional on $x_i$,
\begin{equation}\mathrm{pr}(I_i = 1 | yi , xi ) = \mathrm{pr}(I_i = 1 | xi ). \tag{(3)}\end{equation} 
}

This definition ignores the potential issue with informative sampling which is that it may induce dependence among sample units, as for example in cluster sampling: All inclusion probabilities may be constant, conditioning on $x$ or $x$ and $y$ will not change the probability to be selected, and still, the design is clearly informative.

%\subsection{} \citet{Breidt2009}

%\subsection{Scott and Wild, On the robustness of weighted methods for fitting models to case-control data}

\subsection{Bonn\'{e}ry, Breidt and Coquet} In a series of three papers, Bonn\'{e}ry, Breidt and Coquet explored the properties of the approximated sample distribution as described by Pfeffermann. They proposed a general definition of the sample distribution that can be applied to both with and with replacement selection, and showed non parametric (convergence of kernel density estimators (resp. empirical cumulated distribution function estimators) to a limit sample distribution density (resp cumulated distribution function) and parametric results (convergence of the maximum weighted likelihood estimator.
In their papers, they emphasize that the informative selection process can produce dependence among observations, an issue that is not part of the scope of many papers about informative selection.

Following \cite{Scott1977}, the definition of informative selection used in these paper consists in saying that informative selection occurs when the distribution of what is observed is different from the distribution that would have been obtained from a simple random sampling.  Still the definition provided in these papers is not as complete as the definition of the current paper.

% \citet{Bonnery2012} Consider informative selection of a sample from a finite population. Responses are realized as independent and identically distributed (i.i.d.) random variables with a probability density function (p.d.f.) f , referred to as the superpopulationmodel. The selection is informative in the sense that the sample responses, given that they were selected, are not i.i.d. f . In general, the informative selectionmechanismmay induce dependence among the selected observations. The

\subsection{Miscellanous papers}

%\cite{Nandram2007Bayesian}

\citet{Savitsky2016Bayesian}A sampling design that produces a correlation between selection probabilities and observed values for sampled units is referred to as informative.

\cite[Section 3.5, pp.~41.3--41.4 ]{Berger1988} cites:
Kalbfleisch

\todof{Review papers of Heitjan!}

\subsection{Barnard, Jenkins
and Winsten
(1962)} \todof{read paper, give the reference described} the relationships between  parametric and sample spaces restricting
the use of the likelihood
principle.

\section{Informative designs in Geo statistics}
In geo-statistics, the term "informative design" denotes an "optimal design" in the sense that some desing may not contain any information.

\todof{Give Cressie's citation}
 
\section{On the role of the sample size, of labels, exchangeability and symmetric functions}

A series of papers on the role of labels, e.g. the information carried by labels, if they are identifiable or not (which makes it difficult to use the design variable information) seems out of scope for this paper.
Nevertheless, the notion was studied by Scott and the question is important for survey sampling, as the case of unidentifiable units often occurs. It is then important to propose a definition for informative design that applies in the case of non identifiable designs. This also explains why we took care to define the random variables $R$, the sample to population indexes mapping, that can be observable or latent. We give some notes below on these papers. The paper by \citet{ScottSmith1973} for example was written at a time of a controversy around the role of randomisation, and the role of the labels was important to discuss. It is now clear that a sample index should be non informative, it should just be a record index in a file. What is informative is the potential information related to the index, e.g. the mapping between this index and the population, index as long as the population index is linked to informative auxiliary information, or in the case of latent with replacement sampling.

\subsection{Durbin (1969), Inferential aspects of the randomness of sample size in survey sampling} 
\cite{durbin1969inferential}
is cited by Dawid 1977. It just discusses as title indicates whether one should use the information provided by the sample size, including when some testing is to be made. In Chapter 4, the author assumes that the distribution of the sample sixe depends on one or more unknown parameters  but does not depend on the parameter of interest. Last section gives a discussion (general remarks) on conditional tests or unconditional tests. The authors discuss failed attempts by Fisher to build a theory about conditional tests given the ancilliary. The author also refers to\cite{Basu1964}.

\subsection{ScottSmith1973} \citet{ScottSmith1973} Not much in this paper of interest for our definition.
\citet{ScottSmith1973} mention the fact that sometimes units on the sample are not identifiable. This is an example that can happen today when data is de-identified. This paper mentions that the problem of not identifiable sample units was also described in \cite{Godambe1966a}, \cite{Ericson1969}, \cite{HartleyRao1962}). A whole section is titled ``2. The role of labels''. This paper is not helpful for the definition of informativeness, although it is discussing how to account for the design in specific situations.
The goal of the paper is to focus on the exchangeability of the sample distribution, and to derive some properties from the exchangeability situation. This explains why they started by discussin the role of labels. There is no mention about informative or non informative design in this paper.
\newpage

\section{Old sections}

\subsection{At random transformation and non informative transformation}

For $(\Rv,\Comp\Rv)\isogeq(\Omega,\tribu_\Omega)$, define:
$$\mathrm{Atrandomize}_{\Comp\Rv/\Rv}:\range{P}\to,P\mapsto \left<\left.\Comp\rv\mapsto
 P^{(\Rv,\Comp\rv)\mid \ignoreset_{\Comp\rv,\Comp\Rv,\Rv}}\right|P^{\Comp\Rv}\right>^{\interclass_{\Rv,\Comp\Rv}},$$
 as well as

 $$\mathrm{Atrandomize}_{\Rv/\Comp\Rv}:\range{P}\to,P\mapsto \left<\left.\rv\mapsto
 P^{(\Comp\Rv,\rv)\mid \ignoreset_{\rv,\Rv,\Comp\Rv}}\right|P^{\Rv}\right>^{\interclass_{\Comp\Rv,\Rv}}.$$

If $\left(P\mapsto P^\Rv,P\mapsto(\Comp\rv\mapsto P^{(\rv,\Comp\Rv\mid\Rv=\Rv })\right)\isoeq\range{P}$ and $\range{P}'=\{P^\Comp\Rv\mid P\in\range{P}\}$  then $\range{P}^\star=\{\Atrandomize_{\Comp\Rv/\Rv}(P)\mid P\in\range{P}\}$. If $(\Rv,\Comp\Rv)\isoeq\Omega$, then 
$\mathrm{Atrandomize}_{\Rv/\Comp\Rv}(P)=\mathrm{Atrandomize}_{\Comp\Rv/\Rv}(P)=(P^\Rv\otimes P^{\Comp\Rv})^{\interclass_{\Rv,\Comp\Rv}}$.

\begin{example}
Consider two discrete variables $(\Rv,\Comp\Rv)$ with values in $\{(0,0),(0,1),(1,1)\}$ such that $(\Rv,\Comp\Rv)=\Id_\Omega$. The following tables represent the joint probability

\begin{tabular}{c|cc|cc|cc}
$P(\{\Rv,\Comp\Rv)=(\rv,\Comp\rv)\})$&$(.,0)$&$(.,1)$&$(.,0)$&$(.,1)$&$(.,0)$&$(.,1)$\\\hline
$(0,.)$&$1/6$&$1/2$&$1/3$&$1/3$&$1/3$&$1/2$\\
$(1,.)$&$1/3$&$0$&$1/3$&$1/2$&$1/6$&$0$\\
\end{tabular}

\end{example}

\begin{definition}
$\Comp\Rv$ is at random with respect to $\Rv$ and $\range{P}$ if \begin{enumerate*}
\item $(\Rv,\Comp\Rv)\isogeq (\Omega,\tribu_\Omega)$,
\item $\forall P\in\range{P}$, $\mathrm{Atrandomize}_{\Comp\Rv/\Rv}(P)=P$, and
\item $(P\mapsto P^{\Rv},P\mapsto (\rv\mapsto P^{\Comp\Rv\mid\Rv=\rv}))\isoeq\range{P}$.
\end{enumerate*}
\end{definition}

%\begin{example}
%Consider a population model, with a selection of 3 units with replacement and a population of known size 10.
%Let $\Rv=\Signal$, $\Comp\Rv=\Transfo$.
%For $\rv=(1,1,1)$, $\ignoreset_{\rv,\Rv,\Comp\Rv}=\Omega$, whereas for $\rv=(1,2,3)$, 
%$\ignoreset_{\rv,\Rv,\Comp\Rv}=\{\omega\in\Omega\mid \Samplepopmap(\omega) \text{ is injective}\}$.
%
%The observation is $\Obs=\Transfo[\Signal]$.
%Condition
%
%\end{example}

%\begin{example}
%Consider a population of rabbits on an island with a limited amount of resources $\somesubset$.
%Let $\Rv$ be the average amount of grass rabbits eat in a particular day.
%Let $\Comp\Rv=(\Comp\Rv_1,\Comp\Rv_2)$ be the amount of resources and the number of rabbits.
%The support of $\Rv$ conditionnally to $\Comp\Rv=\Comp\rv$ is $\ignoreset_{\Comp\rv,\Comp\Rv,\Rv}=[0, \Comp\rv_1/\Comp\rv_2]$. The support of $\Comp\Rv$ conditionnally to $\Rv=\rv$ is $\ignoreset_{\rv,\Rv,\Comp\Rv}=\{(\Comp\rv_1,\Comp\rv_2)\mid \Comp\rv_1\geq\rv~\Comp\rv_2)]$.
%
%\end{example}

\subsection{Sufficient conditions for ignorability for likelihood based inference}
Consider $\paramf:\range{P}\to$ such that $\paramf\dependsonly P\mapsto P^\Rv$ and such that $\exists \Comp\paramf$ for which $\paramf,\Comp\paramf\isoeq \range{P}$.
Then the random variable $\Comp\Rv$ is at random for the target $(\paramf,\Comp\paramf)$ by opposition to $\Rv$ given $\Obs$ is observed if 
\begin{enumerate}
\item $\left(P\mapsto P^\Rv,P\mapsto(\Comp\rv\mapsto P^{(\rv,\Comp\Rv\mid\Rv=\Rv })\right)\isoeq\range{P}$
\item $(P\mapsto P^\Rv,$
\item
\item 
\end{enumerate}

$\exists$ a $(\paramf,\Comp\paramf)$-specific sufficient statistic $S$ in the model $\mathscr{P}$ and conditionnally on $\Obs=\obs$ if:
\begin{equation}\label{eq:atrandomV}\forall P\in\range{P},~\exists \sufstat_0\in \Sufstat(\Omega)\text{ such that }P^{S\mid X=x}\text{-a.s.}(s), P^{\Comp\Rv\mid \Sufstat=\sufstat}=\left(\Atrandomize_{\Rv,\Comp\Rv}(P)\right)^{\Comp\Rv\mid \Sufstat=\sufstat_0}\end{equation}

\begin{definition}[``At random'']\label{def:atrandom}
Let $(\paramf,\Comp\paramf)\isogeq \range{P}$, $(\Rv,\Comp\Rv)\isogeq(\Omega,\tribu_\Omega)$.
The random variable $\Comp\Rv$ is at random for the target $(\paramf,\Comp\paramf)$ by opposition to $\Rv$ given $\Obs$ is observed if $\exists$ a $(\paramf,\Comp\paramf)$-specific sufficient statistic $S$ in the model $\mathscr{P}$ and conditionnally on $\Obs=\obs$ if:
\begin{equation}\label{eq:atrandomV}\forall P\in\range{P},~\exists \sufstat_0\in \Sufstat(\Omega)\text{ such that }P^{S\mid X=x}\text{-a.s.}(s), P^{\Comp\Rv\mid \Sufstat=\sufstat}=\left(\mathrm{Atrandomize}_{\Rv,\Comp\Rv}(P)\right)^{\Comp\Rv\mid \Sufstat=\sufstat_0}\end{equation}

The random variable $\Rv$ is at random for the target $(\paramf,\Comp\paramf)$ by opposition to $\Comp\Rv$ given $\Obs$ is observed if $\exists$ a $(\paramf,\Comp\paramf)$-specific sufficient statistic $S$ in the model $\mathscr{P}$ such that \eqref{eq:atrandomV} holds $P^{\Obs}$-a.s.$(\obs)$

\end{definition}

%
%If $\image(\range{P}\to,P\mapsto P^Y)=\image(\range{P}\to,P\mapsto (\ignore_{\Rv,\Comp\Rv}(P))^Y)$, one can define:
%$\paramf_{i}:\ignore_{\Rv,\Comp\Rv}\range{P}\to,P\mapsto P^Y),$
%and ``apply'' $\paramf$ to $\mathrm{Ignore}^{\Obs=\obs}_{\Rv,\Comp\Rv}(P)$ without difficulty as long as $\mathrm{Ignore}^{\Obs=\obs}_{\Rv,\Comp\Rv}$
%Then $\mathrm{ignore}_{\Rv,\Comp\Rv}$
%
%
%
%In the case where $\range{P}$ is dominated, then for a coherent set of density probabilities 
%$\{f_{.\mid.;P}\mid P\in\range{P}\}$, we have
%
%
%
%$\derive P^{\Comp\Rv}\restrict_{\ignoreset_{\rv,\Rv,\Comp\Rv}}$

In general, $P^{(\Rv,\Comp\Rv)}=P^{(\Rv,\Comp\Rv)\mid\Comp\Rv=\rv}.P^\Rv$. 
Considering that $\Rv$ is ``at random'' consists in considering that $P^{(\Rv,\Comp\Rv)\mid\Comp\Rv=\Comp\rv}=P^{\Rv}\otimes\Dirac_{\Comp\rv} $, so that 
$P^{(\Rv,\Comp\Rv)}=P^{(\Rv,\Comp\Rv\mid\Comp\Rv=\rv)}.P^\Rv$. In order to create the ``at random version'' of a model, one needs to define the function
 $$\mathrm{Atrandomize}_{\Comp\Rv/\Rv}:\range{P}\to,P\mapsto \left< \left.\Comp\rv\mapsto
 P^{\Rv\mid \ignoreset_{\Comp\rv,\Comp\Rv,\Rv}}(A)\right|P^{\Comp\Rv}\right>^{\interclass_{\Comp\Rv,\Rv}}$$

Assuming that $P$ is at random for $\Rv$, consists in assuming that $\forall P\in\range{P}$,

  Ignoring the random process that generates $\Rv$ consists in replacing the model 
  $\mathscr{P}$ by $\bigcup_{P\in \range{P}}\mathrm{ignore}_{\Rv,\Comp\Rv}(P))$.
Note that this new model is now parametrized by $\range{P}\times\Rv(\Omega)$, although this parametrization may not be bijective (two parameters may point to the same element).

\permanentcomments{This is what I had in mind before
The random variable $\Rv$ is {\em \bf at ramdom} with respect to $\paramf(P)$ provided we conditionnally on $\Obs=\obs$  if there exists at least one sufficient statistic $S$  with respect to $\paramf(P)$ such that

\begin{equation}\label{cond:atrandom}
\forall P\in \range{P},\ \exists \sufstat_0\in\sufstatf(\range{\Obs})\text{ such that }P^{\Sufstat\mid \Obs=\obs}-\text{a.s.}(\sufstat),\ \density_{\Rv\mid \Sufstat=\sufstat}(\rv)=\density_{\Rv\mid \Sufstat=\sufstat_0}(\rv). 
\end{equation}

The random variable $\Rv$ is {uniformly \em \bf at ramdom} with respect to $\paramf(P)$  given $\Obs$ is observed if there exists at least one sufficient statistic $S$  with respect to $\paramf(P)$  such that 

\begin{equation}\label{cond:atrandom2}
\forall P\in \range{P},\ \exists \sufstat_0\in\sufstatf(\range{\Obs})\text{ such that }P^{\Sufstat\mid \Obs=\obs}-\text{a.s.}(\sufstat),
\  P^{\Rv\mid \Sufstat=\sufstat}= P^{\Rv\mid \Sufstat=\sufstat_0}
\end{equation}

}

\begin{example}[Transformation at random]\label{def:atrandom}
Assume a model where $Y$ is minimal sufficient in $\range{P}$ for $\paramf:P\mapsto P^Y$, and that one observes $\Obs=\Transfo[Y]$, then the transformation $\Transfo$ is said {\textbf {at random}} conditionnally on $\Obs=\obs$ with respect to $\paramf$ if:
\begin{equation}\label{cond:atrandom}
\forall P\in \mathscr{P},\ \exists y^\star\in\range{Y}\text{ such that } P^{Y\mid \Transfo[Y]=\signal^\star}-\text{a.s.}(\signal),\ \density_{\Transfo\mid \Signal=\Signal}(\transfo)=\density_{\Transfo\mid \Sufstat=\sufstat_0;\xiparam}(\transfo). 
\end{equation}

\todof{Make sure one can write $\density_{\Transfo\mid \Sufstat=\sufstat;\xiparam}(\transfo)$ Make it a property or an assupmtion. It should just come from the fact that when one conditions by a sufficient stat, the result is ancillary.

$P_{\theta,\xiparam}^{\Transfo\mid \Sufstat=\sufstat}=P_{\theta,\xiparam}^{\Transfo\mid \sufstatf(Y)=\sufstat}$
}
When the observation is $\Obs=\obsf(\Transfo[Y],T,Z)$, the transformation $\Transfo$ is said (uniformly) at random if condition \eqref{cond:atrandom} holds $ P^{\Obs}$-a.s.$(\obs)$,e.g. if:
\begin{equation}\label{cond:atrandom2}
\forall (\theta,\xiparam)\in \Gamma,\  P_{\theta,\xiparam}^{\Obs}-\text{a.s.}(\obs),\ \exists \sufstat_0\in\range{Y}\text{ such that } P_{\theta,\xiparam}^{\Sufstat\mid \Obs=\obs}-\text{a.s.}(\sufstat),\  P_{\theta,\xiparam}^{\Transfo\mid \Sufstat=\sufstat}= P_{\theta,\xiparam}^{\Transfo\mid \Sufstat=\sufstat_0}. 
\end{equation}

\end{example}

\comments{JS 20190204:Rubin’s I is for missing values, and this should be pointed out.

Your notation is, in another way, different from Rubin’s. That is, Rubin has u0 and u1. Is this clearer than having only u and u1, even though the two choices are equivalent? I have used Rubin’s notation, realizing that you have equivalent ones.

DB 2019 02 21

You are right, I have not said what $u_0$ was in Rubin, and I use it. I need to specify that. 
}

\todof{Explain Rubin $u_0$}

\begin{example}[Missing at random]\label{ex:mar}
\citet[Definition~1, p.584 ]{Rubin1976}, given a selection without replacement from a superpopulation model, 
defines missing data  missing at random [conditionally to the observation of $I=i$ and $\Transfo[Y]=y^\star$]
as $\forall (\theta,\xiparam)\in \Gamma$, $\exists y_0\in\range{Y}$ such that $ P_{\theta,\xiparam}^{Y\mid \Transfo[Y]=y^\star}-\text{a.s.}(y),\ \density_{1-I\mid Y=y;\theta,\xiparam}(1-i)=\density_{1-I\mid Y=y_0;\theta,\xiparam}(1-i)$.
(Symbols $\phi$, $\tilde u$, $u_{(1)}$, $\tilde m$ in \citet{Rubin1976} notations corresponds to $\xiparam$, $y$, $y^\star$, $1-I$ in our notations).
 Note that \citet{Rubin1976} definition is slightly different as he defines $\phi$ as the parameter of $P^{I\mid Y}$, so in our notations, $\phi$ is a function of $(\theta,\xiparam)$.
Note also that Rubin Definition 1 is ambiguous:
{\em ``The missing data are missing at random if for each value of $\phi$, $g_\phi(\tilde m\mid \tilde u)$ takes the same values for all $u_{(0)}$.''}
can be translated into:
$\forall (\theta,\xiparam)\in \Gamma$, $\exists C$ such that
$ P_{\theta,\xiparam}^{Y\mid \Transfo[Y]=y^\star}a.s(y),\ \expected_{\theta,\xiparam}[1-I\mid Y=y]=C$.
but may also mean:
$\exists C$ such that $\forall (\theta,\xiparam)\in \Gamma$, $ P^{Y\mid \Transfo[Y]=y^\star}a.s(y) E_{\theta,\xiparam}[1-I\mid Y=y]=C$, 
which is different. Example 1 of \citet{Rubin1976} though shows that the first interpretation is the right one. 
\end{example}

\cite{Rubin1976} gives a definition of what ignoring the missing data mechanism means under a model dominated by a $\sigma$-finite measure, and the definition of ignorable missing data mechanism may be deduced 
from \cite[Th. 6.2, ]{Rubin1976}: first \cite{Rubin1976} defines what is ignoring the missing data mechanism, then this theorem gives sufficient conditions on the missing data mechanism for the inference when ignoring the missing data mechanism to be equivalent to the inference when not ignoring it: we deduce it corresponds to the definition of an ignorable mechanism. We propose a more general defintion that applies in our general framework:

\todof{I need to check this remark, Domination is there to define densities wrt a dominated measure. Here we could maybe follow Sverdrup and work on the density of $P^Y\mid S$ wrt $P^Y$}

\begin{remark}
 Note that the definition of ``at random'' conditionnally to $T=t$ (Definition \ref{def:atrandom}, condition \eqref{cond:atrandom}) requires to use densities, which requires that the model is dominated.
\end{remark}

\todof{I need to put more examples here.}

\todof{I need to illustrate MAR and MCAR.}

\subsubsection{R-ignorable random variable with respect to $h(P)$}
\begin{definition}[R-Ignorable transformation]\label{def:ignorabletransformation}
The random variable $\Rv$ is {\em \bf uniformly R-ignorable} with respect to $\paramf(P)$ given $\Obs$ is observed if there exist a sufficient statistic $S=s(\Obs)$ for $h(P)$ derived from $\Obs$ and a distinct complement $\bar\Rv$ of $\Rv$ such that 
\begin{equation}\label{con:ignorable}
 \forall P\in \range{P}, P^\Rv-a.s.(\rv),\  P^{S\mid \Rv=\rv}= P^{\Sufstat[\biclass_\Rv(\rv,\Comp\Rv)]}.
\end{equation}
The random variable $\Rv$ is {\em \bf R-ignorable} with respect to $\paramf(P)$ conditionnaly on $\Obs=\obs$ if there exists a sufficient statistic $S=s(\Obs)$ for $h(P)$ derived from $\Obs$ and a distinct complement $\bar\Rv$ of $\Rv$ such that 
\begin{equation}\label{con:ignorable2}
 \forall P\in \range{P},\ \density_{\Sufstat\mid \Rv=\rv}(\sufstat)=\density_{\Sufstat[\biclass_\Rv(\rv,\Comp\Rv)]}(\sufstat).
\end{equation}
{\textbf{Ignoring}} the random variable $\Rv$ consists in doing as if  $ P^{\Rv}-\text{a.s.}(\rv),\  P^{S\mid \Rv=\rv}= P^{\Sufstat[\biclass_\Rv(\rv,\Comp\Rv)]}$ when doing the inference.
Which is equivalent to condider that 
\end{definition}
\begin{remark}
In  ``R-ignoraribility'', ``R'' stands for ``Rubin''. This choice was made to differentiate the concepts of ignorability introduced by \cite{Rubin1976} and the general concept of ignorability.
\end{remark}

\begin{definition}[Noninformative random variable]\label{def:1}

The random variable $\Rv$ is said {\em \bf non-informative} contitional to $\Obs=\obs$ with respect to $\paramf(P)$ if   and 
only if $\exists \Comp\paramf,\ \Comp\Rv$  distinct complements of $\paramf$ and $\Rv$ respectively such that
\begin{enumerate}
 \item $\Rv$ is R-ignorable with respect to $\paramf(P)$ conditionnaly on  $\Obs=\obs$  ({\bf R-Ignorability}),
 \item $P^{\Rv\mid \obsf[\biclass_{\Rv}(\rv,\Comp\Rv)]=\obs}$  does not depend on $\paramf(P)$({\bf + Ancillarity}).
\end{enumerate}
The assertion $P^{\Rv}$  does not depend on $\paramf(P)$ is equivalent to $\exists \Comp\paramf$ a distinct complement of $\paramf$ such that 
$\forall P,P'\in\range{P},\ \Comp\paramf(P)=\Comp\paramf(P')\Rightarrow P^{\Rv}=(P')^{\Rv}$.

The random variable  $\Rv$ is said (uniformly) {\em \bf non-informative} given $\Obs$ is observed with respect to $\paramf(P)$ if for the same distinct complement $\Comp\Rv$ of $\Rv$.
\begin{enumerate}
 \item $\Rv$ is (uniformly) R-ignorable with respect to $\paramf(P)$  ({\bf R-Ignorability})
 \item and if $P^{\Rv}$  does not depend on $\paramf(P)$ ({\bf + Ancillarity}).
\end{enumerate}
The assertion $P^{\Rv}$  does not depend on $\paramf(P)$ is equivalent to $\exists \Comp\paramf$ a distinct complement of $\paramf$ such that 
$\forall P,P'\in\range{P},\ \Comp\paramf(P)=\Comp\paramf(P')\Rightarrow P^{\Rv}=(P')^{\Rv}$.

\end{definition}

\subsection{Definitions of at random, non informative and R-ignorable transformation}\label{sec:definition}

\subsubsection{Ignorable Transformation}

\begin{definition}[Ignorable transformation]\label{def:ignorabletransformation}
When the observation is the couple $(\Transfo[Y],T)$,  the transformation $T$ is said (uniformly) {\em \bf ignorable} with respect to $\theta$ if and only if
\begin{equation}\label{con:ignorable}
 \forall \theta, \xiparam \in\Gamma, P_{\theta,\xiparam}^T-a.s.(t),\  P_{\theta,\xiparam}^{\Transfo[Y]\mid T=t}= P_{\theta}^{t(Y)}.
\end{equation}
When the observation is the couple $(\Transfo[Y],T)$, the transformation $T$ is said {\em \bf ignorable} conditionnaly on $\Transfo[Y]=y^\star$ and $T=t$ with respect to $\theta$ if 
\begin{equation}\label{con:ignorable2}
 \forall \theta, \xiparam \in\Gamma,\ \density_{\Transfo[Y]\mid T=t;\theta,\xiparam}(y^\star)=\density_{t(Y);\theta}(y^\star).
\end{equation}
{\textbf{Ignoring} the transformation} consists in doing as if $ P_{\theta,\xiparam}^{\Transfo[Y]\mid T=t}= P_{\theta}^{t(Y)}$ in the inference.
\end{definition}

\comments{JS 20190204: You refer to example 21 long before it is presented. This isn’t good, but perhaps unavoidable.

DB 20190224 - You are right, I need to fix this}

\begin{remark}[Ignorability]
 The choice of using the term ignorable for defining a transformation that have the properties of Defintion $\ref{def:ignorabletransformation}$ was made to be consistent with \citet{Rubin1976} definition of ignorable non response.
 Note that very regrettably ``$T$ is an ignorable transformation'' in this sense does not mean that ``$T$ is ignorable'' in the original sense of ignorability in statistics.
 An ignorable statistic is any statistic that is independent on at least one sufficient statistic 
 (\citet[p.~142, Exercise 33]{Schervish1995}). Example \ref{ex:ignorableandnonignorable} shows a case where  $T$ can be a non ignorable selection in the sense of Definition \ref{def:ignorabletransformation} and an ignorable statistic in the classical sense: the use of the term ignorable by \citet{Rubin1976} conflicts with the more general meaning of ignorable in statistics.
\end{remark}

\begin{example}[Selecting the maximum]\label{ex: selectmax}
Suppose $Y_k\sim\mathrm{Uniform}([0,\theta])$, and our sample consists of the maximum observed value of $Y$. Here, $R[1]=\argmax_{k\in\Pop}Y_k$, $\Transfo[Y]=Y_{R[1]}=\max_k(Y)$, and $T$ is the function $T:y\mapsto y_{R[1]}$.
$Y$ is not independent of $T$, but $\Transfo[Y]$ and $T$ are independent. The selection is not ignorable as $ P^T$-a.s.(t), $ P^{t(Y)}=\mathrm{Uniform}([0,\theta])\neq P^{\Transfo[Y]}= P^{\max(Y)}$.
\end{example}

\begin{remark}
In Example \ref{ex: selectmax}, $\Transfo[Y]$ and $T$ are independent, and the selection is not ignorable, it just ensures that $ P^T-a.s(t)$, $ P_{\theta,\xiparam}^{\Transfo[Y]\mid T=t}= P_{\theta}^{\Transfo[Y]}$. 
 \end{remark}
\begin{proof}Let $t\in T(\Omega)$, $(\theta, \xiparam)\in\Gamma$. 
 The following is always true: $ P_{\theta,\xiparam}^{\Transfo[Y]\mid T=t}= P_{\theta,\xiparam}^{t(Y)\mid T=t}$. Independence implies that $ P_{\theta,\xiparam}^{\Transfo[Y]\mid T=t}= P_{\theta,\xiparam}^{\Transfo[Y]}$.
 For all $P^{\Transfo[Y]}$-measurable set $A$, We have $ P_{\theta,\xiparam}^{\Transfo[Y]}(A)=\int  P_{\theta,\xiparam}^{t(Y)\mid T=t}(A)\derive P^T(t)$.
 As $T$ and $Y$ are independent, $ P_{\theta,\xiparam}^{Y\mid T=t}= P_{\theta,\xiparam}^{Y}= P_{\theta}^{Y}$. 
 \end{proof}

\begin{remark}
Example \ref{ex: selectmax} is an extreme case of independence of $\Transfo[Y]$ and $T$ were $ P^T-a.s(t), \Transfo[Y]=t(Y)$. A less extreme case consists in taking the sample consisting of the units with the largest values for the design variable $Z$, that we assume positively correlated with $Y$. Once again, $T$ and $\Transfo[Y]$ are independent, and the selection non ignorable. 
 \end{remark}
 
\begin{remark}
 In Example \ref{ex: selectmax}, one could have defined $T$ differently, as the function $T:y\mapsto \max(y)$, but in this case, $T$ is not a selection, according to our definition of a selection, and then we cannot talk about ignorability of the selection. In this case the transformation is ignorable:
 with probability $1$, $T$ is the function max, and $ P^{\Transfo[Y]\mid T=\max}= P^{\max(Y)}$.
 This example shows the subtility of the defintion of ignorable selection and the necessity to first define selection.
\end{remark}

\begin{property}[Ignorability and likelihood]
Assume ignorable selection, the full likelihood of $\theta$ is 
$$\likelihood((\theta,\xiparam)\mapsto\theta;G)(\theta;(y^\star,t))=\sup_\xiparam\{\mathscr{L}((\theta,\xiparam)\mapsto(\theta,\xiparam),t(Y))(\theta,\xiparam;y^\star)\times \mathscr{L}((\theta,\xiparam)\mapsto(\theta,\xiparam);T)(\theta,\xiparam;t)\}.$$
When the transformation is ignorable, it does not mean that ignoring the transformation will result in the same analysis than when not ignoring it: it consists in ignoring the information brought by $T$. The inference is still valid, but one does not use all the information available.
Ignoring the selection consists in replacing $\mathscr{L}((\theta,\xiparam)\mapsto\theta;G)(\theta;(y^\star,t))$ by $\mathscr{L}(\theta\mapsto\theta;t(Y))(\theta;(y^\star))$.
In addition to the ignorabilty condition, a sufficient condition for the two likelihoods to be equivalent is that 
$\mathscr{L}((\theta,\xiparam)\mapsto(\theta,\xiparam);T)(\theta,\xiparam;t)$ does not depend on $\theta$
\end{property}

\subsubsection{Non Informative Transformation}

The notions of at random and ignorable tell us when a transformation {\em can }be ignored. 
To propose a definition of non informative transformation, we must make it consistent with the concept of information of Fisher. It is possible that the transformation can be ignored: the analysis will still be valid, but should not, as by ignoring the transformation, we are discarding some contribution to the likelihood of the transformation. This is also the spirit of the heuristic definition given by \citet{Pfeffermann1998a}: selection is informative when it {\em{has to} be accounted for in the inference process}. Following \citet{Scott1977}, another hint for the definition would be the idea that 
selection is non informative if the ``sample distribution'' \comments{JS 20190204:  “is the same” is not fully defined.

DB 2019 02 24: I am working on that}is the same than the ``population distribution'', which lead to conclude, in a specific framework that simple random sampling was the only uniformly non informative design. This idea was used by \citet{Bonnerytheseen}, and \citet{Bonnery2012} to propose a definition that consisted in saying that the selection was non informative if the likelihood of the observation was the same than the one that would have been obtained with a simple random sampling. The following definition conciliates those approaches.

\comments{JS 20190204: My confusion in the notation may be seen below where T = … includes both parts of the composition function; see Definition 3. It will be much clearer if the two parts of the definition of T are given separately, i.e., the mapping of sample indices and that of Y.}

%Due to the definition of $\theta$, we remark that different possible writings of $ P_\theta^{t(Y)}$ are  $\left( P_{\theta}^{Y}\right)^t$  and $\left( P_{\theta,\xiparam}^{Y}\right)^t$.

\begin{definition}[Noninformative transformation]\label{def:1}
The transformation $T$ is said (uniformly) {\em \bf non-informative} given $\Transfo[Y]$ and $T$ are observed for estimation of $\theta$ if 
\begin{enumerate}\item $T$ is (uniformly) ignorable for estimation of $\theta$  ({\bf Ignorability})
 \item and if $ P_{\theta,\xiparam}^T$  does not depend on $\theta$ ({\bf + T contains no information}).
\end{enumerate}
The transformation $T$ is said {\em \bf non-informative} contitional to $\Transfo[Y]=t(y)$ and $T=t$ with respect to $\theta$ if $T$ is  ignorable with respect to $\theta$ conditionnaly on  $\Transfo[Y]=y^\star$  and 
if $ P_{\theta,\xiparam}^{T\mid \Transfo[Y]=y^\star}$  does not depend on $\theta$ (e.g, $\forall (\theta_0,\xiparam_0)\in \Gamma, \left[(\theta,\xiparam_0)\in\Gamma\Rightarrow 
 P_{\theta_0,\xiparam_0}^{T\mid \Transfo[Y]=y^\star}= P_{\theta,\xiparam_0}^{T\mid \Transfo[Y]=y^\star}\right]$) 
\end{definition}

Those two definitions allow to answer different questions: ``can we ignore the information given by the transformation ?'' and ``should we ignore the information given by the transformation''?
In ignorable transformation, one can ignore the transformation, and  still, the inference will be valid in some sense. In non informative transformation, the transformation is ignorable in the sense of \citet{Rubin1976} and also ancillary. The example below shows a transformation that is both ignorable and informative.

How did we decide to combine the two conditions of ignorability, and absence of information on $\theta$ in the transformation to produce the defintion of informative selection ? We followed the heuristic of existing papers. As discussed in the introduction, when defining many authors focus on the differnce between sample and population distributions to define informative selection. Other authors will prefer to focus on the fact that the inclusion probabilities, even after conditioning on some auxiliary variables, still contain information on the distribution parameters, via their dependence to the study variable. This definition seems to be the one that will get the largest consensus.

\begin{example}[An ignorable transformation can be informative]\label{ex:ignorableandinformative}
 Assume that $Y\sim \mathrm{Bernoulli}(\theta)^{\otimes 2}$, and $T$ is the transformation associated with simple random sampling of size 1 between 2 with probability $\theta$, and the idendity with probability $1-\theta$:
 $$T=\left|\begin{array}{lll}:(y_1,y_2)\mapsto y_1 &\text{with probability}& \theta/2\\
            :(y_1,y_2)\mapsto y_2 &\text{with probability}& \theta/2\\
            :(y_1,y_2)\mapsto (y_1,y_2) &\text{with probability}&1-\theta 
\end{array}\right.$$
 Then the transformation is ignorable and informative:
 for $i \in \{1,2\}$, $ P_\theta^{\Transfo[Y]\mid T=:y\mapsto y_i}= P^{Y_i}$, and $P_\theta^{\Transfo[Y]\mid T=:y\mapsto (y_1,y_2)}= P^{Y_1,Y_2}$. In terms of likelihood,
 $L(\Transfo[Y]=y; T=:y\mapsto y_1)=(\theta/2)~\density_{Y_1;\theta}(y)$: the observation of the selection is informative.
 
 In this example, the parameter space is not separated: $\Gamma\neq\Thetaset\times\Xiset$.  Indeed, the model can be re-written
 $Y\sim \mathrm{Bernoulli}(\theta)^{\otimes 2}$,
 $$T=\left|\begin{array}{lll}:(y_1,y_2)\mapsto y_1 &\text{with probability}& \xiparam/2\\
            :(y_1,y_2)\mapsto y_2 &\text{with probability}& \xiparam/2\\
            :(y_1,y_2)\mapsto (y_1,y_2) &\text{with probability}&1-\xiparam 
\end{array}\right.,$$
and $\Gamma=\{\theta,\xiparam\in \Thetaset\times\Xiset\mid \theta=\xiparam\}$.
 Then the transformation is ignorable and informative:
 for $i \in \{1,2\}$, $ P_\theta^{\Transfo[Y]\mid T=:y\mapsto y_i}= P^{Y_i}$, and $P_\theta^{\Transfo[Y]\mid T=:y\mapsto (y_1,y_2)}= P^{Y_1,Y_2}$. In terms of likelihood,
 $L(\Transfo[Y]=y; T=:y\mapsto y_1)=(\theta/2)~\density_{Y_1;\theta}(y)$: the observation of the selection is informative.

 This situation, because it has no practical interest, was not considered in 
 \cite[Theorem 7.1]{Rubin1976}, that assumes that $[\theta]$ is distinct from $[\xiparam]$ (e.g. $\Gamma=\Thetaset\times\Xiset$) and states that missing at random ensures that the likelihood ratio obtained after ignoring missing data equals the correct likelihood ratio. In practice, $\xiparam$ and $\theta$ are distinct. But here the idea is to propose a general mathematical definition and to explore its aspects. 
\end{example}

%\comments{JS 20190204: I have struggled with the notation here. It will be clearer if you use the notation in Definition 2.3. As an example, R is a function from L(\omega) to U(\omega), i.e., it maps indices. So, I don’t understand the definition that you have used for R = (arg min Yk | k  …, arg max Yk | k ….). Although I got stuck earlier, I doubt that it’s easy to see the explanation in the last two lines.}

\begin{example}[Where $T$ is an informative non ignorable transformation, and $T$ is an ignorable statistic]\label{ex:ignorableandnonignorable}
 Assume that $Y\sim \mathrm{Normal}(0,\sigma^2)^{\otimes N}$, and $T$ is the selection associated with $R=(\argmin(Y_k\mid k\in\{1,\ldots,N\}),\argmax(Y_k\mid k\in\{1,\ldots,N\}))$.
 Then here $\Transfo[Y]=(\min(Y),\max(Y))$ and $\Transfo[Y]$ is independent of $T$. Given $R=(1,5)$, $T$ is the function $t:y\mapsto (y_1,y_5)$, and $ P^{t(Y)}= P^{(Y_1,Y_5)}\neq  P^{\Transfo[Y]\mid T=t}= P^{\min(Y),\max(Y)}$.
 So the selection is non ignorable, and  informative, although $ P_{\theta,\xiparam}^T$ does not depend on $\theta$, as unconditionnally on $Y$, $T$ is a simple random sampling, so the variable $T$ is ancillary, and ignorable as it is independent on the sufficient statistic $\Transfo[Y]$.
\end{example}

\begin{remark}
% In the existing litterature, the definition of informative selection as selections that create a difference between the population distribution and the sample distribution, but also as selections that contain information on the population parameters. The difference between sample and population distribution is given by the ignorabilty or non ignorability of the selection. To take into account the fact that the selection bears information, we need to differntiate ignorability and noninformativeness. 
Example \ref{ex:ignorableandnonignorable} shows that informative selection does not necessarily mean that 
the selection $T$ contains information on $\theta$. Non ignorability is sufficient to get informativeness.
\end{remark}

\begin{remark}[Non informative selection and design-based inference]
As pointed out in \citet{Bonnerytheseen}, under the design-based setting, where $D$, and $Y$ are considered non random, and the parametric space $\Gamma$ is separated (e.g. $\Gamma= \Thetaset\times\Xiset$, the selection is always ignorable. We could  find a definition of informative selection for fixed population models in the litterature: \citet[p.~12]{CasselSarndalWretman1977} defines a non informative design as a design such that design variable depend on study variables, but this definition is inconsistent with the fixed population framework described in the book design variables are not random, so necessarily independent. One could think that dependence here means functional dependence, in which case, it is very limiting as it corresponds to a case of non separation between $\Thetaset$ and $\Xiset$. Our conclusion is that it is not pertinent to consider the notion of informative selection in a pure fixed population model, as by nature of the model, all selections are non informative. For the notion of informativeness to be pertinent, one needs to apply it to a model where at least $Z\mid Y$ is random.
\end{remark}

\begin{property}[Sufficient condition for ignorability:Independence of $T$ and $Y$]
When $\forall (\theta,\xiparam)\in\Gamma$, $T$ and $Y$ are independent, and $G=(T,\Transfo[Y])$ is observed, the transformation is ignorable for the estimation of $\paramf(\theta)=\theta$.
\end{property}

\begin{proof}
 Assume $T$ and $Y$ are independent, let $t\in T(\Omega)$, $(\theta, \xiparam)\in\Gamma$. 
 Then $ P_{\theta,\xiparam}^{\Transfo[Y]\mid T=t}=\left( P_{\theta,\xiparam}^{Y\mid T=t}\right)^t$.
 As $T$ and $Y$ are independent, $ P_{\theta,\xiparam}^{Y\mid T=t}= P_{\theta,\xiparam}^{Y}= P_{\theta}^{Y}$. So $ P_{\theta,\xiparam}^{\Transfo[Y]\mid T=t}= P_{\theta}^{t(Y)}$.
\end{proof}

\begin{remark}[Rubin 1976]
In \cite{Rubin1976} the question of the amount of information brought by $T$ is ignored, the authors only 
consider $G=\Transfo[Y]$, the paper does not answer the question of what to do if we observe $G=(\Transfo[Y],T)$. It shows that the definition of ignorability and informative selection must also depend on $G$. In this paper, we propose a general definition for all possible $G$ (see \label{sec:oijweofij}.
\end{remark}

\begin{property}[Link to likelihood and density]
Assuming that $\Obs=(\Transfo[Y],T)$, that all densities and conditional densities are defined, 
and that $\left\{P^\Obs_{\theta,\xiparam}\right\}_{(\theta,\xiparam)\in\Gamma}$ is dominated by a $\sigma$-finite measure $\dominant$,  
then the following propositions are equivalent:
\begin{enumerate}
 \item The transformation $T$ is non informative for estimation of $\theta$ given $\Obs=(\Transfo[Y],T)$.
 \item $\forall \theta\in\Thetaset, \xiparam\in \Xiset$,  $\dominant$-a.s.$(y^\star,t)$, 
 $$\mathscr{L}((\theta,\xiparam)\mapsto\theta;G)(\theta;(y^\star,t))=\mathscr{L}(\theta\mapsto\theta,t(Y))(\theta;y^\star)$$
 \item $\forall \theta,\xiparam\in \Thetaset\times\Xiset$,  $\dominant$-a.s.$(y^\star,t)$, 
 $\density_{(\Transfo[Y],T);\theta,\xiparam}(y^\star,t)=\density_{(t(Y));\theta}(y^\star)\density_{T;\theta,\xiparam}(t)$ and $\density_{T;\theta,\xiparam}(t)$ does not depend on $\theta$.
 \end{enumerate}

 \end{property}

Note that \citet{Rubin1976} does not specify the nature of the model $( P_{\theta,\xiparam})_{\theta,\xiparam\in \Thetaset\times\Xiset}$, in particular, 
there is no comment relative to the model being dominated by a $\sigma$-finite measure, although this condition is implicitely implied by the use of densities.

%\subsection{Definitions in the Bayesian framework}
%\comments{This section needs to be rewritten.}
%
%\begin{definition}[Noninformative transformation]\label{def:ignorabletransformation}
%The transformation $\Rv$ is said {\em \bf ignorable} with respect to  $\paramf(P,?)$ given $\Obs=\obs$ if and only if the posterior distribution of $\paramf(P,Y)$ is unchanged after ignoring $\Rv$
%$$\priordist^{\paramf(P,Y)\mid \Obs=\obs}\text{-a.s.}(\param), \exists \rv_0, \priordist^{\Rv\mid \Obs=\obs}\text{-a.s.}(\rv),\ 
%\priordist^{\paramf(P,\Id_\Omega)\mid \Obs=\obs,\Rv=\rv}=\priordist^{\paramf(P,\Id_\Omega)\mid \Obs[\biclass(\rv,\Comp\Rv)]=\obs}$$
%\end{definition}

\end{document}